\documentclass[12pt]{amsart}
\usepackage[margin=1in,letterpaper,portrait]{geometry}
\usepackage{ytableau}
\ytableausetup{notabloids}
\ytableausetup{centertableaux}
\usepackage{latexsym,amsmath,amssymb,verbatim, tikz}
\usepackage{graphicx}
\usepackage[colorlinks]{hyperref}
\usepackage{mathrsfs}
\usepackage{amsfonts}
\usepackage{amsthm,youngtab}
\usepackage{graphicx}
\usepackage{caption}
\usepackage{enumerate,enumitem}
\usepackage{tikz-cd}
\usepackage{rotating}
\usepackage[all]{xy}
\usepackage[titletoc, title]{appendix}
\usepackage{amscd}
\usepackage{float}
\usepackage{hyperref} 
\hypersetup{
bookmarksnumbered,
pdfstartview={FitH},
breaklinks=true,
linkcolor=blue,
urlcolor=blue,
citecolor=blue,
bookmarksdepth=2
}	
\usepackage{caption}
\usepackage[utf8x]{inputenc}
\usepackage{dsfont}

\usepackage{stmaryrd} 
\usepackage{pgfplots}
\usepgfplotslibrary{external} 
\usepackage[colorinlistoftodos,bordercolor=orange,backgroundcolor=orange!20,linecolor=orange,textsize=scriptsize]{todonotes}
\makeatletter
\renewcommand{\todo}[2][]{\tikzexternaldisable\@todo[#1]{#2}\tikzexternalenable}
\makeatother

\theoremstyle{plain}
   \newtheorem{thm}{Theorem}[section]
  \newtheorem{theorem}{Theorem}
   \newtheorem{prop}[thm]{Proposition}
   \newtheorem{lemma}[thm]{Lemma}
   \newtheorem{corollary}[thm]{Corollary}
\theoremstyle{definition}
   \newtheorem{defn}[thm]{Definition}
   \newtheorem{example}[thm]{Example}
   \newtheorem{remark}[thm]{Remark}
\newtheorem{Question}{Question}

\numberwithin{equation}{section}

\newcommand{\CC}{{\mathbb {C}}}

\newcommand{\ZZ}{{\mathbb {Z}}}

\newcommand{\C}{\mathcal{C}}

\newcommand{\CQ}{\mathcal{C}_Q}
\newcommand{\YQ}{\mathcal{Y}_Q}
\newcommand{\Yxi}{\mathcal{Y}^{\leq \xi}}
\newcommand{\Cxi}{\mathcal{C}^{\leq \xi}}
\newcommand{\Axi}{\mathcal{A}^{\leq \xi}}
\newcommand{\Ixi}{I^{\leq \xi}}

\newcommand{\CN}{\mathbb{C}[\mathbf{N}]}
\newcommand{\yjh}{\hat{y_j}}

\newcommand{\s}{\mathcal{S}}

\newcommand{\iQ}{\mathbf{i}_Q}
\newcommand{\wiQ}{\widehat{\mathbf{i}}_Q}
\newcommand{\siQ}{\mathcal{S}^{\mathbf{i}_Q}}

\newcommand{\bCQ}{\overline{\mathcal{C}}_Q}
\newcommand{\bAQ}{\overline{\mathcal{A}}_Q}
\newcommand{\AQ}{\mathcal{A}_Q}
\newcommand{\barD}{\overline{D}}

\newcommand{\tchi}{{\widetilde{\chi}}}
\newcommand{\td}{{\widetilde{D}}}

\DeclareMathOperator*{\wt}{wt}

\newlength{\mysizetiny}
\newlength{\mysizesmall}
\newlength{\mysize}
\newlength{\mysizelarge}
\begin{document}

\title{Equivariant multiplicities via representations of quantum affine algebras}

\author{Elie CASBI and Jian-Rong LI}

\address{Elie CASBI, Max-Planck-Institut für Mathematik, Vivatsgasse 7, 53111 Bonn, Germany.} 
\email{casbi@mpim-bonn.mpg.de}

\address{Jian-Rong Li, Faculty of Mathematics, University of Vienna, Oskar-Morgenstern-Platz 1, 1090 Vienna, Austria.} 
\email{lijr07@gmail.com}

\date{}

\maketitle

\begin{abstract}
 For any simply-laced type simple Lie algebra $\mathfrak{g}$ and any height function $\xi$ adapted to an orientation $Q$ of the Dynkin diagram of $\mathfrak{g}$, Hernandez-Leclerc introduced a certain category $\Cxi$ of representations of the quantum affine algebra $U_q(\widehat{\mathfrak{g}})$, as well as a subcategory $\CQ$ of $\Cxi$ whose complexified Grothendieck ring is isomorphic to the coordinate ring $\CN$ of a maximal unipotent subgroup. In this paper, we define an algebraic morphism $\td_{\xi}$ on a torus $\Yxi$ containing the image of $K_0(\Cxi)$ under the truncated $q$-character morphism.  We prove that the restriction of $\td_{\xi}$ to $K_0(\CQ)$ coincides with the morphism $\barD$ recently introduced by Baumann-Kamnitzer-Knutson in their study of equivariant multiplicities of Mirkovi\'c-Vilonen cycles. This is achieved using the T-systems satisfied by the characters of Kirillov-Reshetikhin modules in $\CQ$, as well as certain results by Brundan-Kleshchev-McNamara  on the representation theory of quiver Hecke algebras. 
 This alternative description of $\barD$ allows us to prove a conjecture by the first author on the distinguished values of $\barD$ on the flag minors of $\CN$. 
 We also provide applications of our results from the perspective of Kang-Kashiwara-Kim-Oh's generalized Schur-Weyl duality. Finally, we use Kashiwara-Kim-Oh-Park's recent constructions to define a cluster algebra $\bAQ$ as a subquotient of $K_0(\Cxi)$ naturally containing $\CN$, and suggest the existence of an analogue of the Mirkovi\'c-Vilonen basis in $\bAQ$ on which the values of $\td_{\xi}$ may be interpreted as certain equivariant multiplicities.

\end{abstract}

\setcounter{tocdepth}{1}
\tableofcontents

\section{Introduction}

 Since their introduction by Drinfeld \cite{Drin87} and Jimbo \cite{Jim85}, the quantized enveloping algebras of Lie algebras and Kac-Moody algebras have been intensively studied and were at the heart of numerous important developments in the past decades. The quantum group $U_q(\mathfrak{g})$ associated to a finite-dimensional simple Lie algebra $\mathfrak{g}$ can be viewed as a deformation of the universal enveloping algebra of $\mathfrak{g}$. The construction of remarkable bases of the negative part $U_q(\mathfrak{n})$ arising from a triangular decomposition of $U_q(\mathfrak{g})$ has been one of them, initiated with the construction of the dual canonical basis by Lusztig  \cite{Lusztig} and the upper global basis by Kashiwara \cite{Kashicrystal}. Other bases with good properties were later considered, such as Lusztig's dual semi-canonical basis or the Mirkovi\'c-Vilonen basis arising from the geometric Satake correspondence \cite{MV}. The attempt towards a combinatorial description of the dual canonical basis has been one of the main motivations for the introduction of cluster algebras by Fomin and Zelevinsky \cite{FZ1}. It was proved by Berenstein-Fomin-Zelevinsky \cite{BFZ} that the coordinate ring $\CN$ of a maximal unipotent subgroup of the Lie group $\mathbf{G}$ of $\mathfrak{g}$ has a cluster algebra structure. This cluster algebra has infinitely many seeds in general, but there is a finite family 
 $$ \{ \s^{\mathbf{i}} = \left( (x_1^{\mathbf{i}}, \ldots , x_N^{\mathbf{i}}) , Q^{\mathbf{i}} \right) , \mathbf{i} \in \text{Red}(w_0) \} $$ of distinguished seeds called standard seeds, whose  cluster variables are given by explicit regular functions on $\mathbf{N}$ and whose exchange quiver $Q^{\mathbf{i}}$ is constructed purely combinatorially. They are indexed by the  set $\text{Red}(w_0)$ of all reduced expressions  of the longest element $w_0$ of the Weyl group of $\mathfrak{g}$. The cluster variables $x_1^{\mathbf{i}}, \ldots , x_N^{\mathbf{i}}$ are called the flag minors associated to $\mathbf{i}$.

 In their recent proof of Muthiah's conjecture \cite{Muthiah}, Baumann-Kamnitzer-Knutson \cite{BKK} introduced a remarkable algebra morphism 
 $$ \barD : \CN \longrightarrow  \mathbb{C}(\alpha_1, \ldots , \alpha_n) $$
 essentially via Fourier transforms of the Duistermaat-Heckmann measures (here $\alpha_1, \ldots , \alpha_n$ are formal variables corresponding to the simple roots of $\mathbf{G}$). They proved that the evaluation of $\barD$ on the elements of the Mirkovi\'c-Vilonen basis are related to certain geometric invariants of the corresponding Mirkovi\'c-Vilonen cycles called equivariant multiplicities, defined by Joseph \cite{Joseph}, Rossmann \cite{Rossmann} and later developed by Brion \cite{Brion}. Furthermore, the morphism $\barD$ turns out to be useful to compare good bases of $\CN$: in an appendix of the same work \cite{BKK}, Dranowski, Kamnitzer, and Morton-Ferguson use this morphism $\barD$ to prove that the MV basis and the dual semi-canonical basis are not the same.

The main purpose of the present paper is to extend $\barD$ to a larger algebra naturally containing $\CN$, defined as the complexified Grothendieck ring of a monoidal category $\Cxi$ of finite-dimensional representations of the quantum affine algebra $U_q(\widehat{\mathfrak{g}})$. 
This category was introduced by Hernandez-Leclerc \cite{HLJEMS} and they showed that the Grothendieck ring of the category has a cluster algebra structure. Recently Kashiwara-Kim-Oh-Park \cite{KKOP21} proved that $\Cxi$ provides a monoidal categorification of this cluster algebra in the sense of \cite{HL}.
Restricting our construction to $\CN$ allows us to investigate the behaviour of $\barD$ on the elements of the dual canonical basis using former results by Hernandez-Leclerc \cite{HL15}.  
  Our motivations are two-fold. 
  
   Firstly, the cluster algebra $\CN$ has another monoidal categorification using quiver Hecke algebra \cite{KKKO}. It was proved by the first author in \cite{Casbi3} that when $\mathfrak{g}$ is of type $A_n, n \geq 1$ or $D_4$, the morphism $\barD$ takes distinguished values on the flag minors of $\CN$, similar to its values on the classes of Kleshchev-Ram's strongly homogeneous modules over the quiver Hecke algebras associated to $\mathfrak{g}$ or on the elements of the MV basis corresponding to smooth MV cycles. Certain polynomial identities relating these values for flag minors belonging to the same standard seed were also exhibited and were proved (for all simply-laced types) to be preserved under cluster mutation from one standard seed to another ({{\cite[Theorem 5.6]{Casbi3}}}). But the cases where $\mathfrak{g}$ is of type $D_n, n \geq 5$ or $E_r, r=6,7,8$ were left open ({{\cite[Conjecture 5.5]{Casbi3}}}), and the meaning of these remarkable families of polynomial identities still remained mysterious.

 The second motivation for the present work comes from the fact that  Baumann-Kamnitzer-Knutson's morphism $\barD$ is known to admit natural interpretations in terms of various categorifications of $\CN$. For instance, the evaluation of $\barD$ on the elements of the dual semi-canonical basis can be naturally expressed in terms of the Euler characteristics of certain varieties of representations of the preprojective algebra associated to $\mathfrak{g}$. For an element of the dual canonical basis, viewed as the isomorphism class of a module $M$ over the quiver Hecke algebras associated to $\mathfrak{g}$, a similar formula can be written using the dimensions of the weight subspaces of $M$. However, the dual canonical basis admits another categorification, due to Hernandez-Leclerc \cite{HL15}, which involves certain finite-dimensional representations of the quantum affine algebra $U_q(\widehat{\mathfrak{g}})$. For instance, such representations were used in \cite{CDFL} to study the dual canonical basis of the Grassmannian cluster algebra $\mathbb{C}[{\rm Gr}(k,n)]$. 
   It is thus natural to ask whether the values of $\barD$ can be interpreted in a natural way using this other categorification. 
   
    The present paper provides answers to these questions. Although it is also related to the behaviour of $\barD$ with respect to the cluster structure of $\CN$, our approach involves different ideas from those of \cite{Casbi3}.  Furthermore, whereas the results of \cite{Casbi3} were proved using the representation theory of quiver Hecke algebras, and thus could only make sense on $\CN$, the framework we develop here allows to extend these results to larger cluster algebras and therefore opens new perspectives (see Section~\ref{finalsection} for example). It also yields natural proofs and interpretations of the polynomial identities mentioned above as well as several other conjectural observations made in \cite{Casbi3} (for instance {{\cite[Remark 6.4]{Casbi3}}}).

  \smallskip
  
   Hernandez-Leclerc's categorification of $\CN$ involves a family of monoidal categories $\CQ$ of finite-dimensional representations of $U_q(\widehat{\mathfrak{g}})$, indexed by the orientations $Q$ of the Dynkin diagram of $\mathfrak{g}$. The main result in \cite{HL15} consists in constructing for each choice of $Q$ a ring isomorphism from the Grothendieck ring $K_0(\CQ)$ to $\CN$ inducing a bijective correspondence between the classes of simple objects in $\CQ$ and the elements of the dual canonical basis of $\CN$ ({{\cite[Theorem 6.1]{HL15}}}). 
In \cite{KKOP21}, Kashiwara-Kim-Oh-Park defined a larger monoidal category $\Cxi$ containing $\CQ$ for each choice of height function $\xi$ adapted to $Q$. In the case where $\xi$ corresponds to a sink-source orientation, $\Cxi$ coincides with the category $\C^{-}$ introduced by Hernandez-Leclerc in \cite{HLJEMS}. 
There is an injective ring morphism 
$$ \widetilde{\chi}_q : K_0(\Cxi) \longrightarrow \Yxi := \mathbb{Z}[Y_{i,p}^{\pm 1} , (i,p) \in \Ixi] \enspace \text{where $\Ixi := \{(i,p) ,  i \in I , p \in \xi(i) + 2 \mathbb{Z}_{\leq 0} \}$} $$
 called truncated $q$-character morphism, which is a truncated version of Frenkel-Reshetikhin's $q$-character \cite{FR1}. An important family of simple objects in $\Cxi$ are the Kirillov-Reshetikhin modules $X_{i,p}^{(k)}, (i,p) \in \Ixi,  k \geq 1$, whose (truncated) $q$-characters are known to satisfy certain distinguished identities called $T$-systems (see \cite{H06}). 
    It is shown ({{\cite[Theorem 5.1]{HLJEMS}}}) that the Grothendieck ring of $\Cxi$ has a cluster algebra structure, with an initial seed given by an explicit (infinite) quiver $Q_{\xi}$ ({{\cite[Section 2.1.2]{HLJEMS}}}) with set of vertices $\Ixi$, together with a cluster consisting of the isomorphism classes of the Kirillov-Reshetikhin modules of the form $X_{i,p} := X_{i,p}^{(1 + (\xi(i)-p)/2)} , (i,p) \in \Ixi$. The T-system relations between the characters of the $X_{i,p}^{(k)}$ correspond to the exchange relations of certain sequences of mutations for the cluster structure of $K_0(\Cxi)$. Moreover, denoting by $\iQ$ a reduced expression of $w_0$ adapted to $Q$, the exchange quiver $Q^{\iQ}$ of the standard seed $\siQ$ of $\CN$ can be viewed as a (finite) subquiver of $Q_{\xi}$, and the cluster variables (flag minors) $x_1^{\iQ}, \ldots , x_N^{\iQ}$ of $\siQ$ are identified with the classes of the modules $X_{i,p}, (i,p) \in I_Q$ via the natural embedding $\CN \simeq K_0(\CQ) \hookrightarrow K_0(\Cxi)$ (where $I_Q$ is a finite subset of $\Ixi$).

  
  \smallskip
  
 In this paper we introduce an algebra morphism $\td_{\xi}$ from the complexified torus $ \mathbb{C} \otimes \Yxi$ to the field $\mathbb{C}(\alpha_1, \ldots, \alpha_n)$. Its definition involves  the coefficients of the inverse of the quantum matrix of $\mathfrak{g}$, which are a family of integers $\tilde{C}_{i,j}(m) , i,j \in I, m \in \mathbb{Z}$ appearing in the theory of $q$-characters \cite{FR1} initiated by Frenkel and Reshetikhin, $q,t$-characters initiated by Nakajima \cite{NakaAnn.Math}, and then further developed in \cite{HeAdv,HL15}. We also refer to \cite{FHOO} for recent advances in this area.
 The precise definition of $\td_{\xi}$ is the following:
\begin{align*}
  \forall (i,p) \in \Ixi, \enspace \td_{\xi}(Y_{i,p}) := \prod_{(j,s) \in \Ixi} \left( \epsilon_{j,s} \tau_Q^{(\xi(j)-s)/2}(\gamma_j) \right)^{\tilde{C}_{i,j}(s-p-1) - \tilde{C}_{i,j}(s-p+1)}
\end{align*}
where for each $j \in I$, $\gamma_j$ is the sum of the simple roots $\alpha_i$ such that there exists a path from $i$ to $j$ in $Q$, $\tau_Q$ is the Coxeter transformation associated to $Q$ and $\epsilon_{j,s} \in \{-1,1 \}$ is the unique sign such that $\tau_Q^{(\xi(j)-s)/2}(\gamma_j) \in \epsilon_{j,s} \Phi_{+}$ for every $(j,s) \in \Ixi$. Note that this product is always finite, because $\tilde{C}_{i,j}(m) := 0$ if $m \leq 0$. In Section~\ref{ABCforDtilde}, we investigate the images under $\td_{\xi}$ of the truncated $q$-characters of the Kirillov-Reshetikhin modules $X_{i,p}, (i,p) \in \Ixi$ categorifying the cluster variables of Hernandez-Leclerc's initial seed in $K_0(\Cxi)$. We prove that the rational fractions $\td_{\xi} \left( \tchi_q(X_{i,p}) \right)$ satisfy remarkable properties analogous to those exhibited in \cite{Casbi3} for the values of $\bar{D}$ on the flag minors of $\CN$. 

We also consider the restriction $\td_Q$ of $\td_{\xi}$ to the torus $\YQ :=  \mathbb{Z}[Y_{i,p}^{\pm 1}, (i,p) \in I_Q]$ image of $K_0(\CQ)$ under the truncated $q$-character morphism $\widetilde{\chi}_q$. 
Our first main result is the following:
 
\begin{theorem}[cf. {Theorem \ref{mainthm1}}] \label{thm:mainthm1 in introduction}
For every simply-laced type Lie algebra $\mathfrak{g}$ and for every orientation $Q$ of the Dynkin diagram of $\mathfrak{g}$, the following diagram commutes:
\begin{align*}
\xymatrix{ \CN \ar[rd]_{\barD} \ar[r]^{\simeq} & \mathbb{C} \otimes K_0(\CQ) \ar@{^{(}->}[r]^{\tchi_q}  & \mathbb{C} \otimes \YQ \ar[ld]^{\tilde{D}_Q} \\
{} &  \mathbb{C}(\alpha_i , i \in I)  & {} } 
\end{align*}
\end{theorem}

 A significant part of this paper (Sections~\ref{Dbarcuspidals} and~\ref{DtildeKRmonotonic} as well as the beginning of Section~\ref{sectionproofs}) will be devoted to proving Theorem~\ref{thm:mainthm1 in introduction} in the case of a particular well-chosen orientation $Q_0$ for each simply-laced type (see Section~\ref{Dbarcuspidals}). It is achieved by proving that $\barD$ and $\td_{Q_0}$ agree on the dual root vectors associated to $\mathbf{i}_{Q_0}$, which are known to generate $\CN$ as an algebra. The dual root vectors are categorified on the one hand by the so-called cuspidal representations over quiver Hecke algebras (see \cite{BKM,KR,McNamarafinite}) and by the fundamental representations in $\CQ$ on the other hand. In Section~\ref{Dbarcuspidals}, we use the representation theory of quiver Hecke algebras to provide formulas for the evaluation of $\barD$ on the dual root vectors of $\CN$. Our proof crucially relies on certain results by Brundan-Kleshchev-McNamara \cite{BKM} on cuspidal representations as well as Kleshchev-Ram's construction \cite{KRhom} of (strongly) homogeneous modules. In Section~\ref{DtildeKRmonotonic}, we prove several formulas for the evaluation of $\td_{Q_0}$ on the classes of \textit{all} Kirillov-Reshetikhin modules in $\CQ$ using the T-system relations satisfied by the (truncated) $q$-characters of these modules. As the fundamental representations are particular cases of Kirillov-Reshetikhin modules, we can conclude by comparing with the results obtained in Section~\ref{Dbarcuspidals}. 
 
 The values of $\barD$ (resp. $\td_{Q_0}$) on the dual root vectors of $\CN$ (resp. the classes of Kirillov-Reshetikhin modules in $\C_{Q_0}$) are obtained in Section~\ref{Dbarcuspidals} (resp. Section~\ref{DtildeKRmonotonic}) by considering each simply-laced type. The case of type $A_n$ is in fact contained as a subcase of the type $D_n$ but for the reader's convenience we chose to state the formulas in different subsections for each of these types. We deal with the types $E_r, r=6,7,8$ separately using a computer software. For certain dual root vectors, the results are extremely complicated, which suggests there exists probably no uniform formula for the image of $\td_Q$ on the classes of Kirillov-Reshetikhin modules (or even simply on the dual root vectors) that would hold for any simply-laced type and for an arbitrary orientation $Q$.

\smallskip
  
 Combining Theorem~\ref{thm:mainthm1 in introduction} with {{\cite[Theorem 5.6]{Casbi3}}} allows us to prove the second main result of this paper, which was stated as a conjecture in \cite{Casbi3} ({{\cite[Conjecture 5.5]{Casbi3}}}). 
 
\begin{theorem}[cf. {Theorem \ref{mainthm2}}] \label{thm:mainthm2 in introduction}
Let $\mathfrak{g}$ be a simple Lie algebra of simply-laced type. Then for any reduced expression $\mathbf{i} = (i_1, \ldots, i_N)$ of $w_0$, the flag minors $x_1^{\mathbf{i}} , \ldots , x_N^{\mathbf{i}}$ satisfy $\barD(x_j^{\mathbf{i}}) = 1/P_j^{\mathbf{i}}$ where $P_j^{\mathbf{i}}$ is a product of positive roots. Furthermore, the polynomials $P_1^{\mathbf{i}} , \ldots , P_N^{\mathbf{i}}$ satisfy the identities
\begin{align*}
\forall 1 \leq j \leq N, \enspace  P_j^{\mathbf{i}} P_{j_{-}}^{\mathbf{i}}  = \beta_j \prod_{\substack{l<j<l_{+} \\ i_l \sim i_j}}  P_l^{\mathbf{i}}.
\end{align*}   
where $\beta_j := s_{i_1} \cdots s_{i_{j-1}} (\alpha_{i_j})$ for each $j \in \{1, \ldots , N \}$. 
\end{theorem}

We refer to Sections~\ref{remindCN} and~\ref{remindAR} for precise definitions of the notations involved in this identity. 
We first prove the statement in the case $\mathbf{i} = \mathbf{i}_{Q_0}$ using Theorem~\ref{thm:mainthm1 in introduction}, which provides an efficient way of computing the images under $\barD$ of the flag minors $x_1^{\mathbf{i}_{Q_0}} , \ldots , x_N^{\mathbf{i}_{Q_0}}$ because the Kirillov-Reshetikhin modules $X_{i,p}$ have truncated $q$-characters reduced to a single term. Then  {{\cite[Theorem 5.6]{Casbi3}}} guarantees that the result holds for arbitrary reduced expressions of $w_0$. 
  In the case of reduced expressions of $w_0$ adapted to orientations of the Dynkin diagram of $\mathfrak{g}$, these polynomial identities are now naturally understood as consequences of the well-known recursive relations between the coefficients $\tilde{C}_{i,j}(m)$. Section~\ref{ABCforDtilde} contains further explanations about this, as well as analogous interpretations of various other observations from \cite{Casbi3}, such as {{\cite[Remark 6.4]{Casbi3}}}. 

 \smallskip

 We have been informed that it could be also possible to obtain a geometric proof of Theorem~\ref{thm:mainthm1 in introduction} relying on the geometric Satake correspondence and the results from \cite{BKK}. The idea is to prove that the MV cycles associated to the flag minors of the standard seed $\s^{\iQ}$ satisfy certain smoothness properties, which allow one to compute their equivariant multiplicities, and hence the values of $\barD$ using  {{\cite[Corollary 10.6]{BKK}}}. One could then conclude by combining this with the results of Section~\ref{sectionDtilde} of the present paper.
 
  Our approach has the advantage to provide closed formulas for the evaluation of $\barD$ on a large family of cluster variables in $\CN$ (namely, all the classes of the Kirillov-Reshetikhin modules in $\mathcal{C}_{Q_0}$), several of which correspond to MV cycles that fail to satisfy any smoothness property (this can be seen for example by the fact the numerators may not be equal to $1$ in our formulas). Furthermore, the techniques used in the present paper can be extended in a direct way beyond $\CN$ to obtain formulas for the evaluation of $\td_{\xi}$ on classes of Kirillov-Reshetikhin modules in $K_0(\Cxi)$, which are not described by the geometric Satake correspondence. 

\smallskip
 
 In Section~\ref{sectionappliSW} we present some applications of our results involving Kang-Kashiwara-Kim-Oh's generalized quantum affine Schur-Weyl duality \cite{KKKOSelecta}. 
An element of the dual canonical basis of $\CN$ can be viewed either as the class of a (simple) module $L$ in $\CQ$, or as the class of the corresponding object $\mathcal{F}_Q(L)$ in $R-mod$, where $\mathcal{F}_Q$ denotes Kang-Kashiwara-Kim-Oh's generalized quantum affine Schur-Weyl duality functor. Then Theorem \ref{thm:mainthm1 in introduction} yields the following identity
\begin{equation} \label{eqDtildeintro}
 \sum_{\mathfrak{m}' \preccurlyeq \mathfrak{m}} \dim(L_{\mathfrak{m}'}) \td_Q(\mathfrak{m}') =  \sum_{\mathbf{j}=(j_1, \ldots j_d)} \dim \left( (\mathcal{F}_Q(L))_{\mathbf{j}} \right) \barD_{\mathbf{j}} 
  \end{equation}
  where $\mathfrak{m}'$ are the Laurent monomials in the variables $Y_{i,s}$ appearing in the truncated $q$-character of $L$, and the $\barD_{\mathbf{j}}$ are certain explicit rational fractions for each weight $\mathbf{j}$ of $\mathcal{F}_Q(L)$ (see Section~\ref{sectionappliSW}). 
 In other words, although the objects of the categories $\CQ$ and $R-mod$ are a priori of different natures, the equality~\eqref{eqDtildeintro} provides an unexpected explicit relationship between the respective weight-subspace structures of a representation $L$ of $\CQ$ and the corresponding  module $\mathcal{F}_Q(L)$ in $R-mod$. We provide a concrete illustration of this fact by proving a formula relating the dimensions of $\mathcal{F}_{Q_0}(L)$ and the truncated part of $L$ when $\mathfrak{g}$ is of type $A_n$ and $Q=Q_0$ (Theorem~\ref{thmapplicationtoKKKO-SW}). 
 
 \smallskip
 
  In the final section of this work (Section~\ref{finalsection}), we turn back to the geometric motivations at the origin of the construction of $\barD$ by Baumann-Kamnitzer-Knutson \cite{BKK}. As our results show that $\td_{\xi}$ is an extension of $\barD$ to $\mathbb{C} \otimes K_0(\Cxi)$, it is natural to ask whether certain values of $\td_{\xi}$ may be possibly  related to certain equivariant multiplicities of some closed algebraic varieties. However, it turns out that unlike $\barD$, $\td_{\xi}$ takes trivial values on certain cluster variables of $K_0(\Cxi)$, which seems difficult to understand geometrically. To circumvent this issue, we show that the values of $\td_{\xi}$ on the cluster variables $[X_{i,p}]$ of Hernandez-Leclerc's initial seed in $K_0(\Cxi)$ satisfy certain periodicity properties (Corollary~\ref{corperiodicity}). This is derived from the periodicity of the coefficients $\tilde{C}_{i,j}(m)$ established by Hernandez-Leclerc ({{\cite[Corollary 2.3]{HL15}}}). Therefore one looses essentially no information by restricting $\td_{\xi}$ to a smaller cluster algebra $\AQ$ of finite cluster rank still containing $\CN$ as a subalgebra. In some sense, one can view $\AQ$ as a period of $\td_{\xi}$ and $\CN$ as a half-period of $\td_{\xi}$. Then we prove (Corollary~\ref{corfrozen}) that $\td_{\xi}$ takes trivial values on the frozen variables of $\AQ$. Therefore it factors through a morphism $\barD_Q$ defined on the quotient algebra $\bAQ$. We propose this algebra $\bAQ$ as the most appropriate domain to study $\td_{\xi}$. We ask for the existence of a basis in $\bAQ$ containing the Mirkovi\'c-Vilonen basis of $\CN$, whose elements may be indexed by a family of closed algebraic varieties, such that the values of $\barD_Q$ on the elements of this basis could be interpreted as certain equivariant multiplicities of the corresponding varieties with respect to the action of some torus. In a different direction, the recent results by Kashiwara-Kim-Oh-Park \cite{KKOP21} imply that $\AQ$ admits a monoidal categorification (in the sense of Hernandez-Leclerc \cite{HL}) by a subcategory of $\Cxi$. Therefore, it would be interesting to investigate whether the quotient algebra $\bAQ$ could also be studied via monoidal categorification using the categorical specialization techniques developped by Kang-Kashiwara-Kim \cite{KKKInvent18} and more recently Kashiwara-Kim-Oh-Park \cite{KKOPlocalization}. 

 \smallskip
 Our construction also suggests possible connections with the quantum cluster algebra structures of quantized coordinate rings and more generally of certain quantum Grothendieck rings of Hernandez-Leclerc's categories. Indeed, the expressions $\tilde{C}_{i,j}(s-p-1) - \tilde{C}_{i,j}(s-p+1)$ involved in the definition of $\td_{\xi}$ coincide (up to sign) with the entries of the $t$-commutation matrix  describing Hernandez-Leclerc's quantum torus \cite[Equation (8)]{HL15}. This is known from the works of Bittmann \cite{Lea21} to correspond to the quantum torus of a quantum cluster algebra $K_t(\Cxi)$, which is a non-commutative deformation of $K_0(\Cxi)$ naturally containing the quantized coordinate ring $\mathcal{A}_q(\mathfrak{n})$ (dual of $U_q(\mathfrak{n})$) via an algebra isomorphism $\mathcal{A}_q(\mathfrak{n}) \simeq K_t(\CQ)$ identifying the indeterminates $t$ and $q$.

\smallskip
 
The paper is organized as follows. In Section~\ref{remindHL15} we gather all the necessary reminders about Hernandez and Leclerc's categorification of cluster algebras  and its applications to the study of the coordinate ring $\CN$. In Section~\ref{remindCtilde}, we provide some reminders about the coefficients of the inverses of quantum Cartan matrices and prove a couple of elementary properties which will be useful in the sequel of the paper. In Section~\ref{sectionremindKLR}, we recall the main results from the representation theory of quiver Hecke algebras, in particular certain results from~\cite{BKM}. Section~\ref{remindDbar} is devoted to some reminders about  Baumann-Kamnitzer-Knutson's morphism $\barD$ as well as the main results from~\cite{Casbi3}. In Section~\ref{sectionDtilde}, we introduce the main objects of the present paper, namely the morphisms $\td_{\xi}$ and $\td_Q$, investigate several of their properties and state our first main result Theorem \ref{mainthm1}. In Section~\ref{Dbarcuspidals}, we use the representation theory of quiver Hecke algebras to compute the values of $\barD$ on the dual root vectors of $\CN$ associated to a well-chosen orientation $Q_0$ of the Dynkin diagram of $\mathfrak{g}$. In Section~\ref{DtildeKRmonotonic}, we provide formulas in types $A_n$ and $D_n$ for the values of $\td_{Q_0}$ on the classes of all Kirillov-Reshetikhin modules in $\mathcal{C}_{Q_0}$. Section~\ref{sectionproofs} is devoted to the proofs of the two main results of this paper, Theorems~\ref{mainthm1} and~\ref{mainthm2}. We begin by proving Theorem~\ref{mainthm1} in the case  $Q=Q_0$, which allows us to prove Theorem~\ref{mainthm2}, which can then be used to prove Theorem~\ref{mainthm1} for arbitrary orientations. In Section~\ref{sectionappliSW}, we provide a representation-theoretic interpretation of Theorem~\ref{mainthm1} from the perspective of Kang-Kashiwara-Kim-Oh's generalized quantum affine Schur-Weyl duality, with an application when $\mathfrak{g}$ is of type $A_n$ (Theorem~\ref{thmapplicationtoKKKO-SW}). Finally in Section~\ref{finalsection} we define a cluster algebra $\bAQ$ naturally containing $\CN$ and discuss possible geometric interpretations of the values taken by $\td_{\xi}$ on $\bAQ$.

\subsection*{Acknowledgements}
 We would like to thank Jan Schr\"oer for his encouragements and for many enjoyable discussions. We are grateful to Alexander Kleshchev for his help and his recommendation of using the results in \cite{BKM}, which turned out to play a key role in Section~\ref{Dbarcuspidals}. We also thank Joel Kamnitzer and Pierre Baumann for several precious explanations about their inspiring work \cite{BKK} and Anne Dranowski for many fruitful conversations. Finally we thank Myungho Kim for answering our questions about the recent works \cite{KKOPlocalization,KKOP21} and the anonymous referee for many valuable suggestions and remarks. 

E.C. thanks the Max-Planck-Institut f\"ur Mathematik for the opportunity of working as a postdoctoral fellow under Jan Schr\"oer's mentorship. 
J.L. is supported by the Austrian Science Fund (FWF): M 2633-N32 Meitner Program and P 34602 individual project.

\section{Hernandez-Leclerc's category $\Cxi$}
 \label{remindHL15}

In this section, we recall Hernandez-Leclerc's categorifications of certain cluster algebras via the categories $\Cxi$ introduced in \cite{HLJEMS} (denoted $\C^{-}$ in \cite{HLJEMS}), for each height function $\xi$ on the vertices of the Dynkin diagram of $\mathfrak{g}$. We also recall Hernandez-Leclerc's former construction \cite{HL15} of  categorifications of coordinate rings via subcategories $\CQ$ of $\Cxi$, where $Q$ is an orientation of the Dynkin graph of $\mathfrak{g}$ adapted to $\xi$. 

 \smallskip
 
  \subsection{Coordinate rings and their cluster structures}
  \label{remindCN}
  
  Let $\mathfrak{g}$ be a simple complex Lie algebra of simply-laced type, let $I$ be the set of vertices of the Dynkin diagram of $\mathfrak{g}$, and let $n=|I|$. We denote by $C=(c_{i,j})$ the Cartan matrix of $\mathfrak{g}$ and for any $i,j \in I$ we will write $i \sim j$ for $c_{i,j} = -1$. 
 Let us fix a nilpotent subalgebra $\mathfrak{n}$ arising from a triangular decomposition of $\mathfrak{g}$ and let $\mathbf{N}$ denote the corresponding Lie group. We consider the ring $\CN$ of regular functions on $\mathbf{N}$, which we will refer to as the \textit{coordinate ring} in what follows. We also let $W$ denote the Weyl group of $\mathfrak{g}$ and $w_0$ denote the longest element of $W$. We let $\alpha_1, \ldots , \alpha_n$ (resp. $\omega_1, \ldots , \omega_n$) denote the simple roots (resp. the fundamental weights) of $\mathfrak{g}$. Let $\Gamma_{+} := \bigoplus_i \mathbb{N} \alpha_i$ and let $\Phi_{+} \subset \Gamma_{+}$ denote the set of positive roots of $\mathfrak{g}$. 
 
  \smallskip
   
  The coordinate ring $\CN$ contains a distinguished family of elements $D(u \lambda , v \lambda)$ called \textit{unipotent minors} parametrized by triples $(\lambda,u,v) \in P_{+} \times W \times W$ where $P_{+}$ stands for the set of dominant weights of $\mathfrak{g}$. These unipotent minors always belong to the dual canonical basis of $\CN$ when they are not zero (see for instance {{\cite[Lemma 9.1.1]{KKKO}}}). Two special subsets of unipotent minors will play a central role throughout this paper, both of them depending on a choice of reduced expression of $w_0$. Let $N := \sharp \Phi_{+} = l(w_0)$ and let $\mathbf{i}=(i_1, \ldots , i_N)$ be a reduced expression of $w_0$. On the one hand, we will consider the unipotent minors 
  $$ x_k^{\mathbf{i}} := D(s_{i_1} \cdots s_{i_k} \omega_{i_k}, \omega_{i_k}) \enspace , \enspace  1 \leq k \leq N $$
   which are called the \textit{flag minors} associated to $\mathbf{i}$. On the other hand, the unipotent minors 
   $$ {r_k^{\ast}}^{\mathbf{i}} :=  D(s_{i_1} \cdots s_{i_k} \omega_{i_k}, s_{i_1} \cdots s_{i_{k-1}} \omega_{i_k}) \enspace , \enspace 1 \leq k \leq N $$ 
   are called the \textit{dual root vectors} associated to $\mathbf{i}$. Note that the dual root vectors also belong to the dual PBW basis corresponding to $\mathbf{i}$. 
  
 One of the properties of $\CN$ we will be mostly interested in is its \textit{cluster algebra structure} in the sense of Fomin-Zelevinsky \cite{FZ1}: a cluster algebra is defined as the subalgebra of a field of functions $\mathbb{Q}(x_1, \ldots , x_N)$ generated by a family of distinguished elements called \textit{cluster variables}. These  are obtained by performing a recursive procedure starting from the initial data (called a \textit{seed}) of a $N$-tuple $x_1, \ldots , x_N$ of variables (called a cluster) as well as a quiver with $N$ vertices and without any loop or $2$-cycles (called an exchange quiver). For every $1 \leq k \leq N$, one can define a new generator $x'_k$ given by the \textit{exchange relation}
 \begin{equation} \label{exchangerelation}
  x'_k = \frac{1}{x_k} \left( \prod_{ \text{$j \rightarrow k$ in $Q$}} x_j + \prod_{ \text{$j \leftarrow k$ in $Q$}} x_j  \right) 
   \end{equation}
  as well as a new quiver $Q'_k$, both uniquely determined by $x_1, \ldots , x_N$ and $Q$. This yields a new seed given by $ \s'_k := ((x_1, \ldots , x_{k-1} , x'_k , x_{k+1} , \ldots , N) , Q'_k)$. This procedure is called the \textit{mutation in the direction $k$} of the seed $ \s := ((x_1, \ldots , x_N),Q)$. It has the important property of being involutive, i.e. performing the mutation in the same direction $k$ to the seed $\s'_k$ recovers the seed $\s$. Iterating this for all possible sequences of directions of mutations, we get a (finite or infinite) set of new generators called cluster variables, each of them appearing in several clusters. The rank of a cluster algebra is the cardinality of each of its clusters.
 
  The general theory of cluster algebras developed in particular in \cite{FZ4} has found a large range of applications to various areas of mathematics such as representation theory, Poisson geometry, representations of quivers etc... As far as coordinate rings are concerned, the main result that will be relevant for us is the following.
  
\begin{thm}[Berenstein-Fomin-Zelevinsky \cite{BFZ}, Geiss-Leclerc-Schr\"oer \cite{GLS}]

\begin{enumerate}
  \item The coordinate ring $\CN$ has a cluster algebra structure, of rank equal to the number of positive roots of $\mathfrak{g}$.
  \item For each reduced expression $\mathbf{i}$ of $w_0$, the flag minors $D(s_{i_1} \cdots s_{i_k} \omega_{i_k}, \omega_{i_k}), 1 \leq k \leq N$ form a cluster in $\CN$.
\end{enumerate}
\end{thm}
  
   There is in addition a purely combinatorial way of defining a quiver $Q^{\mathbf{i}}$ with $N$ vertices for each reduced expression $\mathbf{i}$, which yields a seed in $\CN$:
   $$ \s^{\mathbf{i}} = \left( (x_1^{\mathbf{i}}, \ldots , x_N^{\mathbf{i}}) , Q^{\mathbf{i}} \right) . $$
   The seeds $\s^{\mathbf{i}}$ are called the \textit{standard seeds} of $\CN$. Note that different reduced expressions may yield the same seed: this is the case for instance for reduced expressions in the same commutation class. Moreover, cluster mutations from one standard seed to another correspond to performing braid relations on the corresponding reduced expressions.

     \begin{remark}
     \begin{enumerate}
      \item  The cluster structure of $\CN$ can have infinitely many seeds in general (in fact it is always the case, unless $\mathfrak{g}$ is of type $A_n, n \leq 4$ see \cite{GLSAENS})
       \item  The exchange relations associated to the cluster mutations from one standard seed to another are special cases of the \textit{determinantal identities} which play an important role in the work of Geiss-Leclerc-Schr\"oer  (see {{\cite[Proposition 5.4]{GLS}}}) as well as in the work of Fomin-Zelevinsky \cite{FZJAMS}.  
      \end{enumerate}
      \end{remark}
      
 
  \subsection{Auslander-Reiten theory}
\label{remindAR}

 We fix an orientation $Q$ of the Dynkin diagram of $\mathfrak{g}$. Let $\langle \cdot , \cdot \rangle_Q$  denote the Euler-Ringel form of the quiver $Q$, i.e. the unique bilinear form on the free abelian group $\Gamma := \bigoplus_i  \mathbb{Z} \alpha_i$ given on simple roots by $\langle \alpha_i, \alpha_j \rangle_Q := \delta_{i,j} - \sharp \{ \text{$i \rightarrow j$ in $Q$} \}$. Finally we denote by $( \cdot , \cdot)$ the Cartan pairing associated to $\mathfrak{g}$, i.e. the symmetric bilinear form on $\Gamma_{+}$ defined by $(\alpha_i,\alpha_j)=c_{i,j}$. Recall that one has $ (\beta , \gamma) = \langle \beta , \gamma \rangle_Q + \langle \gamma , \beta \rangle_Q$ for any $\beta, \gamma \in \Gamma_{+}$.

  We fix a height function $\xi : I \rightarrow \mathbb{Z}$ adapted to $Q$ i.e. an integer-valued function satisfying
 $$ \xi(j)  = \xi(i) - 1 \quad \text{if there is an arrow $i \rightarrow j$ in $Q$}. $$
For $i \in I$, denote by $s_i(Q)$ the quiver obtained from $Q$ by changing the orientation of every arrow with source $i$ or target $i$. A sequence ${\bf i} = (i_1, \ldots, i_k)$ of elements of $I$ is called {\sl adapted} to $Q$ if $i_1$ is a source of $Q$, $i_2$ is a source of $s_{i_1}(Q)$, $\ldots$, $i_k$ is a source of $s_{i_{k-1}}\cdots s_{i_1}(Q)$. There is a unique Coxeter element in $W$, denoted by $\tau_Q$, having reduced expressions adapted to $Q$. It satisfies $\tau_Q^{h}=1$ where $h$ denotes the dual Coxeter number defined by $h := 2N/n$, where $N$ is the number of positive roots in the root system of $\mathfrak{g}$.

  Most importantly, we also fix a reduced expression $\iQ=(i_1, \ldots , i_N)$ of $w_0$ adapted to $Q$. Recall that such a reduced word always exists and is unique up to commutation. For each $i \in I$, we denote by $n_Q(i)$ the number of occurrences of the letter $i$ in the reduced word $\iQ$. 
For every $i \in I$, let $B(i)$ denote the set of indices $j \in I$ such that there exists a path from $j$ to $i$ in $Q$, and let $\gamma_i := \sum_{j \in B(i)} \alpha_j$. Then $\gamma_i \in \Phi_{+}$ for every $i \in I$, and moreover one has 
  $$ \Phi_{+} = \{ \tau_Q^{r-1}(\gamma_i) , i \in I , 1 \leq r \leq n_Q(i) \} . $$
Following \cite{KKOP21}, we also define an infinite sequence $\wiQ =(i_1, i_2, \ldots )$ of elements of $I$ as follows. 
For each $i \in I$, we let $i^{*}$ denote the unique element of $I$ such that $w_0(\alpha_i) = - \alpha_{i^{*}}$. The map $i \mapsto i^{*}$ is an involution. Then for each $1 \leq k \leq N$ and each $m \geq 0$, we set 
$$ i_{k+Nm} := \begin{cases}
 i_k &\text{if $m$ is even,} \\
 i_k^{*} & \text{if $m$ is odd.}
\end{cases}
$$
It is proved (see {{\cite[Proposition 6.11]{KKOP21}}}) that for any $t \geq 1$, the finite sequence 
\[
(i_t, i_{t+1}, \ldots, i_{t+N-1}) 
\]
is a reduced expression of $w_0$ adapted to the orientation $s_{i_{t-1}} \cdots s_{i_1}(Q)$. 

 We will use the following notation from {{\cite[Equation (4.2)]{KKOP21}}}. For each $t \geq 1$, we set 
 $$ t_{+} := \min \left( \{ t'>t , i_{t'}=i_t \} \cup \{ + \infty \} \right) \enspace \text{and} \enspace t_{-} := \max \left( \{ t'<t , i_{t'}=i_t \} \cup \{ 0 \} \right) . $$
  For $t_+$, the set $\{ t'>t , i_{t'}=i_t \}$ is never empty. For $t_-$, the set $\{ t'<t , i_{t'}=i_t \}$ can be empty so we use the convention $\max \emptyset := 0$.  
 
  We set 
  $$ \Ixi := \{ (i,p) \mid i \in I, p \in \xi(i) + 2 \mathbb{Z}_{\leq 0} \} . $$
   There is a bijection $\varphi : \Ixi \longrightarrow \mathbb{Z}_{\geq 1}$ defined by 
  \begin{equation} \label{defvarphi}
   \varphi(i,p) := 
  \begin{cases}
   \min \{ t \geq 1, i_t=i \} &\text{if $p=\xi(i)$,} \\
   \left( \varphi(i,p+2) \right)_{+} &\text{otherwise.}
  \end{cases}
  \end{equation}
 Equivalently, $\varphi(i,p)$ is the position of the $m$th occurrence of the letter $i$ in $\wiQ$, where $m := (\xi(i)-p+2)/2$.
The inverse of $\varphi$ is given by 
$$ \varphi^{-1}(k) = (i_k, \xi(i_k)-2N_Q(k)+2) \quad \text{where $N_Q(k) := \sharp \{ k' \leq k , i_{k'} = i_k \}$} . $$
We also set 
 \begin{equation} \label{defIQ}
 I_Q :=  \varphi^{-1} \left( \{1, \ldots , N \} \right) = \{ (j,s) \in \Ixi , \enspace \xi(j) \geq s \geq \xi(j)-2n_{Q}(j)+2 \} . 
 \end{equation}
Then following \cite{HL15} one can define a sequence of positive roots $(\beta_k)_{k \geq 1} \in \Phi_{+}^{\mathbb{Z}_{\geq 1}}$ by setting
  \begin{equation} \label{defbetaxi}
  \beta_{\varphi(j, \xi(j))} := \gamma_j \enspace \text{and} \enspace  \beta_{\varphi(j,s-2)} = 
  \begin{cases}
  \tau_{Q} \left( \beta_{\varphi(j,s)} \right) &\text{if $\tau_{Q} \left( \beta_{\varphi(j,s)} \right) \in \Phi_{+}$,} \\
  -  \tau_{Q} \left( \beta_{\varphi(j,s)} \right) &\text{if $\tau_{Q} \left( \beta_{\varphi(j,s)} \right) \in - \Phi_{+}$.}
  \end{cases}
    \end{equation}
  for every $(j,s) \in \Ixi$. We also denote by $\epsilon_{j,s} \in \{-1,1\}$ the unique sign such that 
\[  
\tau^{(\xi(j)-s)/2}(\gamma_j) = \epsilon_{j,s} \beta_{\varphi(j,s)}.
\] 
  We have the following result. 
  
   \begin{prop}[{{\cite[Corollary 2.40]{FOComm.Math.Phys.}}}]  \label{propositionFujitaOh}
   For every $(j,s) \in \Ixi$ one has $\beta_{\varphi(j,s)} = \beta_{\varphi(j^{*},s+h)}$  and  $\epsilon_{j,s} = - \epsilon_{j^{*},s+h}$.
    \end{prop}
    
  In other words one has $\tau^{(\xi(j)-s)/2}(\gamma_j) = - \tau^{(\xi(j^{*})-s-h)/2}(\gamma_{j^{*}})$ for every $(j,s) \in \Ixi$. Consequently, given $(j,s)$ and $(j',s')$ in $\Ixi$, one has $\beta_{\varphi(j,s)} = \beta_{\varphi(j',s')}$ if and only if $j'=j$ and $s' \in s + 2h \mathbb{Z}$ (and in this case $\epsilon_{j,s} = \epsilon_{j',s'}$) or $j'=j^{*}$ and $s' \in s + h + 2h \mathbb{Z}$ (and in this case $\epsilon_{j,s} = - \epsilon_{j',s'}$).
Moreover it is known (see for instance {{\cite[Proposition 6.11 (2)-(iii)]{KKOP21}}} and references therein) that for every $(j,s) \in \Ixi$ one has $\varphi(j,s) = \varphi(j,s+2h) + 2N$. Therefore we have 
 \begin{equation} \label{doubleheart}
 \varphi^{-1}([1,2N]) = \{ (j,s) \in \Ixi \mid \xi(j) \geq s \geq \xi(j)-2h+2 \} . 
 \end{equation}
 Moreover, it is also known that if $(j,s) \in \varphi^{-1}([1,2N])$ then $\epsilon_{j,s}=1$ if and only if $(j,s) \in \varphi^{-1}([1,N])$, and in this case we have $\tau^{(\xi(j)-s)/2}(\gamma_j) =  \beta_t = s_{i_1} \cdots s_{i_{t-1}} \alpha_{i_t}$ where $t := \varphi(j,s)$ (see {{\cite[Proposition 6.11 (2)-(ii)]{KKOP21}}}). 
  
   \begin{remark}
   Comparing with the notations used in \cite{FOComm.Math.Phys.}, our bijection $\varphi$ corresponds to the projection onto the first component of the bijection $\phi_Q$ defined in {{\cite[Section 2.7]{FOComm.Math.Phys.}}}. Moreover $\epsilon_{j,s} = (-1)^k$  where $k$ is the second component of $\phi_Q(j,s)$. 
    \end{remark}
    

 \subsection{Quantum affine algebras and their representations}
\label{remindquantumaffine}

The {\sl quantum affine algebra} $U_q(\widehat{\mathfrak{g}})$ is a Hopf algebra that is a $q$-deformation of the universal enveloping algebra of $\widehat{\mathfrak{g}}$ \cite{Drin87, Jim85}. In this paper, we take $q$ to be a non-zero complex number which is not a root of unity. 


Denote by $\mathcal{P}$ the free abelian group generated by $Y_{i, a}^{\pm 1}$, $i \in I$, $a \in \CC^{\times}$, denote by $\mathcal{P}^+$ the submonoid of $\mathcal{P}$ generated by $Y_{i, a}$, $i \in I$, $a \in \CC^{\times}$. Let $\mathcal{C}$ denote the monoidal category of finite-dimensional representations of the quantum affine algebra $U_q(\widehat{\mathfrak{g}})$.


Any finite dimensional simple object in $\mathcal{C}$ is a highest $l$-weight module with a highest $l$-weight $\mathfrak{m} \in \mathcal{P}^+$, denoted by $L(\mathfrak{m})$ (cf. \cite{CP95a}). The elements in $\mathcal{P}^+$ are called {\sl dominant monomials}.  

Frenkel-Reshetikhin \cite{FR1} introduced the $q$-character map which is an injective ring morphism $\chi_q$ from the Grothendieck ring of $\mathcal{C}$ to $\mathbb{Z}\mathcal{P} = \mathbb{Z}[Y_{i, a}^{\pm 1}]_{i\in I, a\in \mathbb{C}^{\times}}$. For a $U_q(\widehat{\mathfrak{g}})$-module $V$, $\chi_q(V)$ encodes the decomposition of $V$ into common generalized eigenspaces for the action of a large commutative subalgebra of $U_q(\widehat{\mathfrak{g}})$ (the loop-Cartan subalgebra). These generalized eigenspaces are called $l$-weight spaces and generalized eigenvalues are called $l$-weights. One can identify $l$-weights with monomials in $\mathcal{P}$ \cite{FR1}. Then
the $q$-character of a $U_q(\widehat{\mathfrak{g}})$-module $V$ is given by (cf. \cite{FR1})
\begin{align*}
\chi_q(V) = \sum_{  \mathfrak{m} \in \mathcal{P}} \dim(V_{\mathfrak{m}}) \mathfrak{m} \in \mathbb{Z}\mathcal{P},
\end{align*}
where $V_{\mathfrak{m}}$ is the $l$-weight space with $l$-weight $\mathfrak{m}$. 



Let $\mathcal{Q}^+$ be the monoid generated by
\begin{align} \label{eq:Aia}
A_{i, a} = Y_{i, aq}Y_{i, aq^{-1}} \prod_{j \in I , i \sim j} Y_{j,a}^{-1}, \quad i \in I, \ a \in \mathbb{C}^{\times}.
\end{align}
There is a partial order $\preccurlyeq$ on $\mathcal{P}$ (cf. \cite{FR1}) defined by 
$ \mathfrak{m}' \preccurlyeq \mathfrak{m} \text{ if and only if } \mathfrak{m}{\mathfrak{m}'}^{-1} \in \mathcal{Q}^{+}$. 
For any $\mathfrak{m} \in  \mathcal{P}^{+}$, one has 
$$ \chi_q(L(\mathfrak{m})) = \mathfrak{m} \left(1 + \sum_{\mathfrak{m}' \prec \mathfrak{m}} a_{\mathfrak{m}, \mathfrak{m}'} \mathfrak{m}' \right) $$
where only finitely many terms are non zero. 
For $i \in I$, $a \in \mathbb{C}^{\times}$, $k \in \ZZ_{\ge 1}$, the modules 
$$ X_{i,a}^{(k)} := L(Y_{i,a} Y_{i,aq^2} \cdots Y_{i,aq^{2k-2}}) $$
 are called {\sl Kirillov-Reshetikhin modules}. The modules $X_{i,a}^{(1)} = L(Y_{i,a})$ are called {\sl fundamental modules}. 
 
 \subsection{Categorification of cluster algebras}
  \label{remindCQ}

Recall the indeterminates $Y_{i,a}, i \in I, a \in \mathbb{C}^{\times}$ from the previous Section. As in the works of Hernandez-Leclerc \cite{HL,HL15,HLJEMS}, we will only be considering shift parameters $a$ such that $a \in q^{\mathbb{Z}}$ and therefore we will simply write $Y_{i,p}$ instead of $Y_{i,q^p}$ for every $p \in \mathbb{Z}$. In the same way we have 
\begin{equation} \label{defnAip}
A_{i,p} = Y_{i,p+1}Y_{i,p-1} \prod_{j \in I , i \sim j} Y_{j,p}^{-1}, \quad i \in I, p \in \mathbb{Z}.
\end{equation}
Following \cite{KKOP21}, we consider the smallest monoidal  subcategory $\Cxi$ of $\C$ containing all fundamental representations $L(Y_{i,p}), (i,p) \in \Ixi$ and stable under taking subquotients and extensions (this category was denoted $\C^{-}$ in \cite{HLJEMS} where it was introduced for certain choices of height functions $\xi$). The Kirillov-Reshetikhin modules belonging to $\Cxi$ are the $X_{i,p}^{(k)}$ such that  $(i,p) \in \Ixi$ and $ 1 \leq k \leq 1 + (\xi(i)-p)/2$. As shown by the results below from \cite{HL15,HLJEMS}, a special subfamily of these simple objects play a distinguished role from the perspective of the cluster theory, namely 
$$ X_{i,p} := X_{i,p}^{1 + (\xi(i)-p)/2} = L(Y_{i,p} Y_{i,p+2} \cdots Y_{i,\xi(i)}) \quad (i,p) \in \Ixi . $$
Constructed by Hernandez-Leclerc in a former work \cite{HL15}, the category $\CQ$ is  the monoidal subcategory of $\Cxi$ generated by the fundamental representations $L(Y_{i,p}), (i,p) \in I_Q$. The Kirillov-Reshetikhin modules belonging to $\CQ$ are the $X_{i,p}^{(k)}$ such that $(i,p) \in I_Q$ and $ 1 \leq k \leq 1 + (\xi(i)-p)/2$.
 
 Recall the bijection $\varphi$ from Section~\ref{remindAR}. One of the main results of \cite{HL15} is the following:
  
  \begin{thm}[{{\cite[Theorem 6.1]{HL15}}}] \label{thmHL15}
  There is an algebra isomorphism $\mathbb{C} \otimes K_0(\CQ) \simeq \CN$ inducing a bijection from the set of isomorphism classes of simple objects in $\CQ$ to the elements of the dual canonical basis of $\CN$. Furthermore, under this isomorphism, one has 
$$ \forall (i,p) \in I_Q, \enspace  {r_{\varphi(i,p)}^{\ast}}^{\iQ} = [L(Y_{i,p})] \quad \text{and} \quad  x_{\varphi(i,p)}^{\iQ} = [X_{i,p}] . $$ 
   \end{thm}
   
    The next statement deals with the larger category $\Cxi$. 
    
    \begin{thm}[{{\cite[Theorem 5.1]{HLJEMS}}}] \label{thmHL16}
    The complexified Grothendieck ring $\mathbb{C} \otimes  K_0(\Cxi)$ is isomorphic to a cluster algebra $\Axi$, with an initial cluster given by the classes of the Kirillov-Reshetikhin modules $X_{i,p} , (i,p) \in \Ixi$. 
   \end{thm}
   
 The exchange quiver associated to the cluster given by Theorem~\ref{thmHL16} is explicitly constructed in \cite{HLJEMS}. It is checked in \cite{KKOP21} (see {{\cite[Proposition 7.27]{KKOP21}}}) that this quiver is essentially the same as the exchange quivers considered by Berenstein-Fomin-Zelevinsky \cite{BFZ} and Geiss-Leclerc-Schr\"oer \cite{GLS}. Hence by analogy with the standard seeds $\s^{\mathbf{i}}$ in $\CN \simeq \mathbb{C} \otimes K_0(\CQ)$ (see Section~\ref{remindCN}), we will denote by $Q^{\wiQ}$ this exchange quiver, and by $\s^{\wiQ}$ the seed of $\Axi$ given by 
 $$ \s^{\wiQ} = \left( \left( x_1^{\wiQ}, x_2^{\wiQ}, \ldots   \right) , Q^{\wiQ} \right) \quad \text{with $x_t^{\wiQ} = [X_{\varphi^{-1}(t)}]$ for each $t \geq 1$.} $$
 In particular, if $1 \leq t \leq N$, then the flag minor $x_t^{\iQ} \in \CN$ is identified with the cluster variable $x_t^{\wiQ} \in \Axi$ via the injection 
 $$ \CN \simeq \mathbb{C} \otimes K_0(\CQ) \hookrightarrow \mathbb{C} \otimes K_0(\Cxi) \simeq  \Axi . $$

 \subsection{Truncated $q$-characters and $T$-systems}
 \label{remindTsystems}
 
 Using the notations of Section~\ref{remindAR}, we let $\Yxi$ and $\YQ$ denote the subtori of $\mathcal{Y}$ given by
$$ \Yxi := \mathbb{Z}[ Y_{i,p}^{\pm 1} , (i,p) \in \Ixi]  \quad \text{and} \quad \YQ := \mathbb{Z}[ Y_{i,p}^{\pm 1} , (i,p) \in I_Q]  \subset \Yxi .  $$ 
  A useful tool to study the structure of the objects of $\Cxi$ or $\CQ$ is the notion of \textit{truncated $q$-character}, a truncated version of Frenkel-Reshetikhin $q$-character \cite{FR1}. It is an algebra homomorphism 
\begin{align*}
\tchi_q : K_0(\Cxi) \longrightarrow \Yxi 
\end{align*}
such that for every object $M$ in $\Cxi$, the truncated $q$-character $\tchi_q(M)$ is obtained from $\chi_q(M)$ by removing all monomials which do not belong to $\Yxi$. It is proved in \cite{HLJEMS} that $\tilde{\chi}_q$ is injective. This morphism restricts to an embedding 
$$ K_0(\CQ) \longrightarrow \YQ $$
that we still denote $\tchi_q$. It is known that the truncated $q$-characters of the modules $X_{i,p}, (i,p) \in \Ixi$ are reduced to a single term namely their dominant monomial $Y_{i,p} \cdots Y_{i, \xi(i)}$ (this is of course not true anymore for the other Kirillov-Reshetikhin modules).


It is shown in \cite{NakaAnn.Math, H06} that $q$-characters of Kirillov-Reshetikhin modules satisfy $T$-system relations. Therefore the truncated $q$-characters of Kirillov-Reshetikhin modules also satisfy $T$-system relations:
\begin{align} \label{eq:T-system relations}
& \tchi_q(L(X_{i,p}^{(k)})) \tchi_q(L(X_{i,p-2}^{(k)})) = \tchi_q(L(X_{i,p-2}^{(k+1)}))  \tchi_q(L(X_{i,p}^{(k-1)})) + \prod_{j \sim i} \tchi_q(L(X_{j,p-1}^{(k)})).  
\end{align}

\section{Quantum Cartan matrices}
 \label{remindCtilde}
 
  Let $\mathfrak{g}$ be of simply-laced type and let $C(z)$ be the corresponding quantum Cartan matrix, given by
  $$ C_{i,j}(z) := \begin{cases} 
            z+z^{-1} & \text{ if $i=j$,} \\
            -1 & \text{if $i \sim j$,} \\
            0 & \text{otherwise.}
            \end{cases}
     $$
This matrix is invertible and we denote by $\tilde{C}(z)$ its inverse. For each $(i,j) \in I^2$, we let $\tilde{C}_{i,j}(z)$ denote the entry of the matrix $\tilde{C}(z)$ in position $(i,j)$. For every $m \geq 1$ we define $\tilde{C}_{i,j}(m)$ as the coefficient of the term of degree $m$ in  the expansion of $\tilde{C}_{i,j}(z)$ as a power series in $z$, i.e. 
$$ \tilde{C}_{i,j}(z) = \sum_{m \geq 1} \tilde{C}_{i,j}(m) z^m . $$
By convention we extend this definition to all integers by setting $\tilde{C}_{i,j}(m) := 0$ if $m \leq 0$. It is a well-known fact (that can be straightforwardly deduced from the definition) that the $\tilde{C}_{i,j}(m)$ satisfy the following relations:
\begin{equation} \label{relationsforCtilde}
 \begin{cases}
 \tilde{C}_{i,j}(m+1) + \tilde{C}_{i,j}(m-1) - \sum_{k \sim j} \tilde{C}_{i,k}(m) = 0 \quad \text{for any $m \geq 1$} \\
 \tilde{C}_{i,j}(1) = \delta_{i,j}. 
 \end{cases}
 \end{equation}
The following important result  is due to Hernandez-Leclerc \cite{HL15}. We state it using the bijection $\varphi$ and the signs $\epsilon_{i,p}$ introduced in Section~\ref{remindAR}. 

 \begin{thm}[{{\cite[Proposition 2.1]{HL15}}}] \label{thmHLCtilde}
 Let  $(i,p)$ and $(j,s)$ be two elements of $\Ixi$ and assume $s \geq p$. Then one has 
 $$ \tilde{C}_{i,j}(s-p+1) = \epsilon_{i,p} \epsilon_{j,s} \left\langle \beta_{\varphi(i,p)} , \beta_{\varphi(j,s)} \right\rangle_Q . $$ 
  \end{thm}
  
   The following consequence will be also useful for us, especially for the computations we perform in Section~\ref{DtildeKRmonotonic}. It can be straightforwardly deduced  from Theorem~\ref{thmHLCtilde} using the expression of the Cartan pairing in terms of Euler-Ringel forms (see Section~\ref{remindAR}) as well as the well-known identity $ \langle \beta,\gamma \rangle_Q = - \langle \tau_Q^{-1}(\gamma),\beta \rangle_Q$. 
  
   \begin{corollary}[{{\cite[Proposition 3.2]{HL15}}}] \label{cormultiplicityDtilde}
 For any $(i,p), (j,s) \in (\Ixi)^2$ define $\mathcal{N}(i,p;j,s) :=  \tilde{C}_{i,j}(s-p+1) - \tilde{C}_{i,j}(s-p-1)$. Then one has 
 \begin{equation} \label{multiplicityDtilde}
 \mathcal{N}(i,p;j,s) =
\begin{cases}
 \epsilon_{i,p} \epsilon_{j,s}  \left(  \beta_{\varphi(i,p}  , \beta_{\varphi(j,s)} \right) & \text{if $s>p$,} \\
 \delta_{i,j} \delta_{p,s} & \text{otherwise.}
 \end{cases} 
 \end{equation}
    \end{corollary}
    
 \smallskip
 
 We conclude this section with the following elementary property of the coefficients $\tilde{C}_{i,j}(m)$ that will be useful in Sections~\ref{ABCforDtilde} and~\ref{sectionbarDQ}. 
 
  First of all, let us denote by $d(i,j)$ the length of the shortest (non oriented) path from $i$ to $j$ in the Dynkin diagram of $\mathfrak{g}$ (this makes sense as it is a connected acyclic graph). In particular $d(i,i)=0$ for any $i$ and $d(i,j)=1$ if $i \sim j$. 
 
  \begin{lemma} \label{littleLemmaforCtilde}
  Let $i,j \in I$ and let $m \in \mathbb{N}_{\geq 1}$. Assume that $m \leq d(i,j)$. Then one has $\tilde{C}_{i,j}(m)=0$. 
   \end{lemma} 
   
    \begin{proof}
    We prove by strong induction on $m \geq 1$ the statement 
    $$  \forall i,j \quad d(i,j) \geq m  \enspace \Rightarrow \enspace  \tilde{C}_{i,j}(m)=0 . $$
    If $m=1$ then this amounts to prove that if $i \neq j$ then $\tilde{C}_{i,j}(1)=0$. But this follows from the second equality of~\eqref{relationsforCtilde}. Let $m \geq 1$ and  assume the desired statement holds for all $m'$ such that $1 \leq m' \leq m$. Let $i,j \in I$ such that $d(i,j) \geq m+1$. Then the first relation of~\eqref{relationsforCtilde} yields
    $$ \tilde{C}_{i,j}(m+1)= - \tilde{C}_{i,j}(m-1) + \sum_{k \sim j} \tilde{C}_{i,k}(m) . $$
 One has $d(i,j) \geq m+1 > m-1$ hence $\tilde{C}_{i,j}(m-1)=0$ by the induction assumption. Moreover for each $k \sim j$, one has $d(i,k) = d(i,j) \pm 1$ and hence  $d(i,k) \geq m$. Thus the induction assumption again yields $\tilde{C}_{i,k}(m)=0$ for each $k \sim j$. This proves the Lemma. 
    \end{proof}
    
     \begin{remark} \label{rkCtilde}
     With similar arguments, one can also prove  that $\tilde{C}_{i,j}(d(i,j)+1)=1$ for any $i,j \in I$. 
      \end{remark}
      
      \begin{example}
      Consider $\mathfrak{g}$ of type $A_3$. Then the series $\tilde{C}_{i,j}(z), i,j \in \{1,2,3\}$ are given as follows:
       \begin{align*}
 \tilde{C}_{1,1}(z) &= z-z^7+z^9-z^{15} + \cdots \\
 \tilde{C}_{1,2}(z) &= z^2-z^6+z^{10}-z^{14} + \cdots \\
  \tilde{C}_{1,3}(z) &= z^3-z^5+z^{11}-z^{13} + \cdots \\
 \tilde{C}_{2,2}(z) &= z+z^3-z^5-z^7+z^9+z^{11}-z^{13}-z^{15} + \cdots \\
  \tilde{C}_{2,3}(z) &= z^2-z^6+z^{10}-z^{14} + \cdots \\
 \tilde{C}_{3,3}(z) &= z-z^7+z^9-z^{15} + \cdots \\
 \end{align*} 
       \end{example}

\section{Representation theory of quiver Hecke algebras}
  \label{sectionremindKLR}
  
 This section is devoted to some reminders on quiver Hecke algebras and their finite-dimensional representations. We will mostly focus on certain distinguished families of representations, such as the cuspidal modules and the (strongly) homogeneous modules following \cite{BKM,KR,KRhom}. Although a large part of the content of this section remains valid in non-simply-laced types, we will restrict ourselves to the setting of Section~\ref{remindCN}, and refer to \cite{BKM,KR} for a more general exposition. 
 
  \subsection{Quiver Hecke algebras}
  \label{remindKLR}
  
 Let $\mathcal{M}$ denote the set of all finite words over the alphabet $I$. For any such word $\mathbf{j} = (j_1, \ldots , j_d)$, the  \textit{weight} of $\mathbf{j}$ is defined as  
$$ \wt(\mathbf{j}) :=  \sum_{i \in I}  \sharp \{k, j_k = i \} \alpha_i  \in \Gamma_{+} . $$
  Quiver Hecke algebras are defined as a family $\{ R(\beta) , \beta \in \Gamma_{+} \}$ of associative $\mathbb{C}$-algebras indexed by $\Gamma_{+}$. For every $\beta \in \Gamma_{+}$, the algebra $R(\beta)$ is generated by three kind of generators: there are polynomial generators $x_1 , \ldots , x_n$, braiding generators $\tau_1 , \ldots , \tau_{n-1}$, and idempotents $e(\mathbf{j}) , \mathbf{j} \in \text{Seq}(\beta)$ where $\text{Seq}(\beta)$ is the finite subset of $\mathcal{M}$ given by
 $$ \text{Seq}(\beta) := \{ \mathbf{j} \in \mathcal{M} \mid \wt(\mathbf{j}) = \beta \} . $$
 The idempotent generators commute with the polynomial ones and are orthogonal to each other in the sense that $ e(\mathbf{j})e(\mathbf{j}') = \delta_{\mathbf{j},\mathbf{j}'} e(\mathbf{j})$.
For each $\beta \in \Gamma_{+}$, one can consider the category $R(\beta)-mod$ of finite dimensional $R(\beta)$-modules, as well as
 $$ R-mod := \bigoplus_{\beta} R(\beta)-mod . $$
 The category $R-mod$ can be endowed with a structure of a monoidal category via a monoidal product $\circ$ constructed as a parabolic induction. Therefore the Grothendieck group $K_0(R-mod)$ has a ring structure. 
The following results are the main properties of quiver Hecke algebras:

 \begin{thm}[Khovanov-Lauda \cite{KL}, Rouquier \cite{R}] \label{firstthmKLR}
 There is an algebra isomorphism 
 $$  \mathbb{C} \otimes K_0(R-mod) \xrightarrow[]{\simeq} \CN . $$
 \end{thm}
 
 \begin{thm}[Rouquier \cite{R}, Varagnolo-Vasserot \cite{VV}] \label{secondthmKLR}
 The above isomorphism induces a bijection between the set of classes of simple objects in $R-mod$ and the dual canonical basis of $\CN$. 
 \end{thm}
 
  \subsection{Irreducible finite-dimensional representations}
  \label{remindsimpleKLR}
  
  This subsection is devoted to some reminders about the main results of classification of simple objects in the category $R-mod$ associated to a finite-type simple Lie algebra $\mathfrak{g}$, due to Kleshchev-Ram \cite{KR}, McNamara \cite{McNamarafinite} and Brundan-Kleshchev-McNamara \cite{BKM}. We recall in particular the notion of \textit{cuspidal} representation with respect to any fixed convex ordering on the set of positive roots $\Phi_{+}$. 
  
  We assume $<$ is an arbitrary convex ordering on $\Phi_{+}$, and we let $(i_1, \ldots , i_N)$ denote the corresponding reduced expression of $w_0$, the longest element of the Weyl group $W$ of $\mathfrak{g}$.  Then one has $\Phi_{+} = \{ \beta_1 < \cdots < \beta_N \}$ with 
  $$ \beta_k = s_{i_1} \cdots s_{i_{k-1}}(\alpha_{i_k}) $$
  for every $1 \leq k \leq N$. Recall from Section~\ref{remindCN} the dual root vectors $ {r_j^{\ast}}^{\mathbf{i}} \in \CN$ for each $1 \leq j \leq N$. It was proved by McNamara \cite{McNamarafinite} that there exists a family of simple modules $\{ S_{\beta} , \beta \in \Phi_{+} \}$ in $R-mod$, unique up to isomorphism, such that $[S_{\beta_j}] = {r_j^{\ast}}^{\mathbf{i}}$ for every $1 \leq j \leq N$. The module $S_{\beta}$ is called the \textit{cuspidal} module associated to $\beta$ (with respect to the chosen convex ordering $<$ on $\Phi_{+}$).

   Generalizing Leclerc's algorithm ({{\cite[Section 4.3]{Leclerc}}}), Brundan-Kleshchev-McNamara \cite{BKM} describe a procedure producing a word $\mathbf{j}_{\beta} \in  \mathcal{M}$ for every positive root $\beta \in \Phi_{+}$, which we now briefly recall. The crucial tool, that will be useful in the sequel of the present paper, is the notion of \textit{minimal pair}.
   
    \begin{defn}[McNamara \cite{McNamarafinite}, Brundan-Kleshchev-McNamara \cite{BKM}]
    Let $\beta \in \Phi_{+}$ be a positive root. A pair of positive roots $(\gamma, \delta) \in \Phi_{+}^2$ with $\gamma < \delta$ is called a minimal pair for $\beta$ if $\gamma +  \delta = \beta$ and  there is no other pair $(\gamma',\delta')$ such that $\gamma' +  \delta' = \beta$ and  $\gamma < \gamma' < \beta < \delta' < \delta$. 
     \end{defn}
  
  Let us now fix a choice of a minimal pair $(\gamma_{\beta} , \delta_{\beta})$ for each positive root $\beta$. One inductively defines the words $\mathbf{j}_{\beta}$ as follows. For each $i \in I$, set $\mathbf{j}_{\alpha_i} := (i)$. If $\beta \in \Phi_{+}$ is not a simple root, then $\mathbf{j}_{\beta} := \mathbf{j}_{\gamma_{\beta}} \mathbf{j}_{\delta_{\beta}}$. This yields a finite collection of words, in bijection witht the set of positive roots of $\mathfrak{g}$. In the case considered in \cite{KR}, where the order $<$ arises from a total order on $I$, the words  $\mathbf{j}_{\beta_1} , \ldots ,  \mathbf{j}_{\beta_N}$ are called \textit{good Lyndon words}.

We can now state the main classification result.

\begin{thm}[Kleshchev-Ram \cite{KR}, McNamara \cite{McNamarafinite}, Brundan-Kleshchev-McNamara \cite{BKM}] \label{thmrootpartitions}
 There is a bijection between the set of isomorphism classes of  simple objects in $R-mod$ and the set $\mathbb{N}^{\Phi_{+}}$, given by
 $$ \mathbf{c} := (c_1, \ldots , c_N) \in \mathbb{N}^{\Phi_{+}}  \longmapsto  L(\mathbf{c}) := \text{hd} \left( S_{\beta_N}^{\circ c_N} \circ  \cdots \circ S_{\beta_1}^{\circ c_1} \right).
   $$
   Moreover, for each $(c_1, \ldots , c_N)  \in \mathbb{N}^{\Phi_{+}}$, one has 
  $$ \dim_{\mathbb{C}} \left( e({\bf j}_{\beta_N}^{c_N} \cdots {\bf j}_{\beta_1}^{c_1}) L(\mathbf{c}) \right) = 1 . $$
  \end{thm}
  
  In this statement, $\text{hd}(M)$ stands for the head of a module $M$, i.e. the quotient of $M$ by its radical (the intersection of its maximal submodules).
 
    \subsection{Short exact sequences in $R-mod$} 
    \label{remindBKM}
    
  In this paragraph, we recall an important result proved in \cite{BKM} as the \textit{length two property}. It will play a crucial role in Section~\ref{Dbarcuspidals} for our computations of the images of certain dual root vectors under Baumann-Kamnitzer-Knutson's morphism $\barD$, especially in type $D_n , n \geq 4$. 
  
  \smallskip
  
   We fix an arbitrary convex ordering $<$ on $\Phi_{+}$. 
   
    \begin{thm}[{{\cite[Theorem 4.7]{BKM}}}]  \label{thmBKMshortexactseq}
    Let $\beta \in \Phi_{+}$ and let $(\gamma, \delta)$ be a minimal pair for $\beta$. Let $\mathbf{c} = (c_1, \ldots , c_N)$ be the $N$-tuple of integers defined by $c_k := 1$ if $ \beta_k \in \{ \gamma ,\delta \}$ and $c_k := 0$ otherwise. Then one has a short exact sequence in $R-mod$:
$$ 0 \longrightarrow S_{\beta} \longrightarrow S_{\delta} \circ S_{\gamma}  \longrightarrow L(\mathbf{c}) \longrightarrow 0 . $$
    \end{thm}

 \begin{remark}
 This short exact sequence is an ungraded version of Brundan-Kleshchev-McNamara's statement, but it will be sufficient for our purpose. 
 \end{remark}
 
\subsection{Homogeneous modules over quiver Hecke algebras}
\label{remindKRHom3}

 In this paragraph, we recall Kleshchev-Ram's construction of simple homogeneous representations of simply-laced type quiver Hecke algebras. We begin with some reminders on the combinatorics of fully-commutative elements of Weyl groups following Stembridge \cite{Stem}.
 For any $w \in W$, we will denote by $\text{Red}(w)$ the set of all reduced expressions of $w$. 
 
  \bigskip
  
 For $w \in W$ and $\mathbf{i} = (i_1, \ldots , i_N) \in \text{Red}(w)$, one can define an infinite sequence $ \hat{\mathbf{i}} := (i_1, i_2, \ldots)$ exactly as in Section~\ref{remindAR}. Then using the notation $k_{+}$ introduced in Section~\ref{remindAR}, for every $1 \leq k \leq N$ we have that $k_{+} > N$ if and only if $k$ is the position of the last occurrence of $i_k$ in $\mathbf{i}$. 

 The following definition is essentially due to Stembridge \cite{Stem} relying on former constructions by Proctor \cite{P2}. Here we write it in a way suited to simply-laced cases. 
 
  \begin{defn}[Stembridge \cite{Stem}] \label{defminuscule}
 Let $w \in W$. 
  \begin{itemize}
  \item The element $w$ is called fully-commutative if for every $\mathbf{i} \in \text{Red}(w)$ and $1 \leq k \leq N$, one has 
$$ k_{+} \leq N \enspace \Rightarrow \enspace \sharp \{ l \mid k<l<k_{+}, i_k \sim i_l \} \geq 2 . $$   
   \item The element $w$ is called dominant minuscule if for every $\mathbf{i} \in \text{Red}(w)$ and $1 \leq k \leq N$, one has
$$ k_{+} \leq N \enspace \Rightarrow \enspace \sharp \{ l \mid k<l<k_{+}, i_k \sim i_l \} = 2 \quad \text{and} \quad  k_{+} > N \enspace \Leftrightarrow \enspace  \sharp \{ l \mid  l>k, i_k \sim i_l \} = 1 . $$
  \end{itemize}
   \end{defn} 
   
     We will denote by $\mathcal{FC}$ (resp. $\mathcal{M}in^{+}$) the set of fully-commutative (resp. dominant minuscule) elements of $W$. Note that $\mathcal{M}in^{+} \subset \mathcal{FC}$.
   
   \smallskip
   
  We now recall the construction of simple homogeneous representations following Kleshchev-Ram \cite{KRhom}. 
   
    \begin{thm}[{{\cite[Theorem 3.6]{KRhom}}}] \label{thmKRhom}
     For every $w \in \mathcal{FC}$ and for every  $\mathbf{i} \in \text{Red}(w)$, there exists a unique simple module $S(\mathbf{i})$ in $R-mod$ such that $\dim e(\mathbf{i}) S(\mathbf{i}) = 1$. Moreover, if $\mathbf{i}$ (resp. $\mathbf{i}'$) is a reduced expression of $w$ (resp. $w'$), then the modules $S(\mathbf{i})$ and $S(\mathbf{i}')$ are isomorphic in $R-mod$ if and only if $w=w'$. 
     \end{thm}
     
   For each $w \in \mathcal{FC}$, we denote by $S(w)$ the module $S(\mathbf{i})$ for an arbitrary reduced expression $\mathbf{i}$ of $w$. 

  
\begin{remark}
   The modules $S(w), w \in \mathcal{FC}$ are called \textit{homogeneous}. This is due to the fact the quiver Hecke algebras $R(\beta), \beta \in \Gamma_{+}$ carry a natural $\mathbb{Z}$-grading, and the modules $S(w)$ are precisely those which are concentrated in a single degree for this grading. 
   \end{remark}
   
   The following distinguished family of homogeneous representations will be of particular interest for us, especially in Proposition~\ref{lemDbarKLR} below. 
   
   \begin{defn}[Kleshchev-Ram \cite{KRhom}] \label{defstronghomog}
   The simple modules $S(w)$ for $w \in \mathcal{M}in^{+}$ are called  \textit{strongly homogeneous}. 
   \end{defn}
 
  \section{Baumann-Kamnitzer-Knutson's morphism $\barD$}
   \label{remindDbar}
  
  This section is devoted to several reminders on the definition and some of the main properties of the algebraic morphism $\barD$ recently introduced by Baumann-Kamnitzer-Knutson \cite{BKK}. We also recall certain results from the first author \cite{Casbi3} that will be needed in Section~\ref{sectionproofs}, in particular the propagation result {{\cite[Theorem 5.6]{Casbi3}}} (Theorem~\ref{thmpropagation} below) which will be involved in the proof of the second main result of this paper (Theorem~\ref{mainthm2}). 
  
  \subsection{Geometric Satake Correspondence}
 \label{remindMV}
 
  Throughout this section $\mathbf{G}$ denotes a simple simply-connected group, $P$ stands for the weight lattice and $W$ the Weyl group of $\mathbf{G}$. Let $\mathbf{G}^{\vee}$ denote the Langlands dual of $\mathbf{G}$, fix a Borel subgroup $\mathbf{B}^{\vee}$ in $\mathbf{G}^{\vee}$ and a maximal torus $\mathbf{T}^{\vee}$ in $\mathbf{B}^{\vee}$. Furthermore for every $\lambda \in P^{+}$ we let $L(\lambda)$ denote the finite-dimensional irreducible representation of $\mathbf{G}$ of highest weight $\lambda$, and $L(\lambda)_{\mu}$ denote its weight subspace of weight $\mu$ for any $\mu \in P$. 
  

  We set $ \mathcal{O} := \mathbb{C}[[t]]$ and $\mathcal{K} := \mathbb{C}((t))$.
   The affine Grassmannian $Gr_{\mathbf{G}^{\vee}}$ of $\mathbf{G}^{\vee}$ is defined as 
   $$ Gr_{\mathbf{G}^{\vee}} :=  \mathbf{G}^{\vee}(\mathcal{K}) / \mathbf{G}^{\vee}(\mathcal{O}) . $$
   There is a natural action of $\mathbf{T}^{\vee}(\mathbb{C})$ on $Gr_{\mathbf{G}^{\vee}}$ whose locus of fixed points is given by a collection $\{ L_{\mu} , \mu \in P \}$ of points in $Gr_{\mathbf{G}^{\vee}}$ indexed by the weight lattice of $\mathbf{G}$.
 For each $(\lambda,\mu) \in P^{+} \times P$, Mirkovi\'c-Vilonen \cite{MV} constructed  a closed subvariety $\mathcal{MV}^{\lambda,\mu}$ of $Gr_{\mathbf{G}^{\vee}}$ such that there is an isomorphism of vector spaces
  \begin{align} \label{isomMV}
   H_{\bullet}(\mathcal{MV}^{\lambda,\mu}) \simeq L(\lambda)_{\mu} .
   \end{align}
 The irreducible components of $\mathcal{MV}^{\lambda,\mu}$ are called the \textbf{MV cycles} of type $\lambda$ and of weight $\mu$. For every $\lambda \in P^{+}$, the images under the isomorphism~\eqref{isomMV} of the homology classes of all MV cycles of type $\lambda$ and of weight $\mu$ ($\mu \in P$) form a basis of $L(\lambda)$, called the \textit{MV basis of $L(\lambda)$}. Using the classical injections from $L(\lambda)$ to the coordinate ring $\CN$ (see for example {{\cite[Section 2.5]{BKK}}}), one can then build a basis of $\CN$ out of the MV bases of all the simple representations $L(\lambda)$ of $\mathbf{G}$, called the \textit{MV basis of $\CN$}, whose elements are indexed by certain MV cycles called stable MV cycles. We denote by $b_Z$ the element of the MV basis corresponding to the stable MV cycle $Z$.

    \subsection{Equivariant multiplicities}
   \label{remindequivmulti}
  
One of Baumann-Kamnitzer-Knutson's main motivations was Muthiah's conjecture \cite{Muthiah} stating the $W$-equivariance of a certain map  $L(\lambda) \longrightarrow \mathbb{C}(\mathbf{T})$. The proof of \cite{BKK} relies on the notion of \textit{equivariant multiplicities} developped by Brion \cite{Brion} out of former constructions due to Joseph \cite{Joseph} and Rossmann \cite{Rossmann}. 

  Given a closed projective scheme $X$ together with an action of a torus $T$ on $X$, we let $X^{T}$ denote the set of fixed points of this action and $H_{\bullet}^{T}(X)$ denote the $T$-equivariant homology of $X$. It follows from Brion's results \cite{Brion} that the set of homology classes of the points in $X^{T}$ actually forms a basis of $H_{\bullet}^{T}(X)$. Therefore, for any closed subvariety $Y \subset X$ stable under the action of $T$, one can decompose the class of $Y$ on this basis as 
  $$ [Y] = \sum_{p \in X^{T}} \epsilon_{p}^{T}(Y) [\{p\}] . $$
  The coefficient $\epsilon_{p}^{T}(Y)$ is an element of the field $\mathbb{C}(T)$ of functions on $T$ and is called the \textit{equivariant multiplicity} of $Y$ at $p$. Note that $\epsilon_{p}^{T}(Y) = 0$ if $p \notin Y$ (see {{\cite[Theorem 4.2 (i)]{Brion}}}).

  \subsection{The morphism $\barD$}
  \label{remindBKK}
  
 Baumann-Kamnitzer-Knutson \cite{BKK} used this notion of equivariant multiplicity in the study of the MV basis of $\CN$ via Duistermaat-Heckman measures. With the notations of the previous section, we consider $X := Gr_{\mathbf{G}^{\vee}}$ together with the action of the torus $\mathbf{T}^{\vee}(\mathbb{C})$. As recalled above, the set of fixed points of this action is $\{L_{\mu} , \mu \in P \}$.
 
  The definition of $\barD$ goes as follows. It is known (see for instance \cite{GLSAENS,GLS}) that $\CN$ can be identified with the dual (as a Hopf algebra) of $U(\mathfrak{n})$. For any $f \in \CN$ and $e \in U(\mathfrak{n})$, we will denote by $f(e)$ the canonical pairing between $f$ and $e$. Choose a root vector $e_i \in \mathfrak{n}$ of weight $\alpha_i$ for each $i \in I$. 
 Then Baumann-Kamnitzer-Knutson \cite{BKK} define the following map:
 \begin{equation} \label{deftnDbar}
 \begin{array}{cccc}
  \barD : & \CN  & \longrightarrow & \mathbb{C}(T) = \mathbb{C}(\alpha_1, \ldots , \alpha_n) \\
   {}  & f & \longmapsto & \sum_{\mathbf{j}}  f(e_{j_1} \cdots e_{j_d}) \frac{1}{\alpha_{j_1} (\alpha_{j_1} + \alpha_{j_2}) \cdots (\alpha_{j_1} + \cdots \alpha_{j_d})} . 
  \end{array}
  \end{equation}
Although this sum a priori runs over all arbitrary sequences $\mathbf{j}$ of elements of $I$, it is nevertheless finite as $U(\mathfrak{n})$ acts locally nilpotently on $\CN$. 
 The following statement, which is one of the main results of \cite{BKK}, asserts that the evaluation of $\barD$ on an element $b_{Z}$ of the Mirkovi\'c-Vilonen basis can be related to a certain equivariant multiplicity of the corresponding MV cycle $Z$. 
 
   \begin{thm}[{{\cite[Lemma 8.3, Corollary 10.6]{BKK}}}] \label{thmBKK}
    \begin{enumerate}
    \item The map $\barD$ is an algebra morphism. 
    \item  For any $\mu \in  - \Gamma_{+}$ and any stable MV cycle $Z$ of weight  $\mu$, one has 
   $$ \barD(b_{Z}) = \epsilon_{L_{\mu}}^{\mathbf{T}^{\vee}}(Z) . $$
       \end{enumerate}
   \end{thm}

 The morphism $\barD$ provides a useful tool to compare various remarkable  bases of $\CN$. For instance, the definition of $\barD$ can be conveniently reformulated using the categorification of $\CN$ via modules over the quiver Hecke algebras associated to $\mathfrak{g}$ (see Section~\ref{remindKLR}): for any $\beta \in \Gamma_{+}$, any module $M$ in $R(\beta)-mod$ can be decomposed into weight subspaces:
  $$ M = \bigoplus_{\mathbf{j} \in \text{Seq}(\beta)} e(\mathbf{j}) \cdot M $$
  (we refer to Section~\ref{remindKLR} for the notations). Then one has
  \begin{equation} \label{eqnDbarKLR}
  \barD([M]) = \sum_{ \substack{ \mathbf{j} \in \text{Seq}(\beta) \\ \mathbf{j} := (j_1, \ldots , j_d) } } \dim \left(  e(\mathbf{j}) \cdot M) \right)  \frac{1}{\alpha_{j_1} (\alpha_{j_1} + \alpha_{j_2}) \cdots (\alpha_{j_1} + \cdots + \alpha_{j_d})} . 
  \end{equation}
In Section~\ref{Dbarcuspidals} we will use this to investigate the values of $\barD$ on the elements of the dual canonical basis of $\CN$. 
 A similar expression can be written for the evaluation of $\barD$ on the elements of the dual semi-canonical basis of $\CN$ in terms of representations of preprojective algebras. The dimensions of the weight subspaces of modules in $R-mod$ are then replaced by the Euler characteristics of certain type-$\mathbf{j}$ flag varieties in the terminology of Geiss-Leclerc-Schr\"oer \cite{GLSAENS}. As an application of Theorem~\ref{thmBKK}, Dranowski, Kamnitzer, and Morton-Ferguson show in an appendix of \cite{BKK} that the MV basis and the dual semicanonical basis of $\CN$ are not the same by exhibiting elements of these bases satisfying some compatibility condition (see {{\cite[Definition 12.1]{BKK}}}) but where $\barD$ nonetheless takes different values.

 We conclude this paragraph by recalling from \cite{Casbi3} the following remarkable property of Kleshchev-Ram's strongly homogeneous modules in $R-mod$ (see Section~\ref{remindKRHom3}) involving Baumann-Kamnitzer-Knutson's morphism $\barD$. It can be essentially viewed as a representation-theoretic reformulation of Nakada's colored hook formula \cite{Nakada} using the identity~\eqref{eqnDbarKLR}. 
  
\begin{prop}[{{\cite[Proposition 5.1]{Casbi3}}}] \label{lemDbarKLR}
Let $w$ be a dominant minuscule element in $W$ and $S(w)$ the strongly homogeneous simple module in $R-mod$ corresponding to $w$ under the bijection of Theorem~\ref{thmKRhom}. Then one has 
$$ \barD([S(w)]) = \prod_{\beta \in \Phi_{+}^{w}} \frac{1}{\beta} \qquad \text{where $\Phi_{+}^{w} := \Phi_{+} \cap (-w \Phi_{+})$.} $$
 \end{prop}  
 
  \subsection{The values of $\barD$ on the flag minors of $\CN$}
  \label{remindDbarflagminors}
 
 We now recall some setting from the first author's previous work \cite{Casbi3} and in particular the propagation result {{\cite[Theorem 5.6]{Casbi3}}}. 
 Recall from Section~\ref{remindCN} that there is a distinguished family of cluster variables in $\CN$ called flag minors, grouped into clusters indexed by the set of reduced epressions of $w_0$. The main aim of \cite{Casbi3} was to investigate the  values taken by Baumann-Kamnitzer-Knutson's morphism $\barD$ on the flag minors of $\CN$. It was observed in particular that these values seemed to share a similar form with the images under $\barD$ of the elements of the dual canonical basis corresponding to Kleshchev-Ram's strongly homogeneous modules (see Proposition~\ref{lemDbarKLR} above) or the elements of the MV basis associated to smooth MV cycles. 
 
  As in \cite{Casbi3}, we consider the following properties of the flag minors $x_1^{\mathbf{i}} , \ldots , x_N^{\mathbf{i}}$ for any reduced expression $\mathbf{i} = (i_1, \ldots , i_N)$ of $w_0$:
  \begin{enumerate}
   \item[$({\rm A}_{\bf i})$] For every $1 \leq j \leq N$, one has  
   $$ \overline{D}(x_j^{\mathbf{i}}) = 1/P_j^{\mathbf{i}} $$ 
   where $P_j^{\mathbf{i}}$ is a product of positive roots. 
   \item[$({\rm B}_{\bf i})$] For every $1 \leq j \leq N$ one has  
   $$ P_j^{\mathbf{i}} P_{j_{-}}^{\mathbf{i}}  = \beta_j  \prod_{\substack{ l<j<l_{+} \\ i_l \sim i_j}} P_l^{\mathbf{i}} $$  
   where $\beta_j = s_{i_1} \cdots s_{i_{j-1}}(\alpha_{i_j})$. 
   \item[$({\rm C}_{\bf i})$] For every $j$ such that $j_{+} \leq N$ and every $1 \leq i \leq N$, one has 
   $$ [\beta_i ; P_j^{\mathbf{i}}] - [\beta_i ; P_{j_{+}}^{\mathbf{i}}] \leq 1 $$
   where $[\beta;P]$ stands for the multiplicity of the positive root $\beta$ in the polynomial $P$. 
  \end{enumerate}
  
  The following statement was one of the main results of \cite{Casbi3}. It will be involved in the proof of the second main result of this paper (Theorem~\ref{mainthm2}). 
  
   \begin{thm}[{{\cite[Theorem 5.6]{Casbi3}}}] \label{thmpropagation}
 Let $\mathfrak{g}$ be a simple Lie algebra of simply-laced type.  Assume there exists a reduced expression $\mathbf{i}_0$ such that the standard seed $\s^{\mathbf{i}_0}$ satisfies the three properties $(A_{\mathbf{i}_0}), (B_{\mathbf{i}_0}), (C_{\mathbf{i}_0})$. Then for every reduced expression $\mathbf{i}$ of $w_0$, the properties $(A_{\mathbf{i}}), (B_{\mathbf{i}}), (C_{\mathbf{i}})$ hold for the standard seed $\s^{\mathbf{i}}$ of $\CN$. 
  \end{thm}

\section{Definition and properties of the morphism $\td_{\xi}$}
 \label{sectionDtilde}  
  
   In this section we introduce the main object of the present paper, namely the morphism $\td_{\xi} : \mathbb{C} \otimes \Yxi \longrightarrow \mathbb{C}(\alpha_1, \ldots , \alpha_n)$. The definition of $\td_{\xi}$ involves the coefficients of the inverse of the quantum Cartan matrix associated to $\mathfrak{g}$ (see Section~\ref{remindCtilde}). We state the first main result of this paper (Theorem~\ref{mainthm1}), which relates the restriction $\td_Q$ of $\td_{\xi}$ on the subtorus $\mathbb{C} \otimes \YQ$ to Baumann-Kamnitzer-Knutson's morphism $\barD$ introduced in \cite{BKK}. In this framework, we also provide a general formula for the images under $\td_{\xi}$ of the cluster variables of the seed $\s^{\wiQ}$ of $\Axi$ arising from Hernandez-Leclerc's construction \cite{HLJEMS} and we prove that the obtained rational fractions satisfy a family of remarkable properties, analogous to $(A_{\mathbf{i}}), (B_{\mathbf{i}}), (C_{\mathbf{i}})$ from \cite{Casbi3} (see also Section~\ref{remindDbarflagminors} above). 
    
     \subsection{The morphisms $\td_{\xi}$ and $\td_Q$}
      \label{definitionDtilde}
    
   Let $Q$ be an arbitrary orientation of the Dynkin diagram of a simply-laced Lie algebra $\mathfrak{g}$ and let $\xi$ be a height function adapted to $Q$. Recall from Section~\ref{remindHL15} the set $\Ixi$ (resp. $I_Q$), the torus $\Yxi$ (resp. $\YQ$) containing the truncated $q$-characters of all the representations in the category $\Cxi$ (resp. $\CQ$), as well as the bijection $\varphi : \Ixi \longrightarrow \mathbb{Z}_{\geq 1}$. 
   
    We define the algebra morphism $\td_{\xi}$ from the complexified torus $ \mathbb{C} \otimes \Yxi$ to the field $\mathbb{C}(\alpha_1, \ldots, \alpha_n)$ as follows:
    \begin{equation} \label{defDtildexi}
    \forall (i,p) \in \Ixi, \enspace \td_{\xi}(Y_{i,p}) := \prod_{(j,s) \in \Ixi} \beta_{\varphi(j,s)}^{\tilde{C}_{i,j}(s-p-1) - \tilde{C}_{i,j}(s-p+1)} . 
     \end{equation}
As $\tilde{C}_{i,j}(m) := 0$ if $m \leq 0$, only the couples $(j,s) \in \Ixi, s \geq p$ have a non trivial contribution, and hence this product is finite. 
   We also define the morphism $\td_Q$ as the restriction of $\td_{\xi}$ to the complexified torus $\mathbb{C} \otimes \YQ$.

The composition of the truncated $q$-character morphism $\widetilde{\chi}_q$ (see Section~\ref{remindTsystems}) with Hernandez-Leclerc's isomorphisms (see Theorems~\ref{thmHL15} and~\ref{thmHL16}) yields an embedding
 $$\iota : \Axi \longrightarrow \mathbb{C} \otimes  \Yxi $$
  that restricts to an embedding $\CN \longrightarrow \mathbb{C} \otimes \YQ$ following the commutative diagram
 $$
  \xymatrix{ 
   \CN \ar[r]^{\simeq} \ar@{^{(}->}[d] &  \mathbb{C} \otimes K_0(\CQ) \ar@{^{(}->}[d] \ar@{^{(}->}[r]^{\tchi_q}  & \mathbb{C} \otimes \YQ \ar@{^{(}->}[d] \\ 
   \Axi  \ar[r]^{\simeq} &  \mathbb{C} \otimes K_0(\Cxi)  \ar@{^{(}->}[r]^{\tchi_q}  & \mathbb{C} \otimes \Yxi 
  }
 $$
 We now state the first main result of this work. 
 
 \begin{thm} \label{mainthm1}
Let $\mathfrak{g}$ be a simple Lie algebra of simply-laced type and let $Q$ be an arbitrary orientation of the Dynkin diagram of $\mathfrak{g}$. Then the following diagram commutes:
\begin{align*}
\xymatrix{ \CN \ar[rd]_{\barD} \ar[r]^{\simeq} & \mathbb{C} \otimes K_0(\CQ) \ar@{^{(}->}[r]^{\tchi_q}  & \mathbb{C} \otimes \YQ \ar[ld]^{\tilde{D}_Q} \\
{} &  \mathbb{C}(\alpha_i , i \in I)  & {} } 
\end{align*}
\end{thm}

In other words, Hernandez-Leclerc's categorification allows to embed $\CN$ into the torus $\YQ$ via the (truncated) $q$-character morphism; then $\barD$ can be interpreted as the restriction of $\td_Q$ on $\CN$, viewed as a subalgebra of $\YQ$. 
  The proof of this statement will require several steps which we now briefly describe. We begin in Section~\ref{ABCforDtilde} by establishing a family of remarkable properties satisfied by the values of $\td_{\xi}$ on the cluster variables $x_t^{\wiQ}, t \geq 1$ of the initial seed of $\Axi$ constructed by Hernandez-Leclerc \cite{HLJEMS} (see Section~\ref{remindCQ}). The proofs are valid for any simply-laced type and any orientation $Q$. These properties will play a crucial role in Section~\ref{sectionproofs} for the proofs of Theorem~\ref{mainthm1} as well as the second main result of this paper (Theorem~\ref{mainthm2}). 
  Before that, we investigate in detail the case where $Q$ is a certain specific orientation $Q_0$ of the Dynkin diagram of $\mathfrak{g}$. Sections~\ref{Dbarcuspidals} and~\ref{DtildeKRmonotonic} are respectively devoted to providing explicit formulas in types $A_n$ and $D_n$ for the evaluation of $\barD$ and $\td_{Q_0}$ on the dual root vectors of $\CN$ (with respect to the convex ordering on $\Phi_{+}$ corresponding to $\mathbf{i}_{Q_0}$). Note that the case of type $A_n$ is in fact contained as a subcase of the case of type $D_n$ but for we chose to treat them in distinct subsections, for the sake of readibility. We treat the types $E_6, E_7$ and $E_8$ separately.  We then prove in Section~\ref{sectionproofs} that $\barD$ and $\td_{Q_0}$ coincide on $\CN$. Together with the propagation result from the first author's previous work \cite{Casbi3} (Theorem~\ref{thmpropagation} above), this allows us to prove the second main result of this paper (Theorem~\ref{mainthm2}). The proof is valid for any simply-laced type and any orientation $Q$. We conclude the proof of Theorem~\ref{mainthm1} for an arbitrary orientation by combining this with the properties of $\td_{\xi}$ established in Section~\ref{ABCforDtilde}.

   \subsection{Properties of $\td_{\xi}$ and initial seed for $\Axi$}
    \label{ABCforDtilde}
 
 In this subsection, we consider the initial seed $\s^{\wiQ}$ in the cluster algebra $\Axi$ (see Section~\ref{remindCQ}). Recall that the cluster variables of $\s^{\wiQ}$ are given by $x_1^{\wiQ}, x_2^{\wiQ}, \ldots $ with $\iota(x_t^{\wiQ}) = \tchi_q(X_{\varphi^{-1}(t)})$ for each $t \geq 1$ where $\varphi$ is the bijection introduced in Section~\ref{remindAR}. Throughout the rest of this section, we will simply write $x_t$ for $x_t^{\wiQ}$. We prove that the images of these cluster variables under the morphism $\td_{\xi}$ satisfy properties $(A_{\wiQ}), (B_{\wiQ}), (C_{\wiQ})$ analogous to the properties $(A_{\iQ}), (B_{\iQ}), (C_{\iQ})$ from Section~\ref{remindDbarflagminors}, with $\iQ$ replaced by its infinite analogue $\wiQ$ (see Section~\ref{remindAR}). Whereas the latter properties remained mysterious in \cite{Casbi3}, the former are now naturally deduced from the definition of $\td_{xi}$ using the properties of the coefficients $\tilde{C}_{i,j}(m)$ (Section~\ref{remindCtilde}). Note that each property $(A_{\wiQ}), (B_{\wiQ})$ (resp.  $(C_{\wiQ})$) is an infinite system of equalities (resp. inequalities) indexed  by $\mathbb{Z}_{\geq 1}$, whereas  $(A_{\iQ}), (B_{\iQ}), (C_{\iQ})$ were finite systems, indexed by $\{1, \ldots , N \}$.

  
   \begin{lemma} \label{DtildeinitialKRs}
Let $t \geq 1$ and let $(i,p) := \varphi^{-1}(t) \in \Ixi$. Then one has 
   $$ \td_{\xi} \left( \iota (x_t) \right) =  \prod_{(j,s) \in \Ixi} \frac{1}{ \beta_{\varphi(j,s)}^{\tilde{C}_{i,j}(s-p+1)} }  . $$
 \end{lemma}
 
  \begin{proof}
 Let  us fix $t \geq 1$ and $(i,p) := \varphi^{-1}(t) \in \Ixi$. Recall from Section~\ref{remindHL15} that one has $\tchi_q(X_{i,p}) = Y_{i,p} Y_{i,p+2} \cdots Y_{i,\xi(i)}$. Hence applying the definition of $\td_{\xi}$, we get
  \begin{align*}
\td_{\xi} \left( \iota (x_t) \right) = \td_{\xi} \left( \tchi_q(X_{i,p}) \right) &= \td_{\xi}(Y_{i,p} \cdots Y_{i,\xi(i)})  =\td_{\xi} (Y_{i,p}) \cdots \td_Q(Y_{i,\xi(i)}) \\
  &= \prod_{(j,s) \in \Ixi} \beta_{\varphi(j,s)}^{- \left( \mathcal{N}(i,p;j,s) + \cdots + \mathcal{N}(i,\xi(i);j,s) \right)} 
  \end{align*}
  where $\mathcal{N}(i,t;j,s) := \tilde{C}_{i,j}(s-t+1) -  \tilde{C}_{i,j}(s-t-1)$ as in Section~\ref{remindCtilde}. Obviously one has 
  $$ \mathcal{N}(i,p;j,s) + \cdots + \mathcal{N}(i,\xi(i);j,s) =  \tilde{C}_{i,j}(s-p+1) -  \tilde{C}_{i,j}(s-\xi(i)-1) $$
 In order to conclude, it remains to observe that one always has $\xi(p) - \xi(q) \leq d(p,q)$ for any $p,q \in I$ (with equality if and only if $p \in B(q)$ with the notations of Section~\ref{remindAR}) where $d(p,q)$ is the distance function on $I$ defined in Section~\ref{remindCtilde}. In particular one has $s-\xi(i)-1 \leq \xi(j)-\xi(i)-1 < d(i,j)$ for every $(j,s) \in \Ixi$. Thus Lemma~\ref{littleLemmaforCtilde} implies 
$$ \mathcal{N}(i,p;j,s) + \cdots + \mathcal{N}(i,\xi(i);j,s) = \tilde{C}_{i,j}(s-p+1) $$
which proves the Lemma. 
\end{proof}

  We now prove the main statements of this section (Propositions~\ref{AforDtilde}, ~\ref{BforDtilde} and~\ref{CforDtilde}) which can be seen as analogues of $(A_{\iQ})$, $(B_{\iQ})$ and $(C_{\iQ})$ for the seed $\s^{\wiQ}$ in $\Axi$. Note that  Proposition~\ref{AforDtilde} restricts to the variables $x_t, t \leq 2N$. We postpone the case $t>2N$ to Section~\ref{finalsection} (see Corollary~\ref{corperiodicity}), as it is not strictly necessary for the proof of Theorem~\ref{mainthm2}. 
 
 \begin{prop}[\textbf{Property $(A_{\wiQ})$}] \label{AforDtilde}
 Let $t \in \{1, \ldots , 2N \}$. Then one has 
 $$ \td_{\xi} \left( \iota (x_t) \right) = \prod_{\beta \in \Phi_{+}} \frac{1}{\beta^{n_t(\beta)}} \quad \text{with $n_t(\beta)  \in \mathbb{N}$ for every $\beta \in \Phi_{+}$} . $$  
 Moreover, if $n_t(\beta) \neq 0$ then $n_t(\beta) =  | \langle \beta_t , \beta \rangle_Q |$. 
  \end{prop} 
  
   \begin{proof}
    Let $(i,p) := \varphi^{-1}(t)$ and $\gamma := \beta_t$. As $1 \leq t \leq 2N$, we have $ \xi(i) \geq p > \xi(i)-2h$ by~\eqref{doubleheart}. Let $\beta \in \Phi_{+}$. By Lemma~\ref{DtildeinitialKRs} the multiplicity $n_t(\beta)$ of $\beta$ in $\left( \td_{\xi}(\iota (x_t)) \right)^{-1}$ is
 $$ n_t(\beta) =  \sum_{(j,s) \in I_{p, \beta}} \tilde{C}_{i,j}(s-p+1), \qquad  I_{p, \beta} :=   \{ (j,s) \in \Ixi  \mid s \geq p, \beta_{\varphi(j,s)} = \beta \} .  $$
 If $I_{p,\beta} = \emptyset$ then $n_t(\beta) = 0$ and we are done. Otherwise, let $(j,s) \in I_{p, \beta}$. If $(j,s) \notin \varphi^{-1}([1,2N])$ then $s < \xi(j)-2h+2$ by~\eqref{doubleheart} and hence we have $\xi(i)-2h+2 \leq p \leq s < \xi(j)-2h+2$. In particular  we have $ 0 \leq s-p < \xi(j) -\xi(i) \leq d(i,j)$. By Lemma~\ref{littleLemmaforCtilde} this implies $\tilde{C}_{i,j}(s-p+1) =0$. 
 
  On the other hand for any $(j,s) \in I_{p, \beta} \cap \varphi^{-1}([1,2N])$, one has that for any $m \geq 1$, $(j,s-2mh) \notin \varphi^{-1}([1,2N])$ as $s-2mh \leq \xi(j)-2h$, and $(j,s+2mh) \notin \Ixi$ as $s+2mh > \xi(j)$. Consequently Proposition~\ref{propositionFujitaOh}, implies that $I_{p, \beta} \cap \varphi^{-1}([1,2N])$ is either of the form $\{ (j,s) ; (j^{*},s+h) \}$ or of the form $\{(j,s)\}$, for some $(j,s) \in \Ixi$. 
In the first case, Theorem~\ref{thmHLCtilde} implies
 \begin{align*}
  n_t(\beta) &= \tilde{C}_{i,j}(s-p+1) + \tilde{C}_{i,j^{*}}(s+h-p+1) = \epsilon_{j,s} \epsilon_{i,p} \langle \gamma, \beta \rangle_Q + \epsilon_{j^{*},s+h} \epsilon_{i,p} \langle \gamma, \beta \rangle_Q   \\
  &= \epsilon_{j,s} \epsilon_{i,p} \langle \gamma, \beta \rangle_Q - \epsilon_{j,s} \epsilon_{i,p} \langle \gamma, \beta \rangle_Q  = 0 . 
   \end{align*}
 where we used again Proposition~\ref{propositionFujitaOh}. 
 
  If  $I_{p,\beta} \cap \varphi^{-1}([1,2N]) := \{(j,s)\}$ then we distinguish two subcases. If $(i,p)$ and $(j,s)$ belong both to $\varphi^{-1}([1,N])$ (resp. both to $\varphi^{-1}([N+1,2N])$), then  $\epsilon_{j,s} = \epsilon_{i,p}$ and on the other hand the condition $s \geq p$ implies that there are (possibly trivial) morphisms but no extensions from the indecomposable object of dimension vector $\gamma$ to the one of dimension vector $\beta$ in the heart $\text{mod}\mathbb{C}Q$ (resp. $(\text{mod}\mathbb{C}Q)[-1]$) of $\mathcal{D}^{b}(\text{mod}\mathbb{C}Q)$. Hence $\langle \gamma , \beta \rangle_Q \geq 0$ and Theorem~\ref{thmHLCtilde} yields
 $$ n_t(\beta) = \tilde{C}_{i,j}(s-p+1) = \epsilon_{j,s} \epsilon_{i,p} \langle \gamma , \beta \rangle_Q  = \langle \gamma , \beta \rangle_Q \geq 0 . $$
 If on the contrary $(i,p)$ and $(j,s)$ do not belong both to $\varphi^{-1}([1,N])$ or $\varphi^{-1}([N+1,2N])$, then  $\epsilon_{j,s} = - \epsilon_{i,p}$ and on the other hand, in the heart containing $(i,p)$ the unique couple $(j',s')$ such that $\beta_{\varphi(j',s')} = \beta$ necessarily satisfies $s'<p$ (otherwise $I_{p,\beta}$ would be of cardinality $2$). Therefore there are no morphisms  from the indecomposable object of dimension vector $\gamma$ to the one of dimension vector $\beta$ in this heart, which implies $\langle \gamma , \beta \rangle_Q \leq 0$. Theorem~\ref{thmHLCtilde} yields
  $$ n_t(\beta) = \tilde{C}_{i,j}(s-p+1) = \epsilon_{j,s} \epsilon_{i,p} \langle \gamma , \beta \rangle_Q  = - \langle \gamma , \beta \rangle_Q \geq 0 . $$
  This concludes the proof of the Proposition. 
   \end{proof}
 
  \begin{remark}
 In the case of the fundamental modules $X_{i,\xi(i)} = L(Y_{i,\xi(i)})$, the formula of Lemma~\ref{DtildeinitialKRs} can also be rewritten explicitly from the quiver $Q$ as
 $$ \td_{\xi} \left( \iota(x_{\varphi(i,\xi(i))}) \right) =  \td_{\xi}(Y_{i, \xi(i)}) = \prod_{j \in B(i)} \frac{1}{\gamma_j}  $$
 where $B(i)$ denotes the set of indices $j$ such that there is a path from $j$ to $i$ in $Q$ (see Section~\ref{remindAR}). Indeed, one has 
 $$ \tilde{C}_{i,j}(s-\xi(i)+1) =
  \begin{cases}
   1 &\text{if $s=\xi(j)$ and $\xi(j)-\xi(i)=d(i,j)$ by Remark~\ref{rkCtilde},} \\   
   0 & \text{otherwise, by Lemma~\ref{littleLemmaforCtilde}.}
 \end{cases}  
 $$
 As already mentioned in the end of the proof of Lemma~\ref{DtildeinitialKRs}, $\xi(j)-\xi(i)=d(i,j)$ if and only if $j \in B(i)$, which yields the formula.  
  \end{remark}

   \begin{prop}[\textbf{Property $(B_{\wiQ})$}] \label{BforDtilde}
  Let $t \geq 1$ and let $(i,p) := \varphi^{-1}(t) \in \Ixi$. Then one has 
   $$ \td_{\xi} \left( \iota(x_t) \right) \td_{\xi}  \left( \iota(x_{t_{-}}) \right) =  \beta_t^{-1}   \prod_{ \substack{ r<t<r_{+} \\ i_r \sim i  }} \td_{\xi}  \left( \iota (x_r) \right) . $$
   \end{prop}
  
  Recall the notation $N_Q(k)$ from Section~\ref{remindAR}. We will need the following observation:
  
  \begin{lemma} \label{lemalternating}
 Let $t \geq 1$ and let $(i,p) := \varphi^{-1}(t) \in \Ixi$.  Then for any $r \geq 1$ such that $i_r \sim i$ one has 
$$ r < t < r_{+} \enspace \Leftrightarrow \enspace \varphi^{-1}(r) = (i_r,p+1) . $$
   \end{lemma}  
   
   \begin{proof}
   This is a consequence of the well-known fact that reduced expressions of $w_0$ adapted to orientations of Dynkin graphs are \textit{alternating} (see for instance \cite{Zelikson}). This means that each neighbour of the letter $i$ appears exactly once between two consecutive occurrences of $i$ in $\iQ$. Therefore this also holds for the infinite sequence $\wiQ$, as any finite subword of length $N$ of $\wiQ$ is still a reduced expression of $w_0$ adapted to some orientation of the Dynkin graph of $\mathfrak{g}$ (see Section~\ref{remindAR}). 
   Thus there are two possibilities: 
  if the first occurrence of $i_r$ appears before the  first occurrence of $i$, then there is an arrow from $i_r$ to $i$ in $Q$ and hence $\xi(i_r) = \xi(i)+1$; furthermore the $k$th occurrence of $i$ appears between the $k$th and the $k+1$th occurrences of $i_r$. In other words, $r < t < r_{+} \Leftrightarrow N_Q(r)=N_Q(t)$. Thus we get $(\xi(i_r) - \xi(i) +2N_Q(t) -1)/2 = N_Q(t) = N_Q(r)$.
   If on the other hand the first occurrence of $i_r$ appears after the  first occurrence of $i$, then there is an arrow from $i$ to $i_r$ in $Q$ and hence $\xi(i_r) = \xi(i)-1$;  furthermore the $k$th occurrence of $i$ appears between the $k-1$th and the $k$th occurrences of $i_r$. In other words, $r < t < r_{+} \Leftrightarrow N_Q(r)=N_Q(t)-1$. Thus we get $(\xi(i_r) - \xi(i) +2N_Q(t) -1)/2 = N_Q(t)-1= N_Q(r)$. Thus we have proved that 
 $$ (\text{$i_r \sim i$ and $r < t < r_{+}$}) \Leftrightarrow (\text{$i_r \sim i$ and  $N_Q(r) = (\xi(i_r) - \xi(i) +2N_Q(t) -1)/2$}) . $$
 Recalling the definition of $\varphi^{-1}$ from Section~\ref{remindAR}, this is equivalent to 
  $$ \varphi^{-1}(r) = \left( i_r, \xi(i_r) - 2N_Q(r) +2 \right) = \left( i_r, \xi(i) -2N_Q(t) + 3 \right) = \left( i_r, p+1 \right) . $$
  This finishes the proof of the Lemma. 
   \end{proof}
   
   We are now ready to prove Proposition~\ref{BforDtilde}.
  
   \begin{proof}[\textbf{Proof of Proposition~\ref{BforDtilde}}]
    Let us fix $(j,s) \in \Ixi$ and investigate the multiplicity of the positive root $\beta_{\varphi(j,s)}$ in both hand sides of this equality. 
    By Lemma~\ref{lemalternating}, one has 
    $$ \prod_{ \substack{ r<t<r_{+} \\ i_r \sim i  }} \td_{\xi}  \left( \iota (x_r) \right) = \prod_{k \sim i} \td_{\xi}  \left( \iota (x_{\varphi(k,p+1)}) \right) . $$
 Therefore using Lemma~\ref{DtildeinitialKRs}, the multiplicity of $\beta_{\varphi(j,s)}$ in the right hand side can be written as
 $$  - \delta_{i,j} \delta_{p,s} -  \sum_{k \sim i} \tilde{C}_{k,j}(s-p) . $$
  Using the relations~\eqref{relationsforCtilde} with $m=s-p$, this is equal to 
$$- \left( \tilde{C}_{i,j}(s-p+1) + \tilde{C}_{i,j}(s-p-1) \right) . $$
By Lemma~\ref{DtildeinitialKRs}, this is exactly the multiplicity of $\beta_{\varphi(j,s)}$ in the product 
\begin{align*}
\td_Q \left( \iota ( x_{\varphi(i,p)}) \right)  \td_Q \left( \iota ( x_{\varphi(i,p+2)}) \right),
\end{align*}
which is equal to $\td_Q \left( \iota (x_t) \right)  \td_Q \left( \iota (x_{t_{-}}) \right)$ by~\eqref{defvarphi}. 
  Thus the multiplicities of $\beta_{\varphi(j,s)}$ on both hand sides coincide for every $(j,s) \in \Ixi$ which proves the Proposition. 
   \end{proof}
  
 \begin{remark} \label{DtildeAinverse}
Recall from~\eqref{defnAip} the variables $A_{j,s}, j \in I, s \in \mathbb{Z}$.  It is straightforward to check either from Proposition~\ref{BforDtilde} or directly from the definition of $\td_{\xi}$ that for every $(i,p) \in \Ixi$ one has 
 $$ \td_{\xi} \left( A_{i,p-1}^{-1} \right) = \frac{\beta_{\varphi(i,p-2)}}{\beta_{\varphi(i,p)}} = \frac{\beta_{t_{+}}}{\beta_t} $$
 where $t := \varphi(i,p)$. This can be viewed as a generalization of {{\cite[Remark 6.4]{Casbi3}}} as it is known from the works of Hernandez-Leclerc \cite{HLJEMS} that the $A_{i,p-1}^{-1}, (i,p) \in \Ixi$ are exactly the images under $\iota$ of Fomin-Zelevinsky's variables $\yjh$ (see \cite{FZ4}) for the seed $\s^{\wiQ}$ (up to the convention used for the definition of $\yjh$). 
  \end{remark}  
  
 Recall from Section~\ref{remindDbarflagminors} that for any rational fraction $Y$ and any positive root $\beta$, we denote by $[\beta;Y]$ the (algebraic) multiplicity of $\beta$ in $Y$. 
 
  \begin{prop}[\textbf{Property $(C_{\wiQ})$}] \label{CforDtilde}
  Let $(i,p) \in \Ixi$. Then for any $\beta \in \Phi_{+}$, one has $ \mid [\beta ; \td_{\xi} (Y_{i,p})] \mid \leq 1$. 
   \end{prop}  
   
    \begin{proof}
    Let us fix $\beta \in \Phi_{+}$ and let $\gamma := \beta_{\varphi(i,p)}$. By~\eqref{defDtildexi} one has 
 $$ -  [\beta ; \td_{\xi} (Y_{i,p})]  = \sum_{(j,s) \in I_{p,\beta}} \mathcal{N}(i,p;j,s), \qquad \mathcal{N}(i,p;j,s) := \tilde{C}_{i,j}(s-p+1) - \tilde{C}_{i,j}(s-p-1) $$
 where we use the notation $I_{p,\beta}$ from the proof of Proposition~\ref{AforDtilde}.
 
  Let $s$ be the smallest integer such that there exists $j \in I$ with $(j,s) \in I_{p,\beta}$. It follows from Proposition~\ref{propositionFujitaOh} that $I_{p,\beta}$ is either of the form $\{ (j,s), (j^{*},s+h), (j,s+2h), (j^{*},s+3h), \ldots , (j,s+2mh) \}$ for some $m \geq 1$ if $\sharp I_{p,\beta}$ is odd, or of the form $\{ (j,s), (j^{*},s+h), (j,s+2h), (j^{*},s+3h), \ldots , (j,s+2mh), (j^{*},s+(2m+1)h) \}$ for some $m \geq 1$ if $\sharp I_{p,\beta}$ is even. Then Corollary~\ref{cormultiplicityDtilde} yields
 \begin{align*}
  \mathcal{N}(i,p;j^{*},s+(2k-1)h) + \mathcal{N}(i,p;j,s+2kh) &= \epsilon_{i,p} \epsilon_{j^{*},s+(2k-1)h} (\beta, \gamma)+ \epsilon_{i,p} \epsilon_{j,s+2kh} (\beta, \gamma)  \\
   &= 0 
   \end{align*}
 as $\epsilon_{j^{*},s+(2k-1)h} = - \epsilon_{j,s+2kh}$ for any $k \geq 1$. Thus we obtain 
  $$ -  [\beta ; \td_{\xi} (Y_{i,p})] = 
  \begin{cases}
    \mathcal{N}(i,p;j,s) & \text{if $\sharp I_{p, \beta}$ is odd,} \\
     \mathcal{N}(i,p;j,s) - \epsilon_{i,p} \epsilon_{j,s} (\beta , \gamma)  & \text{if $\sharp I_{p, \beta}$ is even.}
  \end{cases}
  $$
  where in the second case we used Corollary~\ref{cormultiplicityDtilde} for $(i,p),(j^{*},s+(2m+1)h)$ and the fact that $\epsilon_{j^{*},s+(2m+1)h} = - \epsilon_{j,s}$. 
  
   From this together with Corollary~\ref{cormultiplicityDtilde}, we obtain that if $s>p$ then $[\beta ; \td_{\xi} (Y_{i,p})]$ is equal to $(\beta , \gamma)$ up to some sign if $\sharp I_{p,\beta}$ is odd, and to $0$ if $\sharp I_{p,\beta}$ is even. Note that in this case one has $\beta \neq \gamma$ (if $\beta= \gamma$ then $(i,p) \in I_{p,\beta}$ and thus $s=p$ by minimality). Then it is elementary to check that for any simply-laced type Lie algebra $\mathfrak{g}$, the Cartan pairing of any two distinct positive roots of $\mathfrak{g}$ is always equal to $-1,0$ or $1$. This can be deduced for instance from Lemmas~\ref{CartanpairingA} and~\ref{CartanpairingD} below respectively for types $A_n$ and $D_n$, and can be checked directly for the types $E_6,E_7$ and $E_8$. 
   
 If on the other hand $s=p$, then Corollary~\ref{cormultiplicityDtilde} implies $\mathcal{N}(i,p;j,s) = \delta_{i,j}$ and moreover by standard Auslander-Reiten theory one has $(\beta, \gamma)=0$ if $i \neq j$ and $(\beta , \gamma)=(\beta,\beta)=2$ if $i=j$. In other words, one has $(\beta , \gamma) = 2 \delta_{i,j}$. Therefore, $-  [\beta ; \td_{\xi} (Y_{i,p})]$ is equal either to $\delta_{i,j}$ if $\sharp I_{p, \beta}$ is odd, or to $- \delta_{i,j}$ if $\sharp I_{p, \beta}$ is even. This concludes the proof of the Proposition. 
    \end{proof}
  
  \section{Evaluation of $\barD$ on the dual root vectors}
  \label{Dbarcuspidals}

 This section is devoted to the computation of the values taken by Baumann-Kamnitzer-Knutson's morphism $\barD$ on the dual root vectors of $\CN$ with respect to the convex ordering on $\Phi_{+}$ corresponding to $\mathbf{i}_{Q_0}$ where $Q_0$ is the orientation of the Dynkin graph of $\mathfrak{g}$ shown in Figure~\ref{figQ0}. We provide explicit formulas in types $A_n , n \geq 1$ and $D_n , n \geq 4$. The computations in types $E_r , r=6,7,8$ are performed separately using a computer software. 

  \subsection{Type $A_n, n \geq 1$}
  \label{DbarcuspidalA}
 
  We consider the case where $\mathfrak{g}$ is of type $A_n , n \geq 1$.  For every $1 \leq i \leq j \leq n$, we set $\alpha_{i,j} := \alpha_i + \cdots + \alpha_j$, and we have $\Phi_{+} = \{ \alpha_{i,j} , 1 \leq i \leq j  \leq n \}$. 
  We choose the following reduced expression of $w_0$:
  $$ (1,2, \ldots , n, 1, 2, \ldots , n-1, \ldots , 1, 2, 1) $$
   which is adapted to the orientation $Q_0$ (the so-called monotonic orientation) of the Dynkin graph of $\mathfrak{g}$ (see Figure~\ref{figQ0}). 
The corresponding convex ordering on $\Phi_{+}$ is given by 
  $$ \alpha_{i,j} < \alpha_{k,l} \Leftrightarrow i<k \enspace \text{or} \enspace \text{$i=k$ and $j<l$} . $$
  Thus, it coincides with the Lyndon ordering arising from the choice of the natural order $1<2 < \cdots <n$ on the index set of simple roots. Consequently, the cuspidal representations are explicitly constructed in {{\cite[Section 8.4]{KR}}}, namely one has $\mathbf{j}_{\alpha_{i,j}} = (i,i+1, \ldots, j)$ and $S_{\alpha_{i,j}}$ is the one-dimensional vector space generated by a single vector on which all the generators of the quiver Hecke algebras $R(\beta), \beta \in \Gamma_{+}$ act by zero, except the idempotent $e(\mathbf{j}_{\alpha_{i,j}})$. Applying Equation~\eqref{eqnDbarKLR} we obtain 
  \begin{equation} \label{DbarofSalpha}
   \barD([S_{\alpha_{i,j}}]) = \frac{1}{\alpha_i(\alpha_i + \alpha_{i+1}) \cdots (\alpha_i + \cdots + \alpha_j)} = \prod_{i \leq k \leq j} \frac{1}{\alpha_{i,k}} . 
   \end{equation}
Alternatively one can note that the word $\mathbf{j}_{\alpha_{i,j}}$ is dominant minuscule, hence $S_{\alpha_{i,j}}$ is strongly homogeneous (see Definition~\ref{defstronghomog}) and we can conclude using Proposition~\ref{lemDbarKLR}.

  \subsection{Type $D_n , n \geq 4$}
  \label{DbarcuspidalD}
  
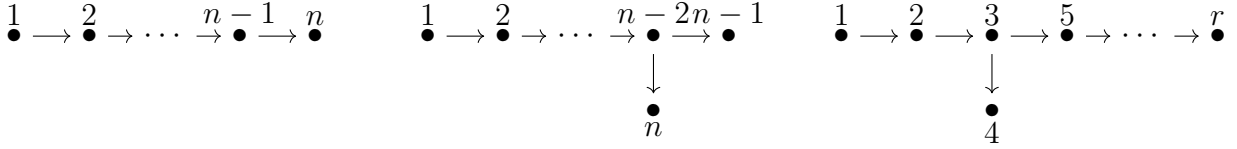
\begin{figure}
 
   \begin{tikzpicture}
    
     \node (a1) at (-5.5,0) {$\bullet$};
     \node (a2) at (-4.5,0) {$\bullet$};
     \node (a3) at (-3.5,0) {$\cdots$};
     \node (a4) at (-2.5,0) {$\bullet$};
     \node (a5) at (-1.5,0) {$\bullet$};
   
   \draw (a1) node[above] {$1$};
   \draw (a2) node[above] {$2$};
   \draw (a4) node[above] {$n-1$};
   \draw (a5) node[above] {$n$};
     
\draw[->] (a1)--(a2);
\draw[->] (a2)--(a3);
\draw[->] (a3)--(a4);
\draw[->] (a4)--(a5);
    
     \node (d1) at (0,0) {$\bullet$};
     \node (d2) at (1,0) {$\bullet$};
     \node (d3) at (2,0) {$\cdots$};
     \node (d4) at (3,0) {$\bullet$};
     \node (d5) at (4,0) {$\bullet$};
     \node (d6) at (3,-1) {$\bullet$};
     
   \draw (d1) node[above] {$1$};
   \draw (d2) node[above] {$2$};
   \draw (d4) node[above] {$n-2$};
   \draw (d5) node[above] {$n-1$};
   \draw (d6) node[below] {$n$};
     
\draw[->] (d1)--(d2);
\draw[->] (d2)--(d3);
\draw[->] (d3)--(d4);
\draw[->] (d4)--(d5);
\draw[->] (d4)--(d6);

     \node (e1) at (5.5,0) {$\bullet$};
     \node (e2) at (6.5,0) {$\bullet$};
     \node (e3) at (7.5,0) {$\bullet$};
     \node (e4) at (7.5,-1) {$\bullet$};
     \node (e5) at (8.5,0) {$\bullet$};
     \node (e6) at (9.5,0) {$\cdots$};
     \node (e7) at (10.5,0) {$\bullet$};
     
   \draw (e1) node[above] {$1$};
   \draw (e2) node[above] {$2$};
   \draw (e3) node[above] {$3$};
   \draw (e4) node[below] {$4$};
    \draw (e5) node[above] {$5$};
   \draw (e7) node[above] {$r$};
     
\draw[->] (e1)--(e2);
\draw[->] (e2)--(e3);
\draw[->] (e3)--(e4);
\draw[->] (e3)--(e5);
\draw[->] (e5)--(e6);
\draw[->] (e6)--(e7);

  \end{tikzpicture}
  
   \caption{The orientation $Q_0$ for each simply-laced type.}
   
    \label{figQ0}
 
 \end{figure}

  We now focus on the case where $\mathfrak{g}$ is of type $D_n, n \geq 4$.
  We will use the following notations: for any $(p,q) \in \{1, \ldots , n-1\}^2$, we set:
 $$ \theta_{p,q} := \alpha_{\min(p,q)} + \cdots + \alpha_{\max(p,q)-1} + 2(\alpha_{\max(p,q)} + \cdots + \alpha_{n-2}) + \alpha_{n-1} + \alpha_n . $$
 and for every $1 \leq p \leq q \leq n$ we set:
 $$ \alpha_{p,q} := 
 \begin{cases} 
  \alpha_p + \cdots + \alpha_q & \text{if $q \leq n-1$,} \\
  \alpha_p + \cdots + \alpha_{n-2} + \alpha_n & \text{if $q=n$.}
 \end{cases}
$$
In particular we have $\alpha_{n-1,n} = \alpha_{n,n} := \alpha_n$. Also note that $\theta_{p,q}$ is a positive root if and only if $p \neq q$. We have 
$$\Phi_{+} = \{ \theta_{p,q}, 1 \leq p < q \leq n-1 \} \sqcup \{ \alpha_{p,q}, 1 \leq p \leq q \leq n \}. $$
We will need to consider the automorphism $\sigma$ of the Dynkin diagram of $\mathfrak{g}$ defined by $\sigma(i) = i$ if $1 \leq i \leq n-2$, $\sigma(n-1)=n$ and $\sigma(n) = n-1$. 
 We choose the following reduced expression of $w_0$:
  $$ (1,2, \ldots , n)^{n-1} $$
   which is adapted to the orientation $Q_0$ of the Dynkin graph of $\mathfrak{g}$ shown in Figure~\ref{figQ0}. 
The corresponding convex ordering on $\Phi_{+}$ is given by 
 \begin{align*}
   \alpha_{i,j} < \alpha_{k,l}  \enspace & \Leftrightarrow \enspace (i<k) \enspace  \text{or} \enspace \left( \text{$i=k$ and ($j < \min(l,n-1)$ or $j=\sigma^{i}(n)$ and $l= \sigma^{i}(n-1)$})  \right), \\
   \alpha_{i,j} < \theta_{p,q}  \enspace & \Leftrightarrow \enspace  i<q \quad  \text{or} \quad \text{$i \leq q$ and $j \leq n-2$,}  \\
   \theta_{p,q} < \theta_{r,s} \enspace & \Leftrightarrow \enspace q<s \quad  \text{or} \quad \text{$q=s$ and $p<r$.}
 \end{align*}
 The cuspidal representations  corresponding to the positive roots $\alpha_{i,j}$  are given in the same way as in type $A_n$, i.e. for every $1 \leq i \leq j \leq n$,  one has $\mathbf{j}_{\alpha_{i,j}} = (i,i+1, \ldots, j)$ and the cuspidal representation $S_{\alpha_{i,j}}$ has a unique non trivial one-dimensional  weight space of weight $\mathbf{j}_{\alpha_{i,j}}$ (if $j=n$ then the word $\mathbf{j}_{\alpha_{i,j}} = (i,i+1, \ldots, j)$ is understood as $(i,i+1, \ldots , n-2,n)$). Thus the conclusion is the same as in the previous paragraph. 
 
  Let us now focus on the cuspidal representations associated to the positive roots $\theta_{p,q} , 1 \leq p < q \leq n-1$. These cuspidal modules turn out to be not homogeneous. It is therefore not possible to compute $\barD([S_{\theta_{p,q}}])$ directly. This is why we use Brundan-Kleshchev-McNamara's distinguished short exact sequences from Theorem~\ref{thmBKMshortexactseq}. 
  
   \begin{lemma} \label{minimalpairfortheta}
   Fix $p,q$ such that $1 \leq p < q \leq n-1$. Then the couple of positive roots $(\alpha_{p, \sigma^p(n-1)} , \alpha_{q, \sigma^p(n)})$ is a minimal pair for $\theta_{p,q}$ with respect to the chosen order $<$ on $\Phi_{+}$. 
   \end{lemma}

 \begin{proof}
 Assume there exists $(\gamma, \delta) \in \Phi_{+}^2$ such that $\gamma + \delta = \theta_{p,q}$ and $\alpha_{p, \sigma^p(n-1)} < \gamma < \theta_{p,q} < \delta <  \alpha_{q, \sigma^p(n)}$. If $\gamma$ is of the form $\theta_{r,s}$, then one has either  $s<q$ or $s=q$ and $r<p$. In both cases, $\theta_{p,q} - \gamma \notin \Gamma_{+}$ which is a contradiction. Thus $\gamma$ is of the form $\alpha_{i,j}$. Then $\alpha_{p, \sigma^p(n-1)} < \gamma$ implies $i>p$. Now if $\delta$ was of the form $\theta_{r,s}$ then $\theta_{p,q} < \delta <  \alpha_{q, \sigma^p(n)}$ implies $r>p$. If on the other hand $\delta = \alpha_{u,v}$ then $\theta_{p,q} < \delta$ implies $u \geq q >p$. In both cases, $\alpha_p$ appears neither in $\gamma$ nor in $\delta$ which contradicts $\gamma + \delta = \theta_{p,q}$ as $\alpha_p$ appears in $\theta_{p,q}$.  
 \end{proof}
  
  \begin{prop} \label{DbarofStheta}
  Fix $p,q$ such that $1 \leq p < q \leq n-1$. Then one has 
  $$ \barD([S_{\theta_{p,q}}]) = \frac{\theta_{p,p}}{\alpha_{q,q} \cdots \alpha_{q,n-2}  \alpha_{p,p} \cdots \alpha_{p,n} \theta_{p,q}}  . $$
   \end{prop}
   
    \begin{proof}
   We assume $p$ is even, the other case being identical. 
    By Lemma~\ref{minimalpairfortheta}, we have a minimal pair for $\theta_{p,q}$ given by $(\alpha_{p,n-1} , \alpha_{q,n})$. By Theorem~\ref{thmBKMshortexactseq}, there is a short exact sequence in $R-mod$:
    \begin{align} \label{SEStheta}
     0 \longrightarrow S_{\theta_{p,q}} \longrightarrow S_{\alpha_{q,n}} \circ S_{\alpha_{p,n-1}}  \longrightarrow L(\mathbf{c}_{p,q}) \longrightarrow 0  
       \end{align}
     where $\mathbf{c}_{p,q}$ is the $N$-tuple of integers whose entries are all zero, except the two entries corresponding to the positive roots $\alpha_{p,n-1}$ and $\alpha_{q,n}$ which are equal to $1$. Hence Theorem~\ref{thmrootpartitions} implies that
     $$ \dim_{\mathbb{C}} \left( e({\bf j}_{\alpha_{q,n}} {\bf j}_{\alpha_{p,n-1}}) L(\mathbf{c}_{p,q}) \right) = 1 . $$
     As recalled above, we have that ${\bf j}_{\alpha_{q,n}} = (q,q+1, \ldots , n-2, n)$ and ${\bf j}_{\alpha_{p,n-1}} = (p,p+1, \ldots , n-2,n-1)$. It is now straightforward to check from Definition~\ref{defminuscule} that the word ${\bf j}_{\alpha_{q,n}} {\bf j}_{\alpha_{p,n-1}}$ is dominant minuscule (i.e. it is a reduced expression of a dominant minuscule element of $W$). Consequently, Theorem~\ref{thmKRhom} implies that $L(\mathbf{c}_{p,q})$ is strongly homogeneous, and hence Proposition~\ref{lemDbarKLR} yields
     $$ \barD([L(\mathbf{c}_{p,q})]) = \prod_{\beta \in \Phi_{+}^{w_{p,q}}} \frac{1}{\beta} $$
     where $w_{p,q} := s_q s_{q+1} \cdots s_{n-2} s_n s_p s_{p+1} \cdots s_{n-2} s_{n-1}$. 
     For every $q \leq s \leq n-2$, we have $s_q \cdots s_{s-1} (\alpha_s) = \alpha_{q,s}$. We also have that $s_q \cdots s_{n-2} (\alpha_n) = \alpha_{q,n}$. For every $p \leq r < n-2$, we have 
  $$ s_q \cdots s_{n-2} s_n s_p \cdots s_{r-1} (\alpha_r) = s_q \cdots s_{n-2} (\alpha_{p,r}) = 
   \begin{cases}  
    \alpha_{p,r} & \text{if $r<q-1$,} \\
   \alpha_{p,r+1}  & \text{if $r \geq q-1$.}
   \end{cases}
 $$
  Finally the last two positive roots are  $s_q \cdots s_{n-2} s_n s_p \cdots s_{n-3} (\alpha_{n-2}) = s_q \cdots s_{n-2} (\alpha_{p,n}) = \alpha_{p,n}$ and   $s_q \cdots s_{n-2} s_n s_p \cdots s_{n-2} (\alpha_{n-1})  = s_q \cdots s_{n-2} s_n (\alpha_{p,n-1}) =  s_q \cdots s_{n-2} (\theta_{p,n-1}) = \theta_{p,q}$. Thus we have 
     $$ \barD([L(\mathbf{c}_{p,q})]) = \frac{1}{\alpha_{p,p} \cdots \alpha_{p,q-2} \alpha_{p,q} \cdots  \alpha_{p,n} \alpha_{q,q} \cdots \alpha_{q, n-2} \alpha_{q,n} \theta_{p,q}} . $$
     As $\barD$ is an algebra morphism on $\CN \simeq K_0(R-mod)$, the short exact sequence~\eqref{SEStheta} yields
    \begin{align*}
     \barD([S_{\theta_{p,q}}]) &= \barD([S_{\alpha_{p,n-1}}]) \barD([S_{\alpha_{q,n}}]) - \barD([L(\mathbf{c}_{p,q})])  \\
      &= \frac{1}{\alpha_{p,p} \cdots \alpha_{p,n-1} \alpha_{q,q} \cdots \alpha_{q,n-2} \alpha_{q,n}} - \frac{1}{\alpha_{p,p} \cdots \alpha_{p,q-2} \alpha_{p,q} \cdots  \alpha_{p,n} \alpha_{q,q} \cdots \alpha_{q, n-2} \alpha_{q,n} \theta_{p,q}} \\
      &= \frac{\theta_{p,q} \alpha_{p,n} - \alpha_{p,q-1} \alpha_{p,n-1}}{\alpha_{p,p} \cdots \alpha_{p,n-1} \alpha_{p,n} \alpha_{q,q} \cdots \alpha_{q,n-2} \alpha_{q,n} \theta_{p,q}} = \frac{\theta_{p,q} \alpha_{q,n} + \alpha_{p,q-1} (\theta_{p,q} - \alpha_{p,n-1})}{\alpha_{p,p} \cdots \alpha_{p,n-1} \alpha_{p,n} \alpha_{q,q} \cdots \alpha_{q,n-2} \alpha_{q,n} \theta_{p,q}} \\
      &=  \frac{\theta_{p,q} \alpha_{q,n} + \alpha_{p,q-1} \alpha_{q,n}}{\alpha_{p,p} \cdots \alpha_{p,n-1} \alpha_{p,n} \alpha_{q,q} \cdots \alpha_{q,n-2} \alpha_{q,n} \theta_{p,q}} =   \frac{\theta_{p,p}}{\alpha_{p,p} \cdots \alpha_{p,n-1} \alpha_{p,n} \alpha_{q,q} \cdots \alpha_{q,n-2} \theta_{p,q}} 
        \end{align*}
        which concludes the proof. 
    \end{proof}

      \subsection{Types $E_6,E_7,E_8$}
   \label{DbarcuspidalE}

 When $\mathfrak{g}$ is a simple Lie algebra of type $E_r , r =6,7,8$ we use a computer software to compute the values of $\barD$ on the dual root vectors of $\CN$. This is based on Brundan-Kleshchev-McNamara's algorithm {{\cite[Theorem 4.2]{BKM}}} (see also \cite{Leclerc}) which allows to recursively determine the graded characters of the cuspidal representations  with respect to the convex ordering on $\Phi_{+}$ corresponding to $\mathbf{i}_{Q_0}$. Then we can simply forget the grading and compute the values of $\barD$ using the equality~\eqref{eqnDbarKLR}.

   On the contrary to the types $A_n$ and $D_n$, the formulas we get cannot be written in a simple way (especially the numerators cannot be factored as products of linear terms). We refer to Section~\ref{DtildeKRE} for more comments about this. 
    
\section{Evaluation of $\td_{Q_0}$ on the classes of Kirillov-Reshetikhin-modules}
\label{DtildeKRmonotonic}

In this section, we provide explicit formulas for the evaluation of $\td_{Q_0}$ on the classes of Kirillov-Reshetikhin-modules of $\mathcal{C}_{Q_0}$ when $\mathfrak{g}$ is of type $A_n, n \geq 1$ or $D_n, n \geq 4$ and $Q_0$ is the orientation of the Dynkin displayed in Figure~\ref{figQ0}. This yields in particular explicit formulas for the evaluation of $\td_{Q_0}$ on the classes of fundamental representations, which are known from the work of Hernandez-Leclerc \cite{HL15} to categorify the dual root vectors of $\CN$. 
For the classical types, we obtain formulas that could be written in a uniform way (the formulas for the type $A_n , n \geq 1$ are special cases contained in the ones for the types $D_n$). However, for the sake of readability, we prefer dealing with the two cases in different subsections. 
We discuss the types $E_r , r=6,7,8$ separately in Section~\ref{DtildeKRE}.

  Throughout this section, we will simply write $\td$ (resp. $\tau$) for $\td_{Q_0}$ (resp. $\tau_{Q_0}$) as there will be no ambiguity. 
  The techniques used in this Section can be naturally extended and applied to classes of Kirillov-Reshetikhin modules not belonging to $\CQ$. Therefore, one could obtain formulas for the evaluation of $\td$ on a large family of cluster variables in cluster algebras strictly larger than $\CN$, such as the algebra $\bAQ$ we introduce in Section~\ref{finalsection} or even the whole cluster algebra $K_0(\Cxi)$ itself.


\subsection{Type $A_n, n \geq 1$}
\label{DtildeKRA}


 We have $n_{Q_0}(i) = n-i+1$ for each $i \in I$ and $I_{Q_0} = \{(i, \xi(i)-2r+2) , i \in I , 1 \leq r \leq n-i+1 \}$.
For  each $(i,s) \in I_{Q_0}$ with $r := \frac{\xi_i-s+2}{2}$ and for each $k \in \{1, \ldots, r \}$,  we set   
\begin{equation} \label{eq:formula for tilde D of KR modules of type An}
    D_{i,s}^{(k)} := \prod_{ \substack{ r-k+1 \leq  p  \leq r  \\ r \leq q \leq r+i-1 } } \frac{1}{\alpha_{p,q}} \quad \text{and} \quad d_{i,s}^{(k)} := \frac{D_{i,s}^{(k)}}{D_{i,s+2}^{(k-1)}}  
    \end{equation}
  where $D_{i,s}^{(0)} := 1$ for any $i,s$. We also denote $D_{i,s} := D_{i,s}^{(r)}$ and $d_{i,s} := d_{i,s}^{(r)}$. Note that one has $d_{i,s} = D_{i,s} / D_{i,s+2}$. 
 The aim is to prove that the value of $\td$ on the class of the Kirillov-Reshetikhin module $X_{i,s}^{(k)}$ is given by $D_{i,s}^{(k)}$ (Theorem~\ref{thmKRmonotonicA}).   We first prove this in the case $r=k$ (Corollary~\ref{corDtildeinitialtypeA}). By Theorem~\ref{thmHL15}, the corresponding Kirillov-Reshetikhin modules $X_{i,s}$ categorify the cluster variables of the standard seed $\s^{\mathbf{i}_{Q_0}}$. We then prove the general formula using the $T$-systems from Section~\ref{remindTsystems}, which correspond to certain cluster mutations in $\CN$. 
  We begin with the following elementary lemma. 
 
  \begin{lemma} \label{CartanpairingA}
  For every $1 \leq i \leq p \leq n$ and $1 \leq j \leq q \leq n$ one has  
$$ (\alpha_{k,i} , \alpha_{l,j}) = \delta_{k,l} + \delta_{i,j} - \delta_{k-1,j} - \delta_{l-1,i} . $$
  \end{lemma}
  
   \begin{prop} \label{lem:tilde D of Yit in type A}
   For every $(i,s) \in I_{Q_0}$, one has $\td(Y_{i,s}) = d_{i,s}$. 
   \end{prop}
   
    \begin{proof}
    Let us fix $(i,s) \in I_{Q_0}$ and set $r := \frac{\xi_i-s+2}{2}$ as above. For this choice of orientation, we have $n_{Q_0}(i)=n-i+1$ for each $i \in I$ and $\xi(j)-\xi(i)=i-j$ for every $i,j \in I$. Let $(j,t) \in \Ixi$  and assume $(j,t) \notin I_{Q_0}$. Thus $t \leq \xi(j)-2n_{Q_0}(j)$ by~\eqref{defIQ}. If $t \geq s$ then $\xi(i)-2n_{Q_0}(i)<s \leq t \leq \xi(j)-2n_{Q_0}(j)$ and in particular $t-s < \xi(j) - \xi(i) -2(n_{Q_0}(j)-n_{Q_0}(i)) = i-j-2(i-j)=j-i=d(i,j)$. Hence Lemma~\ref{littleLemmaforCtilde} implies that $\tilde{C}_{i,j}(t-s+1) = \tilde{C}_{i,j}(t-s-1) = 0$. This is also the case if $t<s$ as $\tilde{C}_{i,j}(m)=0$ if $m \leq 0$. Therefore, we can rewrite~\eqref{defDtildexi} as 
    $$ \td_{Q_0}(Y_{i,s}) = \prod_{(j,t) \in I_{Q_0}} \beta_{\varphi(j,s)}^{\tilde{C}_{i,j}(t-s-1) - \tilde{C}_{i,j}(t-s+1)} = \prod_{(j,t) \in I_{Q_0}} \left( \tau^{(\xi(j)-t)/2}(\gamma_j)   \right)^{\tilde{C}_{i,j}(t-s-1) - \tilde{C}_{i,j}(t-s+1)} . $$
    Moreover, one has $\tau^{r-1}(\gamma_i) = \alpha_{r,r+i-1}$ and  $\tau^{l-1}(\gamma_j) = \alpha_{l,l+j-1}$ for every $j \in \{1, \ldots , n\}$ and $l \in \{1, \ldots , n-j\}$. Therefore, Corollary~\ref{cormultiplicityDtilde}  together with Lemma~\ref{CartanpairingA} yield
   \begin{align*}
    \td(Y_{i,s}) &= \frac{1}{\alpha_{r,r+i-1}}  \prod_{2(l-r) < i-j} \alpha_{l,l+j-1}^{- \delta_{r,l} - \delta_{r+i-1,l+j-1} + \delta_{r-1,l+j-1} + \delta_{r+i,l}}  \\
      &= \frac{1}{\alpha_{r,r+i-1}}  \prod_{1 \leq j < i} \frac{1}{\alpha_{r,r+j-1}} \prod_{i<j \leq r+i-1} \frac{1}{\alpha_{r+i-j,r+i-1}} \prod_{1 \leq j \leq r-1} \alpha_{r-j,r-1}  \\
      &= \frac{\prod_{1 \leq p \leq r-1} \alpha_{p,r-1}}{\left( \prod_{1 \leq p \leq r-1} \alpha_{p,r+i-1} \right) \left( \prod_{r \leq q \leq r+i-2} \alpha_{r,q} \right)} \\
      &= \left( \prod_{ \substack{ 1 \leq p \ \leq r-1 \\ r-1 \leq q \leq i+r-2] } } \alpha_{p,q} \right) \left(  \prod_{ \substack{ 1 \leq p \leq r \\ r \leq q \leq i+r-1 } } \alpha_{p,q} \right)^{-1} = \frac{D_{i,r,\xi_i-2r+2}}{D_{i,r-1,\xi_i-2r+4}} = d_{i,r,\xi_i-2r+2}.
     \end{align*}
    \end{proof}
    
     \begin{corollary} \label{corDtildeinitialtypeA}
  For every $(i,s) \in I_{Q_0}$, one has $\td(\tchi_q(X_{i,s})) = D_{i,s}$. 
      \end{corollary}
      
       \begin{proof}
    Recall from  Section~\ref{remindTsystems} that the truncated $q$-character of the Kirillov-Reshetikhin module $X_{i,s}$ has only term, namely the dominant monomial $Y_{i,s} Y_{i,s+2} \cdots Y_{i,\xi(i)}$. Therefore, Proposition~\ref{lem:tilde D of Yit in type A} implies
       $$ \td(\tchi_q(X_{i,s})) = \td(Y_{i,s}) \td(Y_{i,s+2}) \cdots \td(Y_{i,\xi(i)}) = d_{i,s} d_{i,s+2} \cdots d_{i,\xi(i)} = D_{i,s} . $$
       \end{proof}

\begin{thm} \label{thmKRmonotonicA}
 For every $(i,s) \in I_{Q_0}$ and every $k \in  \{ 1, \ldots, r \}$, we have $\td( \tchi_q(X_{i,s}^{(k)}) ) = D_{i,s}^{(k)}$. 
\end{thm}

\begin{proof}

By Corollary~\ref{corDtildeinitialtypeA}, the statement is true if $r=k$ i.e. for the Kirillov-Reshetikhin modules corresponding to the cluster variables of the (standard) seed $\s^{\mathbf{i}_{Q_0}}$. 
As $\td$ is an algebra morphism, the rational fractions $\td(X_{i,s}^{(k)})$ satisfy the $T$-systems~\eqref{eq:T-system relations}. Since the solution to the T-system with a given initial condition is unique, it suffices to show that $D_{i,s}^{(k)}$ satisfies the recursive equations:
$$ D_{i,s}^{(k)} D_{i,s-2}^{(k)} = D_{i,s-2}^{(k+1)} D_{i,s}^{(k-1)} + \prod_{j \sim i} D_{j,s-1}^{(k)}. $$
Recall that $r = \frac{\xi_i-s+2}{2}$. By definition of the $D_{i,s}^{(k)}$, we have
\begin{align*}
& D_{i,s}^{(k)} D_{i,s-2}^{(k)} - D_{i,s-2}^{(k+1)} D_{i,s}^{(k-1)} \\
& = \prod_{ \substack{ p \in [r-k+1, r] \\ q \in [r,r+i-1] }} \frac{1}{\alpha_{p,q}}  \prod_{ \substack{ p \in [r-k+2,r+1] \\ q \in [r+1, r+i]}} \frac{1}{\alpha_{p,q}} \quad  -  \enspace \prod_{ \substack{ p \in [r-k+1,r+1] \\ q \in [r+1,r+i] }} \frac{1}{\alpha_{p,q}}  \prod_{ \substack{ p \in [r-k+2,r] \\ q \in [r,r+i-1] }} \frac{1}{\alpha_{p,q}} \\
& = \prod_{ \substack{ p \in [r-k+2, r] \\ q \in [r,r+i-1] } } \frac{1}{\alpha_{p,q}}  \prod_{  \substack{ p \in [r-k+2, r+1] \\ q \in [r+1,r+i] } } \frac{1}{\alpha_{p,q}} \left( \frac{1}{\prod_{q \in [r, r+i-1]} \alpha_{r-k+1, q}} - \frac{1}{\prod_{q \in [r+1, r+i]} \alpha_{r-k+1, q}} \right) \\
& = \prod_{ \substack{ p \in [r-k+2, r] \\ q \in [r,r+i-1] } } \frac{1}{\alpha_{p,q}}  \prod_{  \substack{ p \in [r-k+2, r+1] \\ q \in [r+1,r+i] } } \frac{1}{\alpha_{p,q}} \left( \frac{ \alpha_{r+1,r+i} }{ \alpha_{r-k+1, r+i}  \prod_{q \in [r, r+i-1]} \alpha_{r-k+1, q} } \right) \\
& = \prod_{ \substack{ p \in [r-k+1, r] \\ q \in [r,r+i] } } \frac{1}{ \alpha_{p,q} }  \prod_{\substack{ p \in [r-k+2, r+1] \\ q \in [r+1,r+i-1] } } \frac{1}{\alpha_{p,q}}  =  D_{i+1,s-1}^{(k)} D_{i-1,s-1}^{(k)}.
\end{align*}

\end{proof}

\subsection{Type $D_n, n \geq 4$}
\label{DtildeKRD}

 We have $n_{Q_0}(i)=n-1$ for each $i \in I$ and $I_{Q_0} = \{(i, \xi(i)-2r+2) , i \in I , 1 \leq r \leq n-1\}$. For each $(i,s) \in I_{Q_0}$ with $r := \frac{\xi_i-s+2}{2}$ and for each $k \in \{1, \ldots, r \}$, we
 denote  $r' := r+i-n+1$, $r'' := \max(r'-k+1, 0)$, and $r''' := r-k+1$. We define 
\begin{align} \label{eq:formula for tilde D of KR modules of type Dn}
& D_{i,s}^{(k)} := 
\begin{cases}
\left( \prod_{ \substack{p \in [r''', r] \\ q \in [r, n-2+\min(0,r')]} } \frac{1}{\alpha_{p,q}} \right) \left( \prod_{ p \in [r'', r'] }  \frac{ \theta_{p,p}}{ \left( \prod_{q \in [r''', r] } \theta_{p,q} \right)  \left( \prod_{q \in [r', n]}  \alpha_{p,q} \right) } \right), & 1 \leq i \leq  n-2, \\
\left( \prod_{ \substack{p \in [r''',r] \\ q \in [r, n-2]} } \frac{1}{\alpha_{p,q}} \right)   \frac{1}{ \left( \prod_{ p \in [r''',r] } \alpha_{p, \sigma^{r-1}(i)} \right) \left( \prod_{r''' \le p < q \le r} \theta_{p,q} \right) }, & i \in \{n-1,n\},
\end{cases}
\end{align}
where we used the convention that $\theta_{p,q}=1$ if $p=0$ and $\alpha_{p,q}=1$ if $p=0$. 
As in the previous subsection, we also set $d_{i,s}^{(k)} := \frac{D_{i,s}^{(k)}}{D_{i,s+2}^{(k-1)}}$. We also denote $D_{i,s} := D_{i,s}^{(r)}$ and $d_{i,s} := d_{i,s}^{(r)}$. Note that one has $d_{i,s} = D_{i,s} / D_{i,s+2}$.

  \begin{lemma} \label{CartanpairingD}
  For every $1 \leq i \leq j \leq n$ with $i \leq n-1$, $1 \leq p<q \leq n-1$, $1 \leq r<s \leq n-1$, we have that 
  \begin{enumerate}
  \item $(\alpha_{k,i} , \alpha_{l,j}) =
   \begin{cases}
   \delta_{k,l} + 2 \delta_{i,j} -1  & \text{if both $i$ and $j$ are in $\{n-1,n\}$,} \\
   \delta_{k,l} + \delta_{i,j} - \delta_{k-1,j} - \delta_{l-1,i} & \text{otherwise.}
   \end{cases} $ 
  \item  $(\theta_{p,q}, \alpha_{l,j}) = \delta_{l,p} + \delta_{l,q} - \delta_{j+1,p} - \delta_{j+1,q}$. 
  \item $(\theta_{p,q} , \theta_{r,s}) = \delta_{s,p} + \delta_{s,q} + \delta_{r,p} + \delta_{r,q}$. 
  \end{enumerate}
   \end{lemma}

   \begin{prop} \label{lem:tildeD satisfies the initial condition}
   For every $(i,s) \in I_{Q_0}$, one has $\td(Y_{i,s}) = d_{i,s}$. 
   \end{prop}

 \begin{proof}
 Similar arguments as in the proof of Proposition~\ref{lem:tilde D of Yit in type A} show that~\eqref{defDtildexi} can be written as 
 $$ \td_{Q_0}(Y_{i,s}) = \prod_{(j,t) \in I_{Q_0}} \left( \tau^{(\xi(j)-t)/2}(\gamma_j)   \right)^{\tilde{C}_{i,j}(t-s-1) - \tilde{C}_{i,j}(t-s+1)} . $$
 First of all, note that for our choice of orientation we have 
 $$ \tau^{l-1}(\gamma_j) = 
  \begin{cases}
   \alpha_{l,l+j-1} &\text{if $j<n-l$,} \\
   \theta_{i+l-n+1,l} &\text{if $n-l \leq j \leq n-2$,} \\
   \alpha_{l, \sigma^{l-1}(j)} &\text{if $j \in \{n-1,n\}$.}
   \end{cases}
   $$
   We distinguish three distinct cases. 
 
   \textbf{Case 1:} $i \leq n-2$ and $r \leq n-i-1$. In this case, we have $\tau^{r-1}(\gamma_i) = \alpha_{r,r+i-1}$ and the proof is identical to the proof of Proposition~\ref{lem:tilde D of Yit in type A}.
  
   \textbf{Case 2:} $i \leq n-2$ and $r \in \{n-i , \ldots , n-1 \}$. Recall that $r' := i+r-n+1$.  Corollary~\ref{cormultiplicityDtilde} together with Lemma~\ref{CartanpairingD} yield
   \begin{align*}
    \td(Y_{i,s}) &= \frac{1}{\theta_{r',r}}  \prod_{  \substack{ j<n-l \\ 2(l-r) < i-j} } \alpha_{l,l+j-1}^{- \delta_{l,r'} - \delta_{l,r} + \delta_{l+j,r'} + \delta_{l+j,r}}   \\
     & \qquad \qquad \times \prod_{  \substack{ n-l \leq j \leq n-2 \\ 2(l-r) < i-j } } \theta_{j+l-n+1,l}^{-(\delta_{l,r'} + \delta_{l,r} + \delta_{j+l-n+1,r'} + \delta_{j+l-n+1,r})}   \prod_{2(l-r)<i-n+1} \alpha_{l,n}^{-\delta_{l,r'} - \delta_{l,r} + \delta_{l,n}} \\
     &= \frac{1}{\theta_{r',r}}   \prod_{1 \leq j \leq n-r'} \frac{1}{\alpha_{r',r'+j-1}}   \prod_{1 \leq j < n-r} \frac{1}{\alpha_{r,r+j-1}}  \prod_{1 \leq j \leq r'-1} \alpha_{r'-j,r'-1} \prod_{1 \leq j \leq r-1} \alpha_{r-j,r-1}  \\
      & \qquad \qquad \times \prod_{n-r' \leq j \leq n-2} \frac{1}{\theta_{j+r'-n+1,r'}} \prod_{1 \leq l' <r'} \frac{1}{\theta_{l',r}}  \prod_{r'<l<r} \frac{1}{\theta_{r',l}}  \times  \frac{1}{\alpha_{r',n}}  \\
      &= \frac{1}{\theta_{r',r}}   \prod_{r' \leq q \leq n} \frac{1}{\alpha_{r',q}}   \prod_{r \leq q \leq n-2} \frac{1}{\alpha_{r,q}}  \prod_{1 \leq j \leq r'-1} \alpha_{j,r'-1} \prod_{1 \leq j \leq r-1} \alpha_{j,r-1}  \\
      & \qquad \qquad \times \prod_{1 \leq j \leq r'-1} \frac{1}{\theta_{j,r'}} \prod_{1 \leq p \leq r'-1} \frac{1}{\theta_{p,r}}  \prod_{r'<l<r} \frac{1}{\theta_{r',l}}  \\
      &= \frac{\theta_{r', r'} \left( \prod_{ p \in [r'-1] } \alpha_{p,r'-1} \right) \left(  \prod_{ p \in [r-1] }  \alpha_{p,r-1} \right) }{ \left( \prod_{ q \in [r-1] } \theta_{r', q} \right) \left( \prod_{ p \in [r'] } \theta_{p,r} \right) \left( \prod_{ q \in [r', n] } \alpha_{r', q} \right) \left( \prod_{ q \in [r, n-2]} \alpha_{r,q} \right) }. 
     \end{align*}
     On the other hand, we have 
$$ D_{i,s} = \frac{ \prod_{ p \in [i+r-n+1] } \theta_{p,p} }{ \left( \prod_{\substack{ p \in [i+r-n+1] \\ q \in [r] }} \theta_{p,q}  \right) \left( \prod_{ \substack{ p \in [i+r-n+1] \\ q \in [i+r-n+1, n] } } \alpha_{p,q} \right) \left( \prod_{ \substack{ p \in [r] \\ q \in [r, n-2] } } \alpha_{p,q} \right) }. $$
Similarly we have 
$$ D_{i,s+2} = \frac{ \prod_{ p \in [i+r-n] } \theta_{p,p} }{ \left( \prod_{\substack{ p \in [i+r-n]  \\ q \in [r-1]  }} \theta_{p,q}) \right) \left( \prod_{ \substack{ p \in [i+r-n] \\ q \in [i+r-n,n] } } \alpha_{p,q} \right) \left( \prod_{ \substack{ p \in [r-1] \\ q \in [r-1,n-2] } } \alpha_{p,q} \right) }. $$
Therefore 
$$ d_{i,s} = \frac{D_{i,s}}{D_{i,s+2}} =  \frac{\theta_{r', r'} \left( \prod_{ p \in [r'-1] } \alpha_{p,r'-1} \right) \left(  \prod_{ p \in [r-1] }  \alpha_{p,r-1} \right) }{ \left( \prod_{ q \in [r-1] } \theta_{r', q} \right) \left( \prod_{ p \in [r'] } \theta_{p,r} \right) \left( \prod_{ q \in [r', n] } \alpha_{r', q} \right) \left( \prod_{ q \in [r, n-2]} \alpha_{r,q} \right) }. $$
This proves the desired statement in the case $i \leq n-2, n-i \leq r \leq n-1$. 

\textbf{Case 3.} $i \in \{n-1,n\}$. For simplicity, we assume $r$ is odd, as the proof is identical in the other case. 
Corollary~\ref{cormultiplicityDtilde} together with Lemma~\ref{CartanpairingD} yield
  \begin{align*}
    \td(Y_{i,s}) &= \frac{1}{\alpha_{r,i}}  \prod_{  \substack{ j<n-l \\ 2(l-r) < n-1-j} } \alpha_{l,l+j-1}^{- \delta_{l,r} - \delta_{l+j-1,i} + \delta_{r-1,l+j-1} + \delta_{l-1,i}}   \\
     & \qquad \qquad \times \prod_{  \substack{ n-l \leq j \leq n-2 \\ 2(l-r) < n-1-j } } \theta_{j+l-n+1,l}^{-\delta_{r,j+l-n+1} - \delta_{r,l} + \delta_{i,j+l-n} + \delta_{i,l-1}}   \prod_{2(l-r)<0} \alpha_{l,i}^{-(\delta_{r,l} +1)} \alpha_{l, \sigma(i)}^{-(\delta_{r,l} -1)} \\
     &= \frac{1}{\alpha_{r,i}}  \prod_{1 \leq j < n-r} \frac{1}{\alpha_{r,r+j-1}}  \prod_{1 \leq l \leq r-1} \alpha_{l,r-1}  \prod_{1 \leq l' \leq r-1} \frac{1}{\theta_{l',r}}   \times \prod_{1 \leq l \leq r-1} \frac{\alpha_{l, \sigma(i)}}{\alpha_{l,i}} \\
     &=  \frac{ \left( \prod_{ l \in [r-1] } \alpha_{l,r-1} \right)   \left( \prod_{ l \in [r-1] } \alpha_{l,  \sigma^r(i)} \right) }{ \left( \prod_{ q \in [r, n-2] } \alpha_{r,q} \right) \left( \prod_{ l \in [r] } \alpha_{l, \sigma^{r-1}(i)} \right) \left( \prod_{ p \in [r-1] } \theta_{p,r} \right) } . 
     \end{align*}
     On the other hand, we have that 
 $$ D_{i,s} =  \frac{1}{ \left( \prod_{ \substack{ p \in [r] \\ q \in [r,n-2] } } \alpha_{p,q} \right) \left( \prod_{ p \in [r]  }  \alpha_{p,i}  \right) \left( \prod_{1 \le p < q \le r} \theta_{p,q} \right) }. $$
Since $k-1$ is even, we have
$$ D_{i,s+2} =  \frac{1}{ \left( \prod_{ \substack{ p \in [r-1] \\ q \in [r-1, n-2] }  } \alpha_{p,q} \right) \left( \prod_{ p \in [r-1]  }  \alpha_{p,\sigma(i)} \right) \left( \prod_{1 \le p < q \le r-1} \theta_{p,q} \right) }. $$
Therefore
$$ d_{i,s} = \frac{D_{i,s}}{D_{i,s+2}} =  \frac{ \left( \prod_{ p \in [r-1]  } \alpha_{p,r-1} \right)   \left( \prod_{ p \in [r-1]  }  \alpha_{p,\sigma(i)}   \right) }{ \left( \prod_{ q \in [r,n-2] } \alpha_{r,q} \right) \left( \prod_{ p \in [r] } \alpha_{p,i} \right) \left( \prod_{ p \in [r-1] } \theta_{p,r} \right) }. $$
 \end{proof}
This concludes the proof. 


   \begin{corollary} \label{corDtildeinitialtypeD}
  For every $(i,s) \in I_{Q_0}$ one has $\td( \tchi_q(X_{i,s})) = D_{i,s}$. 
      \end{corollary}
      
       \begin{proof}
  The proof is identical to the proof of Corollary~\ref{corDtildeinitialtypeA}, using Proposition~\ref{lem:tildeD satisfies the initial condition}.  
       \end{proof}

\begin{thm} \label{thmKRmonotonicD}
For every $(i,s) \in I_{Q_0}$ and every $k \in \{1, \ldots , r \}$ we have $\td(\tchi_q(X_{i,s}^{(k)})) = D_{i,s}^{(k)}$.
\end{thm}

\begin{proof}
By Corollary~\ref{corDtildeinitialtypeD}, the statement is true if $r=k$. Similarly to the proof of Theorem~\ref{thmKRmonotonicA}, we prove that the fractions $D_{i,s}^{(k)}$ satisfy the recursive equations:
\begin{align*}
& D_{i,s}^{(k)} D_{i,s-2}^{(k)} = D_{i,s-2}^{(k+1)} D_{i,s}^{(k-1)} + \prod_{j \sim i} D_{j,s-1}^{(k)}.
\end{align*}

{\bf Case 1.} $i \in [n-2]$, $r \in [n-i-2]$. The proof of this case is the same as the case of type $A_n$.

{\bf Case 2.} $i \in [n-2]$, $r \in [n-i-1,n-2]$. In this case, $r'=i+r-n+1 \ge 0$. We will prove the case where $r$ is odd. The case where $r$ is even is similar.

We have that
\begin{align*}
& D_{i,s}^{\left(k\right)} = \frac{\prod_{p \in [r'', r']} \theta_{p,p}}{ \left(\prod_{ \substack{ p \in [r'', r'] \\ q \in [r''', r] } } \theta_{p,q} \right) \left(\prod_{ \substack{ p \in [r'', r'] \\ q \in [r', n] } } \alpha_{p,q} \right) \left(\prod_{ \substack{ p \in [r''', r] \\ q \in [r, n-2 ] } } \alpha_{p,q} \right) }.
\end{align*}
Since $\frac{\xi_i - (s-2)+2}{2} = r+1$, we have that
\begin{align*}
D_{i,s-2}^{\left(k\right)} = \frac{\prod_{p \in [r''+1, r'+1]} \theta_{p,p}}{ \left(\prod_{ \substack{ p \in [r''+1, r'+1] \\ q \in [r'''+1, r+1] } } \theta_{p,q} \right) \left(\prod_{ \substack{ p \in [r''+1, r'+1] \\ q \in [r'+1, n] } } \alpha_{p,q} \right) \left(\prod_{ \substack{ p \in [r'''+1, r+1] \\ q \in [r+1, n-2 ] } } \alpha_{p,q} \right) }.
\end{align*}
Similarly we have
\begin{align*}
D_{i,s-2}^{\left(k+1\right)} = \frac{\prod_{p \in [r'', r'+1]} \theta_{p,p}}{ \left(\prod_{ \substack{ p \in [r'', r'+1] \\ q \in [r''', r+1] } } \theta_{p,q} \right) \left(\prod_{ \substack{ p \in [r'', r'+1] \\ q \in [r'+1, n] } } \alpha_{p,q} \right) \left(\prod_{ \substack{ p \in [r''', r+1] \\ q \in [r+1, n-2 ] } } \alpha_{p,q} \right) }.
\end{align*}
We have
\begin{align*}
D_{i,s}^{\left(k-1\right)} = \frac{\prod_{p \in [r''+1, r']} \theta_{p,p}}{ \left(\prod_{ \substack{ p \in [r''+1, r'] \\ q \in [r'''+1, r] } } \theta_{p,q} \right) \left(\prod_{ \substack{ p \in [r''+1, r'] \\ q \in [r', n] } } \alpha_{p,q} \right) \left(\prod_{ \substack{ p \in [r'''+1, r] \\ q \in [r, n-2 ] } } \alpha_{p,q} \right) }.
\end{align*}
Since $\frac{\xi_{i-1} - (s-1)+2}{2} = \frac{\xi_{i} +1 - (s-1)+2}{2} = r+1$, we have that
\begin{align*}
D_{i-1,s-1}^{\left(k\right)} = \frac{\prod_{p \in [r'', r']} \theta_{p,p}}{ \left(\prod_{ \substack{ p \in [r'', r'] \\ q \in [r'''+1, r+1] } } \theta_{p,q} \right) \left(\prod_{ \substack{ p \in [r'', r'] \\ q \in [r', n] } } \alpha_{p,q} \right) \left(\prod_{ \substack{ p \in [r'''+1, r+1] \\ q \in [r+1, n-2 ] } } \alpha_{p,q} \right) }.
\end{align*}
Divide $D_{i,s}^{\left(k\right)} D_{i,s-2}^{\left(k\right)} - D_{i,s-2}^{\left(k+1\right)} D_{i,s}^{\left(k-1\right)}$ by 
\begin{align*}
& \frac{\prod_{p \in [r''+1, r']} \theta_{p,p}}{ \left(\prod_{ \substack{ p \in [r''+1, r'] \\ q \in [r'''+1, r] } } \theta_{p,q} \right) \left(\prod_{ \substack{ p \in [r''+1, r'] \\ q \in [r', n] } } \alpha_{p,q} \right) \left(\prod_{ \substack{ p \in [r'''+1, r] \\ q \in [r, n-2 ] } } \alpha_{p,q} \right) } \times
\\
& \times \frac{\prod_{p \in [r''+1, r'+1]} \theta_{p,p}}{ \left(\prod_{ \substack{ p \in [r''+1, r'+1] \\ q \in [r'''+1, r+1] } } \theta_{p,q} \right) \left(\prod_{ \substack{ p \in [r''+1, r'+1] \\ q \in [r'+1, n] } } \alpha_{p,q} \right) \left(\prod_{ \substack{ p \in [r'''+1, r+1] \\ q \in [r+1, n-2 ] } } \alpha_{p,q} \right) },
\end{align*}
we obtain
\begin{align*}
& \frac{\theta_{r'',r''}}{\left(\prod_{p \in [r'',r']} \theta_{p,r'''}\right) \left(\prod_{q \in [r'''+1, r]} \theta_{r'', q} \right) \left(\prod_{q \in [r'+1, n]} \alpha_{r'', q}\right) \left(\prod_{q \in [r+1, n-2 ]} \alpha_{r''', q}\right) } \times  \\
& \times \left( \frac{1}{\alpha_{r'', r'} \alpha_{r''',r}} - \frac{1}{\theta_{r'+1, r'''} \theta_{r'', r+1} }  \right) \\
& = \frac{\theta_{r'',r''}}{\left(\prod_{p \in [r'',r']} \theta_{p,r'''}\right) \left(\prod_{q \in [r'''+1, r]} \theta_{r'', q} \right) \left(\prod_{q \in [r'+1, n]} \alpha_{r'', q}\right) \left(\prod_{q \in [r+1, n-2 ]} \alpha_{r''', q}\right) } \times \\
& \times \frac{ \theta_{r'',r'''} \theta_{r'+1,r+1} }{\alpha_{r'', r'} \alpha_{r''',r}\theta_{r'+1, r'''} \theta_{r'', r+1} } \\
& = \frac{\theta_{r'',r''}  \theta_{r'+1, r+1}}{ \left(\prod_{p \in [r''+1,r'+1]} \theta_{p,r'''}\right) \left(\prod_{q \in [r'''+1, r+1]} \theta_{r'', q} \right) \left(\prod_{q \in [r', n]} \alpha_{r'', q}\right) \left(\prod_{q \in [r, n-2 ]} \alpha_{r''', q}\right) }.
\end{align*}
{\bf Subcase 2.1.} $i \in [n-3]$. Since $\frac{\xi_{i+1} - \left(s-1\right)+2}{2} = \frac{\xi_{i} -1 - \left(s-1\right)+2}{2} = r$, we have that
\begin{align*}
D_{i+1,s-1}^{\left(k\right)} = \frac{\prod_{p \in [r''+1, r'+1]} \theta_{p,p}}{ \left(\prod_{ \substack{ p \in [r''+1, r'+1] \\ q \in [r''', r] } } \theta_{p,q} \right) \left(\prod_{ \substack{ p \in [r''+1, r'+1] \\ q \in [r'+1, n] } } \alpha_{p,q} \right) \left(\prod_{ \substack{ p \in [r''', r] \\ q \in [r, n-2 ] } } \alpha_{p,q} \right) }.
\end{align*}
Therefore $D_{i,s}^{(k)} D_{i,s-2}^{(k)} - D_{i,s-2}^{(k+1)} D_{i,s}^{(k-1)}=D_{i-1,s-1}^{(k)} D_{i+1,s-1}^{(k)}$.

{\bf Subcase 2.2.} $i =n-2$. In this case, $r' = i+r-n+1 = r-1 \geq 0$, $r''=\max(r'-k+1, 0) = \max(r-k,0)= r-k = r'-k+1= r'''-1$, $r'''=r-k+1 \geq 1$. 

We have $\frac{\xi_{n-1} - (s-1)+2}{2} = \frac{\xi_{n-2} -1 - (s-1)+2}{2} = r$. Since $r$ is odd, we have
\begin{align*}
D_{n-1,s-1}^{(k)} =  \frac{1}{ (\prod_{\substack{p \in [r''', r] \\ q \in [r, n-2]}} \alpha_{p,q} ) ( \prod_{p \in [r''', r]} \alpha_{p,n-1} ) (\prod_{r''' \le p < q \le r} \theta_{p,q}) }.
\end{align*}
We have $\frac{\xi_{n} - (s-1)+2}{2} = \frac{\xi_{n-2} -1 - (s-1)+2}{2} = r$. Since $r$ is odd, we have
\begin{align*}
D_{n,s-1}^{\left(k\right)} =  \frac{1}{ \left(\prod_{\substack{p \in [r''', r] \\ q \in [r, n-2]}} \alpha_{p,q} \right) \left( \prod_{p \in [r''', r]} \alpha_{p,n} \right) \left(\prod_{r''' \le p < q \le r} \theta_{p,q}\right) }.
\end{align*}
Using $r''+1 = r'''$ and 
$$
\prod_{\substack{ p \in [r''+1, r] \\ q \in [r''+1, r+1] }} \theta_{p,q} = \prod_{\substack{ p \in [r''', r] \\ q \in [r''+1, r+1] }} \theta_{p,q} = \frac{ \left( \prod_{r''' \le p < q \le r} \theta_{p,q}\right)^2 \left( \prod_{p \in [r''',r]} \theta_{p,r+1} \right) \left( \prod_{p \in [r''', r]} \theta_{p,p} \right) }{ \prod_{q \in [r''', r]} \theta_{r''',q} },
$$
we have that
$$D_{n-2,s}^{(k)} D_{n-2,s-2}^{(k)} - D_{n-2,s-2}^{(k+1)} D_{n-2,s}^{(k-1)}=D_{n-3,s-1}^{(k)} D_{n-1,s-1}^{(k)} D_{n,s-1}^{(k)}.$$
{\bf Case 3.} $i \in [n-1, n]$, $r \in [n-2]$. We will prove the case where $r$ is odd. The case where $r$ is even is similar.

We have that
\begin{align*}
& D_{i,s}^{\left(k\right)} = \frac{1}{ \left(\prod_{\substack{ p \in [r''',r] \\ q \in [r, n-2] }} \alpha_{p,q}\right) \left(\prod_{p \in [r''',r]} \alpha_{p,i} \right) \left(\prod_{r''' \leq p < q \leq r} \theta_{p,q}\right)   }.
\end{align*}
Since $\frac{\xi_i - \left(s-2\right)+2}{2} = r+1$ and $r+1$ is even, we have that
\begin{align*}
D_{i,s-2}^{\left(k\right)} =  \frac{1}{ \left(\prod_{\substack{ p \in [r'''+1,r+1] \\ q \in [r+1, n-2] }} \alpha_{p,q}\right) \left(\prod_{p \in [r'''+1,r+1]} \alpha_{p, i'} \right) \left(\prod_{r'''+1 \le p < q \le r+1} \theta_{p,q}\right)   },
\end{align*}
and
\begin{align*}
D_{i,s-2}^{\left(k+1\right)} = \frac{1}{ \left(\prod_{\substack{ p \in [r''',r+1] \\ q \in [r+1, n-2] }} \alpha_{p,q}\right) \left(\prod_{p \in [r''',r+1]} \alpha_{p,i'} \right) \left(\prod_{r''' \le p < q \le r+1} \theta_{p,q}\right)   }.
\end{align*}
We also have
\begin{align*}
D_{i,s}^{\left(k-1\right)} = \frac{1}{ \left(\prod_{\substack{ p \in [r'''+1,r] \\ q \in [r, n-2] }} \alpha_{p,q}\right) \left(\prod_{p \in [r'''+1,r]}  \alpha_{p,i} \right) \left(\prod_{r'''+1 \le p < q \le r} \theta_{p,q}\right)   }.
\end{align*}
We have $\frac{\xi_{n-2} - (s-1)+2}{2} = \frac{\xi_{i} +1 - (s-1)+2}{2} = r+1$. 
Therefore the integers $r', r'', r'''$ corresponding to the couple $(n-2,s-1)$ are respectively given by 
\begin{align*}
& (n-2)+(r+1)-n+1=r, \\
& \max(r'_{n-2, s-1}-k+1, 0)=\max(r-k+1, 0) = r-k+1 = r''', \\
&  r_{n-2,s-1}-k+1=r-k+2 = r'''+1. 
\end{align*}
It follows that
\begin{align*}
D_{n-2,s-1}^{\left(k\right)} = \frac{\prod_{p \in [r''', r]} \theta_{p,p}}{ \left(\prod_{ \substack{ p \in [r''', r] \\ q \in [r'''+1, r+1] } } \theta_{p,q} \right) \left(\prod_{ \substack{ p \in [r''', r] \\ q \in [r, n] } } \alpha_{p,q} \right) \left(\prod_{ \substack{ p \in [r'''+1, r+1] \\ q \in [r+1, n-2] } } \alpha_{p,q} \right) }.
\end{align*}
Divide $D_{i,s}^{(k)} D_{i,s-2}^{(k)} - D_{i,-2}^{(k+1)} D_{i,s}^{(k-1)}$ by 
\begin{align*}
& \frac{1}{ \left(\prod_{ \substack{ p \in [r'''+1, r] \\ q \in [r, n-2] } } \alpha_{p,q} \right) \left(\prod_{   p \in [r'''+1, r]   } \alpha_{p,i} \right) \left(\prod_{ r'''+1 \le p < q \le r } \theta_{p,q} \right) } \times
\\
& \times \frac{1}{ \left(\prod_{ \substack{ p \in [r'''+1, r+1] \\ q \in [r+1, n-2] } } \alpha_{p,q} \right) \left(\prod_{ p \in [r'''+1, r+1]   }  \alpha_{p,i'} \right) \left(\prod_{ r'''+1 \le p < q \le r+1 } \theta_{p,q} \right) },
\end{align*}
we obtain
\begin{align*}
& \frac{1}{ \left(\prod_{q \in [r+1, n-2]} \alpha_{r''',q}\right) \left(\prod_{q \in [r'''+1, r]} \theta_{r''',q}\right) } \left( \frac{1}{\alpha_{r''',i} \alpha_{r''',r}}  - \frac{1}{\alpha_{r''',i'} \theta_{r''',r+1}} \right) \\
& = \frac{\alpha_{r+1, i'} \theta_{r''',r'''}}{ \left(\prod_{q \in [r, n]} \alpha_{r''',q}\right) \left(\prod_{q \in [r'''+1, r+1]} \theta_{r''',q}\right) }.
\end{align*}
Using $ \prod_{\substack{p \in [r''',r] \\ q \in [r'''+1, r+1]}} \theta_{p,q} = \left( \prod_{r''' \le p < q \le r+1} \theta_{p,q} \right) \left( \prod_{r'''+1 \le p < q \le r}  \theta_{p,q}\right) \left(\prod_{p \in [r'''+1,r]} \theta_{p,p}\right) $,
we conclude that
$$D_{i,s}^{(k)} D_{i,s-2}^{(k)} - D_{i,s-2}^{(k+1)} D_{i,s}^{(k-1)}=D_{n-2,s-1}^{(k)}.$$

\end{proof}

\subsection{Types $E_6,E_7,E_8$}
 \label{DtildeKRE}

Assume $\mathfrak{g}$ is of type $E_r, r=6,7,8$ and $Q=Q_0$ is the orientation considered in Figure~\ref{figQ0}. Using SageMath \cite{Sage}, we use the $T$-systems~\eqref{eq:T-system relations} to compute the images under $\td$ of the truncated $q$-characters of Kirillov-Reshetikhin modules, and thus in particular fundamental representations.  

For every fundamental module $L(Y_{i,s})$ in $\mathcal{C}_Q$, we computed the graded character of the corresponding cuspidal representation $\mathcal{F}(L(Y_{i,s})) = S_{\beta_{\varphi(i,s)}}$ of the quiver Hecke algebra using the algorithm in \cite{BKM, Leclerc}. Then we obtain the corresponding ungraded character and we apply the map $\overline{D}$ to the resulting character. Since $\td(L(Y_{i,s}))$ and $\overline{D}(S_{\beta_{\varphi(i,s)}})$ are rational functions in $\alpha_1, \ldots, \alpha_n$, to check that they are equal, it suffices to check that they are equal for a few choices of numbers for $\alpha_1, \ldots, \alpha_n$. In this way, we verified that $\td(L(Y_{i,s})) = \overline{D}(S_{\beta_{\varphi(i,s)}})$. 

As we already mentioned in Section~\ref{DbarcuspidalE}, the formulas of $\td_Q$ of Kirillov-Reshetikhin modules in type $E$ can be very complicated and it does not seem possible to write them in a form similar to~\eqref{eq:formula for tilde D of KR modules of type An} or~\eqref{eq:formula for tilde D of KR modules of type Dn}. This motivates us to believe that there is no general formula for the images under $\td_Q$ of the truncated $q$-characters of Kirillov-Reshetikhin modules (or even fundamental modules) for arbitrary orientations $Q$ in any simply-laced type. 

The SageMath program to verify the above can be found in the link: \url{https://drive.google.com/drive/folders/1jXW8WG0p_01GkEqYU9s8tLBOvlmlZnT6?usp=sharing}.



  \section{Proofs of the main results}
   \label{sectionproofs}
  
   This section is devoted to the proofs of Theorems~\ref{mainthm1} and~\ref{mainthm2}. We proceed in the following way: we begin by proving Theorem~\ref{mainthm1} in the particular case of the  orientation $Q_0$, combining the results obtained in Sections~\ref{Dbarcuspidals} and~\ref{DtildeKRmonotonic}. This allows us to prove Theorem~\ref{mainthm2} for the standard seed $\s^{\mathbf{i}_{Q_0}}$. The statement for arbitrary reduced expressions of $w_0$ then follows from Theorem~\ref{thmpropagation}. Finally, we prove Theorem~\ref{mainthm1} in full generality i.e. for an arbitrary orientation $Q$ of the Dynkin diagram of $\mathfrak{g}$, by using Theorem~\ref{mainthm2} with $\mathbf{i} = \iQ$. 
   
   \begin{proof}[\textbf{Proof of Theorem~\ref{mainthm1} in the case of the orientation $Q_0$.}]
   
  Let us fix the orientation $Q_0$ of the Dynkin diagram of $\mathfrak{g}$ as in Figure~\ref{figQ0}. Recall from Sections~\ref{remindCQ} and~\ref{remindsimpleKLR} that the dual root vectors associated to $\mathbf{i}_{Q_0}$ are categorified on the one hand by the fundamental representations of $\C_{Q_0}$ (see Theorem~\ref{thmHL15}) and on the other hand by the cuspidal representations of $R-mod$ (for the convex ordering on $\Phi_{+}$ corresponding to $\mathbf{i}_{Q_0}$, see Section~\ref{Dbarcuspidals}). More precisely, for any $i \in I$ and any $1 \leq r \leq n_{Q_0}(i)$, we have
   $$ [L(Y_{i,\xi(i)-2(r-1)})] = [S_{\tau^{r-1}(\gamma_i)}] $$ 
   in $\CN$.  Moreover, the dual root vectors generate $\CN$ as an algebra. Thus, in order to prove that $\td_{Q_0}$ and $\barD$ coincide, it suffices to prove that they agree on the dual root vectors. 
   
    Recall (see for example Section~\ref{remindquantumaffine}) that the fundamental representation $L(Y_{i,\xi(i)-2(r-1)})$ is the Kirillov-Reshetikhin module $X_{i,\xi(i)-2(r-1)}^{(1)}$. Thus we apply the formulas obtained in Section~\ref{DtildeKRmonotonic} with $k=1$. When $\mathfrak{g}$ is of type $A_n , n \geq 1$, Theorem~\ref{thmKRmonotonicA} yields
   $$ \td_{Q_0} \left(  \tchi_q(L(Y_{i,\xi(i)-2(r-1)})) \right) = D_{i,(\xi(i)-2(r-1))}^{(1)} = \prod_{r \leq q  \leq r+i-1} \frac{1}{\alpha_{r,q}} . $$
   On the other hand, Equation~\eqref{DbarofSalpha} yields
   $$ \barD([S_{\tau^{r-1}(\gamma_i)}]) = \barD([S_{\alpha_{r,r+i-1}}]) =  \prod_{r \leq q  \leq r+i-1} \frac{1}{\alpha_{r,q}}. $$
   This proves the Theorem in type $A_n$. 
    When $\mathfrak{g}$ is of type $D_n , n \geq 4$, Theorem~\ref{thmKRmonotonicD} yields
  \begin{align*}
  & \td_{Q_0} \left( \tchi_q(L(Y_{i,\xi(i) -2(r-1)})) \right) = D_{i,(\xi(i)-2(r-1))}^{(1)}  \\
    &= 
\begin{cases}
\prod_{r \leq q  \leq r+i-1} \frac{1}{\alpha_{r,q}} & \text{if $i \leq n-2$ and $r \leq n-i-1$,} \\
 \frac{\theta_{r',r'}}{\theta_{r',r}} \times \prod_{r' \leq q  \leq n} \frac{1}{\alpha_{r',q}} \times  \prod_{r \leq q  \leq n-2} \frac{1}{\alpha_{r,q}}  & \text{if $i \leq n-2$ and $r \geq n-i$,} \\
  \frac{1}{\alpha_{r,\sigma^{r-1}(i)}} \times  \prod_{r \leq q  \leq n-2} \frac{1}{\alpha_{r,q}} &\text{if $i \in \{n-1,n\}$.}
 \end{cases}   
      \end{align*}
     On the other hand, if $i \leq n-2$ and $r \leq n-i-1$, then one has $\tau^{r-1}(\gamma_i) = \alpha_{r,r+i-1}$ and the conclusion is the same as in type $A_n$; if $i \leq n-2$ and $r \geq n-i$, then $\tau^{r-1}(\gamma_i) = \theta_{r',r}$ and Proposition~\ref{DbarofStheta} yields
     $$ \barD([S_{\theta_{r',r}}]) =  \frac{\theta_{r',r'}}{\alpha_{r,r} \cdots \alpha_{r,n-2}  \alpha_{r',r'} \cdots \alpha_{r',n} \theta_{r',r}}  . $$
     This coincides with the above expression of $\td_{Q_0} \left( \tchi_q(X_{i,\xi(i)-2(r-1)}^{(1)}) \right)$ in this case. 
     Finally, if $i \in \{n-1,n\}$, then $\tau^{r-1}(\gamma_i) = \alpha_{r,\sigma^{r-1}(i)}$ and thus we have 
     $$ \barD([S_{\tau^{r-1}(\gamma_i)}]) = \barD([S_{\alpha_{r,\sigma^{r-1}(i)}}]) =   \frac{1}{\alpha_{r,\sigma^{r-1}(i)}} \times  \prod_{r \leq q  \leq n-2} \frac{1}{\alpha_{r,q}}. $$
     This proves the Theorem in type $D_n$. 
     For the types $E_6, E_7$ and $E_8$ we check by computer that the respective values given by $\barD$ (see Section~\ref{DbarcuspidalE}) and $\td_{Q_0}$ (see Section~\ref{DtildeKRE}) agree on the dual root vectors. 
     
   \end{proof}
   
 We can now use the properties of $\td_{\xi}$ established in Section~\ref{ABCforDtilde} to prove the second main result of the present paper, which was stated as a Conjecture in \cite{Casbi3}  ({{\cite[Conjecture 5.5]{Casbi3}}}). 
 
  \begin{thm} \label{mainthm2}
  Let $\mathfrak{g}$ be a simple Lie algebra of simply-laced type. Then for any reduced expression $\mathbf{i}$ of $w_0$, the flag minors $x_1^{\mathbf{i}} , \ldots , x_N^{\mathbf{i}}$ satisfy $\barD(x_j^{\mathbf{i}}) = 1/P_j^{\mathbf{i}}$ where $P_j^{\mathbf{i}}$ is a product of positive roots. Furthermore, one has $[\beta;P_j^{\mathbf{i}}] - [\beta;P_{j_{+}}^{\mathbf{i}}] \leq 1$ for any $\beta \in \Phi_{+}$ and any $j$ such that $j_{+} \leq N$, and the polynomials $P_1^{\mathbf{i}} , \ldots , P_N^{\mathbf{i}}$ satisfy the identities
  $$ \forall 1 \leq j \leq N, \enspace  P_j^{\mathbf{i}} P_{j_{-}}^{\mathbf{i}}  = \beta_j \prod_{\substack{l<j<l_{+} \\ i_l \sim i_j}}  P_l^{\mathbf{i}} . $$

   \end{thm}

  \begin{proof}
  
  We begin by proving the desired statement for $\mathbf{i} = \mathbf{i}_{Q_0}$. We  deduce the relations $(A_{\mathbf{i}_{Q_0}}), (B_{\mathbf{i}_{Q_0}}), (C_{\mathbf{i}_{Q_0}})$ respectively from Propositions~\ref{AforDtilde}, ~\ref{BforDtilde} and~\ref{CforDtilde}. Recall from Section~\ref{definitionDtilde} the natural embedding  $\CN  \hookrightarrow \Axi$. Recall also that for any $Q$,  the flag minors $x_j^{\iQ} ,1 \leq j \leq N$ are identified under this embedding with the cluster variables $x_t^{\wiQ} , 1 \leq  t \leq N$ of Hernandez-Leclerc's initial seed in $\Axi$ (see Section~\ref{remindCQ}). Thus in the proof below, we will use the notation  $x_j^{\mathbf{i}_{Q_0}}$  for both the flag minor of the standard seed $\s^{\mathbf{i}_{Q_0}}$ in $\CN$ and its image in $\Axi$.

      Denote $\mathbf{i}_{Q_0} = (i_1, \ldots , i_N)$ and let $j \in \{1, \ldots, N \}$. Then by Theorem~\ref{mainthm1} with $Q=Q_0$ proved above, one has 
  $$  \barD(x_j^{\mathbf{i}_{Q_0}}) = \td_{Q_0} \left( \iota(x_j^{\mathbf{i}_{Q_0}}) \right) . $$
 Therefore by Proposition~\ref{AforDtilde} we have 
  $$ \barD(x_j^{\mathbf{i}_{Q_0}}) = \prod_{\beta \in \Phi_{+}} \frac{1}{\beta^{n_j(\beta)}} $$
  where $n_j(\beta)$ is a nonnegative integer for each $\beta \in \Phi_{+}$. This proves that the relation $(A_{\mathbf{i}_{Q_0}})$ holds.

For the relation  $(B_{\mathbf{i}_{Q_0}})$, we denote $P_j := \left( \barD(x_j^{\mathbf{i}_{Q_0}}) \right)^{-1}$ for each $1 \leq j \leq N$. Then using Theorem~\ref{mainthm1} with $Q=Q_0$, we have 
  \begin{align*}
  P_j P_{j_{-}}  & = \left( \barD(x_j^{\mathbf{i}_{Q_0}}) \barD(x_{j_{-}}^{\mathbf{i}_{Q_0}}) \right)^{-1} =  \left(\td_{Q_0}(\iota(x_j^{\mathbf{i}_{Q_0}}))\td_{Q_0}(\iota(x_{j_{-}}^{\mathbf{i}_{Q_0}})) \right)^{-1} \\
  &= \beta_j \prod_{ \substack{ r<t<r_{+} \\ i_r \sim i  }} \td_{Q_0}( \iota (x_r^{\mathbf{i}_{Q_0}}))^{-1} \enspace \text{by Proposition~\ref{BforDtilde}} \\
  &= \beta_j  \prod_{ \substack{ r<t<r_{+} \\ i_r \sim i  }}  \barD(x_r^{\mathbf{i}_{Q_0}})^{-1} = \beta_j \prod_{ \substack{ r<t<r_{+} \\ i_r \sim i  }} P_r \enspace  \text{using again Theorem~\ref{mainthm1} for $Q_0$}.
  \end{align*}

  For the relation $(C_{\mathbf{i}_{Q_0}})$, let $j \in \{1, \ldots , N\}$ such that $j_{+} \leq N$ and let $(i,p) := \varphi^{-1}(j) \in I_{Q_0}$. By~\eqref{defvarphi} we have $j_{+} = \varphi(i,p-2)$. Thus applying Theorem~\ref{mainthm1} for $Q_0$ we get
 \begin{align*}
 P_{j_{+}} &= \barD(x_{j_{+}}^{\mathbf{i}_{Q_0}})^{-1}= \td_{Q_0} \left( \iota (x_{j_{+}}^{\mathbf{i}_{Q_0}}) \right)^{-1} =  \td_{Q_0} \left( \tchi_q(X_{i,p-2}) \right)^{-1} \\
 &= \td_{Q_0}(Y_{i,p-2})^{-1} \td_{Q_0} \left( \tchi_q(X_{i,p}) \right)^{-1} = \td_{Q_0}(Y_{i,p-2})^{-1}  \barD(x_j^{\mathbf{i}_{Q_0}})^{-1} = \td_{Q_0}(Y_{i,p-2})^{-1} P_j.
 \end{align*}
  Hence for each $\beta \in \Phi_{+}$, one has $[\beta ; P_j] - [\beta ; P_{j_{+}}] = [\beta ; \td_{Q_0}(Y_{i,p-2})]$. The conclusion follows from Proposition~\ref{CforDtilde}. 
  
 \smallskip  
  
   We have proved the desired statement in the case $\mathbf{i} = \mathbf{i}_{Q_0}$. 
  The conclusion for arbitrary reduced expressions of $w_0$ is provided by Theorem~\ref{thmpropagation} ({{\cite[Theorem 5.6]{Casbi3}}}), which ensures that the properties  $(A_{\mathbf{i}}), (B_{\mathbf{i}}), (C_{\mathbf{i}})$ hold for each standard seed $\s^{\mathbf{i}}$ of $\CN$.  This finishes the proof of Theorem~\ref{mainthm2}.

   \end{proof}

    \begin{remark}
Alternatively, the relation $(A_{\mathbf{i}_{Q_0}})$ can also be deduced from Corollaries~\ref{corDtildeinitialtypeA} and~\ref{corDtildeinitialtypeD} when $\mathfrak{g}$ is of type $A_n , n \geq 1$ or $D_n , n \geq 4$, and can be checked by computer when $\mathfrak{g}$ is of type $E_6, E_7$ or $E_8$. For the type $A_n$, we recover the formulas of {{\cite[Lemma 7.2]{Casbi3}}} which were there obtained using certain results from \cite{Casbi}.
   \end{remark}

  Now we can use Theorem~\ref{mainthm2} to prove Theorem~\ref{mainthm1} in full generality i.e. for an arbitrary orientation of the Dynkin diagram of $\mathfrak{g}$.
  
  \begin{proof}[\textbf{Proof of Theorem~\ref{mainthm1}: the general case.}]
 Let $Q$ be an arbitrary orientation of the Dynkin diagram of $\mathfrak{g}$ and let us fix $\iQ$ a reduced expression of $w_0$ adapted to $Q$. By Theorem~\ref{mainthm2}, the standard seed $\s^{\iQ}$ of $\CN$ satisfies Properties $(A_{\iQ}), (B_{\iQ})$ and $(C_{\iQ})$. So we have 
 $$ \forall 1 \leq j \leq N, \enspace \barD(x_j^{\iQ}) \barD(x_{j_{-}}^{\iQ})  =  \beta_j^{-1} \prod_{l<j<l_{+}}  \barD(x_l^{\iQ}) . $$
 On the other hand, by Proposition~\ref{BforDtilde} the rational fractions $\td_Q \left( \iota (x_j^{\iQ}) \right) , 1 \leq j \leq N$ satisfy the same relations. Thus by a straightforward induction we have $\barD(x_j^{\iQ}) = \td_Q \left( \iota (x_j^{\iQ}) \right)$ for each $1 \leq j \leq N$. 
 As $\barD$ and $\td_Q \circ \iota$ are both algebra morphisms and the ring $\CN$ has a cluster structure with a seed given by $\s^{\iQ}$, this implies that $\td_Q \circ \iota = \barD$ on the whole algebra $\CN$. 
   \end{proof}

 \section{Application to the generalized quantum affine Schur-Weyl duality}
  \label{sectionappliSW}
 
 In this section we provide a representation-theoretic interpretation of Theorem~\ref{mainthm1} from the perspective of Kang-Kashiwara-Kim-Oh's generalized quantum affine Schur-Weyl duality \cite{KKKOSelecta}.

For any simply-laced type Lie algebra $\mathfrak{g}$ and for any orientation $Q$ of the Dynkin graph of $\mathfrak{g}$, Kang-Kashiwara-Kim-Oh \cite{KKKOSelecta} defined a monoidal functor $\mathcal{F}_Q$ from the category $R-mod$ of finite-dimensional  modules over the quiver Hecke algebras associated to $\mathfrak{g}$ (see Section~\ref{remindKLR}) to the category $\CQ$. This functor $\mathcal{F}_Q$, called the generalized quantum affine Schur-Weyl duality functor was moreover proved by Fujita \cite{Fujita} to be an equivalence of categories. However, the structures of the objects themselves are a priori very different. For instance the objects in $R-mod$ carry a natural $\mathbb{Z}$-grading (see Section~\ref{remindKLR}) which is not the case for the objects of $\CQ$. On the other hand, the classes of the representations in $\CQ$ can be described via Frenkel-Reshetikhin's (truncated) $q$-character \cite{FR1} which allows to perform computations in certain tori (such as $\YQ$), whereas the characters of the objects in $R-mod$  take values in the shuffle algebra, which is much more difficult to tackle with. 
 
   Theorem~\ref{mainthm1} yields a surprising connection between the weight subspaces decompositions of $M$ and $\mathcal{F}_Q(M)$ for every object $M$ in $\CQ$. Indeed, by Theorem~\ref{mainthm1}, we have 
$$ \td_Q(\tchi_q(M)) = \barD([\mathcal{F}_Q(M)]) . $$
   By definition, the truncated $q$-character of $M$ encodes the dimensions of certain of the loop weight spaces of $M$ (see Section~\ref{remindquantumaffine}). Hence, recalling Equation~\eqref{eqnDbarKLR}, the previous equality can be written as 
\begin{equation} \label{eqnDbarDtilde}
  \sum_{\mathfrak{m}' \preccurlyeq \mathfrak{m}} \dim(M_{\mathfrak{m}'}) \td_Q(\mathfrak{m}') =  \sum_{\mathbf{j}=(j_1, \ldots j_d)} \dim \left( (\mathcal{F}_Q(M))_{\mathbf{j}} \right) \barD_{\mathbf{j}} 
  \end{equation}
  where for each $\mathbf{j} := (j_1, \ldots , j_d)$,
  $$\barD_{\mathbf{j}} := \frac{1}{\alpha_{j_1} (\alpha_{j_1} + \alpha_{j_2}) \cdots (\alpha_{j_1} + \cdots + \alpha_{j_d})} . $$
  The sum on the left hand-side runs over all monomials $\mathfrak{m}' \in \YQ$ that are smaller than $\mathfrak{m}$ for the Nakajima ordering (see Section~\ref{remindquantumaffine}), and on the right hand-side $(\mathcal{F}_Q(M))_{\mathbf{j}} := e(\mathbf{j}) \cdot \mathcal{F}_Q(M)$ denotes the weight subspace given by the action of the idempotent $e(\mathbf{j})$ on $\mathcal{F}_Q(M)$. Equation~\eqref{eqnDbarDtilde} is an explicit identity between rational fractions in $\mathbb{C}(\alpha_1, \ldots , \alpha_n)$ involving the dimensions of the weight subspaces of a representation of $\CQ$ on the one hand and those of the corresponding object in $R-mod$ via the generalized Schur-Weyl duality functor on the other hand.

    We now provide a concrete illustration of this fact. For any object $M$ in $\CQ$ with $\tchi_q(M) := \sum_{\mathfrak{m}} a_{\mathfrak{m}} \mathfrak{m}$, we set $\widetilde{\dim}_{\mathbb{C}}(M) := \sum_{\mathfrak{m}} a_{\mathfrak{m}}$. This can be viewed as a  truncated dimension of $M$, in the sense that it gives the sum of the dimensions of the weight subspaces of $M$ that are not killed by the truncation. 
    
     \begin{thm} \label{thmapplicationtoKKKO-SW}
     Assume $\mathfrak{g}$ is of type $A_n , n \geq 1$ and consider the monotonic orientation $Q_0$ of the Dynkin diagram of $\mathfrak{g}$ as in Figure~\ref{figQ0}. Let $M$ be a simple object in $\C_{Q_0}$ and let $\mathfrak{m} := \prod_{i \in I , 1 \leq r \leq n-i+1} Y_{i,\xi(i)-2(r-1)}^{m_{i,r}}$ denote the corresponding dominant monomial. Then one has 
     $$ \frac{\dim_{\mathbb{C}} \left( \mathcal{F}_{Q_0}(M) \right)}{\widetilde{\dim}_{\mathbb{C}}(M)} = \left( \sum_{ \substack{ 1 \leq i \leq n \\ 1 \leq r \leq n-i+1}} i \cdot m_{i,r} \right) ! \prod_{ \substack{1 \leq i \leq n \\ 1 \leq r \leq n-i+1}} \left( \frac{(r-1)!}{(r+i-1)!} \right)^{m_{i,r}}  . $$
      \end{thm}
      
       \begin{proof}
        By Remark~\ref{DtildeAinverse}, we have  
 $$ \td_{Q_0}(A_{i,\xi(i)-2r+1}^{-1}) = \frac{\beta_{\varphi(i,\xi(i)-2r)}}{\beta_{\varphi(i,\xi(i)-2r+2)}} =  \frac{\tau^{r}(\gamma_i)}{\tau^{r-1}(\gamma_i)} = \frac{\alpha_{r+1,r+i}}{\alpha_{r,r+i-1}} $$
      for each $i \in I$ and $1 \leq r < n_{Q_0}(i)$ (see the proof of Proposition~\ref{lem:tilde D of Yit in type A}). Moreover, the positive roots $\alpha_{r,r+i-1}$ and $\alpha_{r+1,r+i}$ are segments of same length $i$, for every $r \geq 1$. Hence we obtain
       $$ \td_{Q_0}(A_{i,\xi(i)-2r+1}^{-1}) \mid_{\alpha_1 = \cdots = \alpha_n = 1} = 1 $$
       for every $i \in I$ and $1 \leq r < n_{Q_0}(i)$. Therefore we have 
   \begin{align*}
    \td_{Q_0} \left( \tchi_q(M) \right)  \mid_{\alpha_1 = \cdots = \alpha_n = 1} \enspace & = \widetilde{\dim}_{\mathbb{C}}(M) \cdot  \td_{Q_0}(\mathfrak{m}) \mid_{\alpha_1 = \cdots = \alpha_n = 1} \\
    & = \widetilde{\dim}_{\mathbb{C}}(M)  \cdot  \prod_{ \substack{ 1 \leq i \leq n \\ 1 \leq r \leq n-i+1}} \left( \td_{Q_0}(Y_{i, \xi(i)-2(r-1)}) \mid_{\alpha_1 = \cdots = \alpha_n = 1}  \right)^{m_{i,r}} .
    \end{align*}
       Now it follows from Proposition~\ref{lem:tilde D of Yit in type A} that for every $i \in I$ and $1 \leq r \leq n-i+1$ one has
       $$ \td_{Q_0}(Y_{i, \xi(i)-2(r-1)}) = \frac{\prod_{1 \leq p \leq r-1 \leq q \leq r+i-2} \alpha_{p,q}}{\prod_{1 \leq p \leq r \leq q \leq r+i-1} \alpha_{p,q}} = \prod_{r \leq q \leq r+i-1} \frac{1}{\alpha_{r,q}} \prod_{1 \leq p \leq r-1} \frac{\alpha_{p,r-1}}{\alpha_{p,r+i-1}} . $$
       Specializing $\alpha_1, \ldots , \alpha_n$ to $1$, this yields
     $$   \td_{Q_0}(Y_{i, \xi(i)-2(r-1)}) \mid_{\alpha_1 = \cdots = \alpha_n = 1}  \enspace = \frac{(r-1)!}{(r+i-1)!} . $$
     Hence we have 
    $$  \td_{Q_0} \left( \tchi_q(M) \right) \mid_{\alpha_1 = \cdots = \alpha_n = 1} \enspace = \prod_{ \substack{ 1 \leq i \leq n \\ 1 \leq r \leq n-i+1}}  \left( \frac{(r-1)!}{(r+i-1)!} \right)^{m_{i,r}}  \cdot \enspace  \widetilde{\dim}_{\mathbb{C}}(M) .$$
    On the other hand, specializing the equality~\eqref{eqnDbarKLR}, we get 
    $$  \barD([\mathcal{F}_{Q_0}(M)]) \mid_{\alpha_1 = \cdots = \alpha_n = 1} = \frac{1}{d!} \dim_{\mathbb{C}}(\mathcal{F}_{Q_0}(M)) $$
    where $d$ is the length of the (unique) element $\beta \in \Gamma_{+}$ such that $\mathcal{F}_{Q_0}(M) \in R(\beta)-mod$. It follows from Kang-Kashiwara-Kim-Oh's construction \cite{KKKOSelecta} that $d = \sum_{i,r} m_{i,r} \mid \tau^{r-1}(\gamma_i) \mid =  \sum_{i,r} i \cdot m_{i,r}$. Therefore Equation~\eqref{eqnDbarDtilde} yields
   $$ \dim_{\mathbb{C}}(\mathcal{F}_{Q_0}(M)) =  \left( \sum_{ \substack{1 \leq i \leq n \\ 1 \leq r \leq n-i+1}} i \cdot m_{i,r} \right) !  \prod_{ \substack{ 1 \leq i \leq n \\ 1 \leq r \leq n-i+1}}  \left( \frac{(r-1)!}{(r+i-1)!} \right)^{m_{i,r}}  \cdot \enspace \widetilde{\dim}_{\mathbb{C}}(M) . $$
       \end{proof}

 \section{Perspectives towards a Mirkovi\'c-Vilonen basis for new cluster algebras}
  \label{finalsection}
  
   In this section we open perspectives relating the morphism $\td_{\xi}$ to the geometric motivations underlying Baumann-Kamnitzer-Knutson's constructions \cite{BKK}. For this purpose, we introduce a cluster algebra $\bAQ$ as a subquotient of $\Axi$ naturally containing $\CN$ and prove that $\td_{\xi}$ descends to a morphism $\barD_Q : \bAQ \rightarrow \mathbb{C}(\alpha_1, \ldots, \alpha_n)$ extending $\barD$. We suggest the existence of a basis in $\bAQ$ containing the Mirkovi\'c-Vilonen basis of $\CN$ where the values of $\barD_Q$ may be interpreted as the equivariant multiplicities of certain closed algebraic varieties in the spirit of Theorem~\ref{thmBKK}. We also point out possible developments via monoidal categorifications of cluster algebras relying on Kashiwara-Kim-Oh-Park's recent advances \cite{KKOPlocalization, KKOP21}. 
   
  \subsection{The cluster algebra $\bAQ$}
   \label{sectionbarDQ}
 
 In this paragraph, we define the cluster algebra $\bAQ$ and show that $\td_{\xi}$ yields a well-defined morphism $\barD_Q : \bAQ \rightarrow \mathbb{C}(\alpha_1, \ldots, \alpha_n)$ extending $\barD$. For this purpose, we prove a technical property of $\td_{\xi}$ (Proposition~\ref{finalprop}) implying that the values of $\td_{\xi}$ on the initial cluster variables $x_t^{\wiQ} , t \geq 1$ of $\Axi$ satisfy certain periodicity properties (Corollary~\ref{corperiodicity}). This mainly relies on the periodicity of the coefficients $\tilde{C}_{i,j}(m)$ established by Hernandez-Leclerc ({{\cite[Corollary 2.3]{HL15}}}). We then prove that the quotient map $\barD_Q$ is well-defined (Corollary~\ref{corfrozen}). 

 Kashiwara-Kim-Oh-Park \cite{KKOP21} recently introduced for each $1 \leq a \leq b \leq + \infty$ a monoidal subcategory $\C^{[a,b]}$ of $\Cxi$ defined as the smallest subcategory of $\Cxi$ containing all the fundamental representations $L(Y_{i,p})$ for $(i,p) \in \varphi^{-1}([a,b])$ and stable under extensions, subquotients and monoidal products. Obviously $\C^{[a,b]}$ can be naturally viewed as a monoidal subcategory of $\C^{[a',b']}$ if $[a,b] \subset [a',b']$.
 
 Here we will be focusing on the category $\C^{[1,2N]}$(recall that $N$ denotes the number of positive roots of $\mathfrak{g}$). It follows from the results in \cite{KKOP21} that the Grothendieck ring $K_0(\C^{[1,2N]})$ has a cluster algebra structure whose  frozen variables are identified with the classes of the Kirillov-Reshetikhin modules $X_{i,p}$ such that $ \left( \varphi(i,p) \right)_{+} > 2N$. These are in bijection with $I$ via $I \ni i \mapsto X_{i, p_i}$ with $p_i := \xi(i)-2h+2$ for each $i \in I$ (where $h$ is the dual Coxeter number of $\mathfrak{g}$, see Section~\ref{remindAR}). 
  We define the cluster algebra $\bAQ$ in the following way:
  $$ \bAQ := K_0(\C^{[1,2N]}) / \left( [X_{i,p_i}] -1 , i \in I \right) . $$
The algebra $\bAQ$ has a cluster algebra structure of rank $2N-n$ with no frozen variables. The set of isomorphism classes of Kirillov-Reshetikhin modules $X_{i,p} , (i,p) \in \varphi^{-1}([1,2N-n])$ forms a cluster in $\bAQ$. The coordinate ring $\CN \simeq \mathbb{C} \otimes K_0(\CQ)$ is naturally embedded into $\bAQ$ via $x_t^{\iQ} \longmapsto [X_{\varphi^{-1}(t)}] , t \in \{1, \ldots , N\}$ as illustrated in Figure~\ref{fig2} below.

  \begin{prop}  \label{finalprop}
  Let $(i,p) \in \Ixi$ such that $p \leq \xi(i)-2h+2$. Then one has 
  $$ \td_{\xi}(Y_{i,p}Y_{i,p+2} \cdots Y_{i,p+2h-2}) = 1 . $$
   \end{prop}
 
 \begin{proof}
 Recall the notation $\mathcal{N}(i,p;j,s)$ from Section~\ref{remindCtilde}. Applying the definition of $\td_{\xi}$ (see~\eqref{defDtildexi}) we have 
 \begin{align*}
  \td_{\xi}(Y_{i,p}Y_{i,p+2} \cdots Y_{i,p+2h-2}) &=  \td_{\xi}(Y_{i,p}) \td_{\xi}(Y_{i,p+2}) \cdots \td_{\xi}( Y_{i,p+2h-2}) \\
  &= \prod_{(j,s) \in \Ixi} \beta_{\varphi(j,s)}^{ \mathcal{N}(i,p;j,s) + \mathcal{N}(i,p+2;j,s) + \cdots + \mathcal{N}(i,p+2h-2;j,s)} \\
  &=  \prod_{(j,s) \in \Ixi} \beta_{\varphi(j,s)}^{\tilde{C}_{i,j}(s-p+1) - \tilde{C}_{i,j}(s-p-2h+1)} .
 \end{align*}
 If $(j,s) \in \Ixi$ is such that $s \geq p+2h$ then {{\cite[Corollary 2.3]{HL15}}} implies $\tilde{C}_{i,j}(s-p-2h+1) = \tilde{C}_{i,j}(s-p+1)$. Recalling moreover that $\tilde{C}_{i,j}(m)=0$ if $m \leq 0$, we can thus rewrite the above expression as 
$$
 \td_{\xi}(Y_{i,p}Y_{i,p+2} \cdots Y_{i,p+2h-2}) =  \prod_{ \substack{ (j,s) \in \Ixi \\ p \leq s < p+2h}} \beta_{\varphi(j,s)}^{\tilde{C}_{i,j}(s-p+1)}  = \prod_{\beta \in \Phi_{+}} \beta^{m(\beta)}
 $$
 where
 $$ 
  m(\beta) := \sum_{(j,s) \in J_{p,\beta}} \tilde{C}_{i,j}(s-p+1), \qquad J_{p,\beta} := \{ (j,s) \in \Ixi \mid p \leq s < p+2h , \beta_{\varphi(j,s)} = \beta \}
 $$
  for each $\beta \in \Phi_{+}$. 
 It follows from Proposition~\ref{propositionFujitaOh} that $J_{p, \beta}$ is non empty.  Moreover, if $(j,s) \in J_{p,\beta}$ then $s+2h \geq p+2h$ and $s-2h<p$. Hence $(j, s \pm 2h) \notin J_{p,\beta}$ and similarly for all the $(j, s \pm 2mh)$ for any $m \in \mathbb{Z} \setminus \{0\}$. Therefore Proposition~\ref{propositionFujitaOh} implies $\sharp J_{p,\beta} \leq 2$ and in case of equality we have $J_{p,\beta} = \{(j,s) ; (j^{*},s+h) \}$ for some $(j,s) \in \Ixi$. We now fix $\beta \in \Phi_{+}$ and prove that $m(\beta)=0$. We distinguish two cases. 
 
  \textbf{Case 1: $\sharp J_{p,\beta}=2$.} Applying Theorem~\ref{thmHLCtilde} we get 
 \begin{align*}
 m(\beta) &= \tilde{C}_{i,j}(s-p+1) + \tilde{C}_{i,j^{*}}(s+h-p+1) = \epsilon_{i,p} \epsilon_{j,s} \langle  \beta_{\varphi(i,p)} , \beta \rangle_Q + \epsilon_{i,p} \epsilon_{j^{*},s+h} \langle  \beta_{\varphi(i,p)} , \beta \rangle_Q.
 \end{align*}
 As $\epsilon_{j,s} = - \epsilon_{j^{*},s+h}$ by Proposition~\ref{propositionFujitaOh}, we get $m(\beta)=0$.

  \textbf{Case 2: $\sharp J_{p,\beta}=1$.} Let us write $J_{p,\beta} := \{(j,s)\}$. Then $m(\beta) = \tilde{C}_{i,j}(s-p+1) =  \epsilon_{i,p} \epsilon_{j,s} \langle  \beta_{\varphi(i,p)} , \beta \rangle_Q$ by Theorem~\ref{thmHLCtilde}. 
  
  On the other hand, by Proposition~\ref{propositionFujitaOh} one has $\beta_{\varphi(j^{*},s-h)} = \beta$. As $(j^{*},s-h) \notin J_{p,\beta}$, one must have $s-h<p$. This implies $s+h < p+2h$. As $(j^{*},s+h) \notin J_{p,\beta}$ this is possible only if $(j^{*},s+h) \notin \Ixi$ i.e. $s+h > \xi(j^{*})$. Therefore we have 
  $$ \xi(j^{*}) < s+h < p+2h \leq \xi(i) .$$
  In particular, $p+2h > s+h$ and $(p+2h) - (s+h) < \xi(i) - \xi(j^{*}) \leq d(i,j^{*})$. Thus Lemma~\ref{littleLemmaforCtilde} yields $\tilde{C}_{i,j^{*}}((p+2h)-(s+h)+1) = 0$. As $p+2h > s+h$ we can again apply Theorem~\ref{thmHLCtilde} and we get 
  $$ 0 =  \tilde{C}_{i,j^{*}}((p+2h)-(s+h)+1)  =  \epsilon_{i,p+2h} \epsilon_{j^{*},s+h} \langle  \beta_{\varphi(i,p+2h)} , \beta \rangle_Q  =  \epsilon_{i,p} \epsilon_{j^{*},s+h} \langle  \beta_{\varphi(i,p)} , \beta \rangle_Q $$
  by Proposition~\ref{propositionFujitaOh}. 
 Thus $\langle  \beta_{\varphi(i,p)} , \beta \rangle_Q = 0$ and hence $m(\beta) = 0$ as well. 
 
  This concludes the proof of the Proposition. 
 \end{proof}
 
 \begin{corollary} \label{corperiodicity}
  For any $t \geq 1$, one has $\td_{\xi} \left( \iota(x_{t+2N}) \right) = \td_{\xi} \left( \iota(x_t) \right)$. 
   \end{corollary} 
 
  \begin{proof}
  Let $(i,p) := \varphi^{-1}(t+2N)$. Then we have $\varphi^{-1}(t)=(i,p+2h)$ (see Section~\ref{remindAR}). We can write 
   \begin{align*}
   \td_{\xi} \left( \iota(x_{t+2N}) \right)  &= \td_{\xi} \left( \tchi_q (X_{i,p}) \right) = \td_{\xi}(Y_{i,p} Y_{i,p+2} \cdots Y_{i,\xi(i)})  \\
   &= \td_{\xi}(Y_{i,p} Y_{i,p+2} \cdots Y_{i,p+2h-2}) \cdot  \td_{\xi}(Y_{i,p+2h} Y_{i,p+2} \cdots Y_{i,\xi(i)})  \\
   &=  \td_{\xi}(Y_{i,p+2h} Y_{i,p+2} \cdots Y_{i,\xi(i)})  \quad  \text{by Proposition~\ref{finalprop}} \\
   &= \td_{\xi} \left( \tchi_q(X_{i,p+2h}) \right) =  \td_{\xi} \left( \iota(x_t) \right). 
      \end{align*}
  \end{proof}
 
  In particular, this implies that the statement of Proposition~\ref{AforDtilde} actually holds for all $t \geq 1$.

  \begin{corollary} \label{corfrozen}
   The morphism $\td_{\xi}$ factors into an algebra morphism 
\begin{align*}
\barD_Q : \bAQ \longrightarrow \mathbb{C}(\alpha_1, \ldots , \alpha_n).
\end{align*}    
   \end{corollary}

  \begin{proof}
  Let $i \in I$ and $p_i := \xi(i)-2h+2$. Applying Proposition~\ref{finalprop} with $p=p_i$, we obtain 
  $$  \td_{\xi} \left( \tchi_q(X_{i,p_i}) \right)  = \td_{\xi}(Y_{i,p_i} Y_{i,p_i+2} \cdots Y_{i,\xi(i)}) =1 . $$
By construction of $\bAQ$, this shows that $\td_{\xi}$ yields a morphism $\barD_Q : \bAQ \longrightarrow \mathbb{C}(\alpha_1, \ldots , \alpha_n)$. 
  \end{proof}

    \begin{remark} \label{rkquotient}
 Corollary~\ref{corperiodicity} shows that most of the information of the morphism $\td_{\xi}$ on $\Axi$ is actually contained in its restriction to $K_0(\C^{[1,2N]})$. The motivation for considering the quotient $\bAQ$ comes from the geometric perspective explained in Section~\ref{eqmforDtilde} below: the trivial values of $\td_{\xi}$ have to be discarded if one wants to interpret the images of $\td_{\xi}$ as  equivariant multiplicities of certain closed algebraic varieties as in Theorem~\ref{thmBKK}. 
 
  The cluster algebra $\bAQ$ contains a seed whose exchange quiver can be viewed as a finite part of Hernandez-Leclerc's quiver $Q^{\wiQ}$, (strictly) containing the exchange quiver of the standard seed $\s^{\iQ}$ of $\CN$. However, unlike $\CN$ or other cluster algebras of the form $K_0(\C^{[1,M]}) , M \geq 1$, the cluster algebra $\bAQ$ does not have any frozen variable. 
   \end{remark}

  \subsection{An example in type $A_3$}
  
   In this section, we study a detailed example of the main features of the present paper when $\mathfrak{g}$ is of type $A_3$ and $\xi$ is a height function adapted to a sink-source orientation of the corresponding Dynkin diagram. 
   
   \smallskip
   
    We choose the height function $\xi : I \longrightarrow \mathbb{Z}$ given by $\xi(1)=\xi(3)=-1$ and $\xi(2)=0$. The corresponding orientation $Q$ of the type $A_3$ Dynkin graph is given by the following sink-source orientation 
    $$ \xymatrix{ Q : & 1 & 2 \ar[l] \ar[r] & 3 } $$
The corresponding Coxeter transformation is given by $\tau_Q = s_2s_1s_3$. We choose the reduced expression $\iQ = (2,1,3,2,1,3)$ of $w_0$, which is clearly adapted to $Q$. The infinite sequence $\wiQ$ is given by 
$$ \wiQ = (2,1,3,2,1,3,2,1,3,2,1,3, \ldots ) . $$
We have $\Ixi = \{ (1, -(2k+1)) , (2, -2k) , (3, -(2k+1)) , k \in \mathbb{Z}_{\geq 0} \}$ and the bijection $\varphi$ is given by $\varphi(1,-(2k+1)) = 2+3k , \varphi(2,-2k) = 1+3k, \varphi(3,-(2k+1)) = 3+3k$.
The exchange quiver $Q^{\wiQ}$ of the initial seed $\s^{\wiQ}$ considered by Hernandez-Leclerc is given by the graph denoted $G^{-}$ in {{\cite[Figure 1]{HLJEMS}}}. In Figure~\ref{fig} we reproduce this quiver, were we put at the node $(i,p)$ the inverse of the value of $\td_{\xi}$ on the cluster variable $x_{\varphi(i,p)}=[X_{i,p}]$.
 Figure~\ref{fig2} provides the exchange quivers and cluster variables (in terms of classes of Kirillov-Reshetikhin modules) for the respective initial seeds of the cluster algebras $\bAQ$ and $\CN$. The picture for $\CN$ is contained in the one of $\bAQ$ in an obvious way. 
  
 \begin{figure} 

 \begin{tikzpicture}
 \begin{scope}[scale=1.3]

    \node (a1) at (-2.2,-1) {$\alpha_2(\alpha_1 + \alpha_2)$};
  \node (a3) at (-2.2,-3) {$\alpha_3(\alpha_2 + \alpha_3)(\alpha_1 + \alpha_2 + \alpha_3)$};
  \node (a5) at (-2.2,-5) {$\alpha_1$};
    \node (a7) at (-2.2,-7) {$1$};
     \node (a9) at (-2.2,-9) {$\alpha_2(\alpha_1 + \alpha_2)$};
      \node (a10) at (-2.2,-10) {$\vdots$};
  \node (b0) at (0,0) {$\alpha_2$};
   \node (b2) at (0,-2) {{\footnotesize $\alpha_2(\alpha_1 + \alpha_2)(\alpha_2 + \alpha_3)(\alpha_1 + \alpha_2 + \alpha_3)$}};
\node (b4) at (0,-4) {$\alpha_1 \alpha_3(\alpha_1 + \alpha_2 + \alpha_3)$};
  \node (b6) at (0,-6) {$1$};
    \node (b8) at (0,-8) {$\alpha_2$};
      \node (b9) at (0,-9) {$\vdots$};
   \node (c1) at (2.2,-1) {$\alpha_2(\alpha_2 + \alpha_3)$};
  \node (c3) at (2.2,-3) {$\alpha_1(\alpha_1 + \alpha_2)(\alpha_1 + \alpha_2 + \alpha_3)$};
  \node (c5) at (2.2,-5) {$\alpha_3$};
    \node (c7) at (2.2,-7) {$1$};
    \node (c9) at (2.2,-9) {$\alpha_2(\alpha_2 + \alpha_3)$};
      \node (c10) at (2.2,-10) {$\vdots$};
     
\draw[->] (a3)--(a1); 
\draw[->] (a5)--(a3);
\draw[->] (a7)--(a5);
\draw[->] (a9)--(a7);
\draw[->] (b2)--(b0); 
\draw[->] (b4)--(b2);
\draw[->] (b6)--(b4);
\draw[->] (b8)--(b6);
\draw[->] (c3)--(c1); 
\draw[->] (c5)--(c3);
\draw[->] (c7)--(c5);
\draw[->] (c9)--(c7);
\draw[->] (b0)--(a1); 
\draw[->] (b0)--(c1);
\draw[->] (b2)--(a3); 
\draw[->] (b2)--(c3);
\draw[->] (b4)--(a5);
\draw[->] (b4)--(c5);
\draw[->] (b6)--(a7); 
\draw[->] (b6)--(c7); 
\draw[->] (b8)--(a9); 
\draw[->] (b8)--(c9); 
\draw[->] (a1)--(b2); 
\draw[->] (c1)--(b2);
\draw[->] (a3)--(b4); 
\draw[->] (c3)--(b4);
\draw[->] (a5)--(b6); 
\draw[->] (c5)--(b6);
\draw[->] (a7)--(b8); 
\draw[->] (c7)--(b8);

\end{scope}
 
 \end{tikzpicture}

 \caption{Values of $\td_{\xi}^{-1}$ on the cluster variables of the seed $\s^{\wiQ}$ of $\Axi$ in type $A_3$.}
 
 \label{fig}

 \end{figure}

 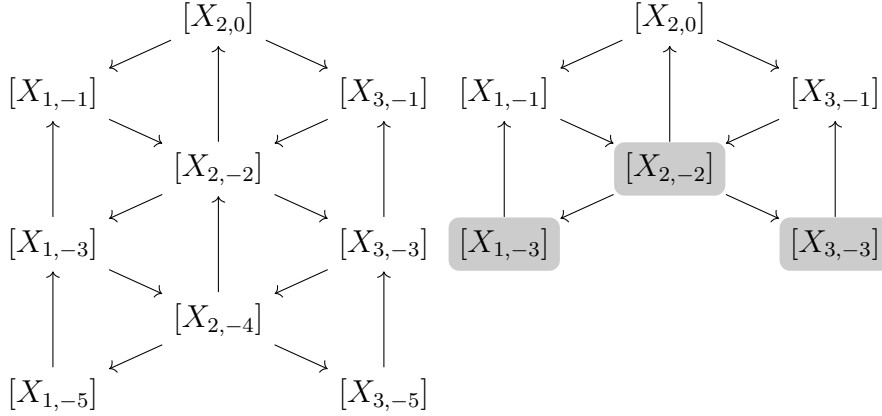
\begin{figure} 

 \begin{tikzpicture}
 \begin{scope}[scale=1]
 
  \tikzstyle{frozen}=[rectangle,rounded corners,fill=black!20]

  \node (a1) at (-4.2,-1) {$[X_{1,-1}]$};
  \node (a3) at (-4.2,-3) {$[X_{1,-3}]$};
  \node (a5) at (-4.2,-5) {$[X_{1,-5}]$};
  \node (b0) at (-2,0) {$[X_{2,0}]$};
  \node (b2) at (-2,-2) {$[X_{2,-2}]$};
  \node (b4) at (-2,-4) {$[X_{2,-4}]$};
  \node (c1) at (0.2,-1) {$[X_{3,-1}]$};
  \node (c3) at (0.2,-3) {$[X_{3,-3}]$};
  \node (c5) at (0.2,-5) {$[X_{3,-5}]$};
 
\draw[->] (a3)--(a1); 
\draw[->] (a5)--(a3);
\draw[->] (b2)--(b0); 
\draw[->] (b4)--(b2);
\draw[->] (c3)--(c1); 
\draw[->] (c5)--(c3);
\draw[->] (b0)--(a1); 
\draw[->] (b0)--(c1);
\draw[->] (b2)--(a3); 
\draw[->] (b2)--(c3);
\draw[->] (b4)--(a5);
\draw[->] (b4)--(c5);
\draw[->] (a1)--(b2); 
\draw[->] (c1)--(b2);
\draw[->] (a3)--(b4); 
\draw[->] (c3)--(b4);

    \node (a21) at (1.8,-1) {$[X_{1,-1}]$};
  \node[frozen] (a23) at (1.8,-3) {$[X_{1,-3}]$};
  \node (b20) at (4,0) {$[X_{2,0}]$};
   \node[frozen] (b22) at (4,-2) {$[X_{2,-2}]$};
   \node (c21) at (6.2,-1) {$[X_{3,-1}]$};
  \node[frozen] (c23) at (6.2,-3) {$[X_{3,-3}]$};
 
\draw[->] (a23)--(a21); 
\draw[->] (b22)--(b20); 
\draw[->] (c23)--(c21); 
\draw[->] (b20)--(a21); 
\draw[->] (b20)--(c21);
\draw[->] (b22)--(a23); 
\draw[->] (b22)--(c23);
\draw[->] (a21)--(b22); 
\draw[->] (c21)--(b22);

\end{scope}
 
 \end{tikzpicture}

 \caption{Initial seeds for the cluster structures of $\bAQ$ (left) and $\CN$ (right) in type $A_3$. The variables in grey boxes are frozen.}
 
 \label{fig2}

 \end{figure}

   \subsection{Towards a monoidal categorification of $\bAQ$}
  It is proved in \cite{KKOP21} that $\C^{[1,2N]}$ is in fact a monoidal categorification of a cluster algebra in the sense of \cite{HL}, i.e. the classes of simple objects in $\C^{[1,2N]}$ belong to the set of cluster monomials in $K_0(\C^{[1,2N]})$. 
  In a previous work \cite{KKOPlocalization}, Kashiwara-Kim-Oh-Park introduced the notion of \textit{commuting family of (graded) braiders} in certain categories of modules over quiver Hecke algebras which were known from \cite{KKKO} to provide monoidal categorifications of cluster algebras (namely the unipotent cells of $\CN$). In \cite{KKOPlocalization}, it is shown that the simple objects corresponding to the frozen variables of these cluster structures form a commuting family of braiders. This allows to construct new monoidal categories by specializing these simple objects to the unit object, following former constructions by Kang-Kashiwara-Kim \cite{KKKInvent18}.
   Therefore, it would be interesting to investigate whether the simple modules $X_{i, p_i} , i \in I$ categorifying the frozen variables in $K_0(\C^{[1,2N]})$ are commuting braiders. This would yield a monoidal category $\bCQ := \C^{[1,2N]}[X_{i,p_i} \simeq \mathbf{1} , \allowbreak i \in I]$ such that $\bAQ = K_0(\bCQ)$.

\subsection{Towards a Mirkovi\'c-Vilonen basis for $\bAQ$}
 \label{eqmforDtilde}

  The morphism $\barD_Q$ defined in Section~\ref{sectionbarDQ} obviously coincides with $\td_Q$ on $\CN$ (viewed as a subalgebra of $\bAQ$). Thus by Theorem~\ref{mainthm1} it also coincides with Baumann-Kamnitzer-Knutson's morphism $\barD$ on $\CN$. We now provide evidences that the morphism $\barD_Q$ on $\bAQ$ can take values not belonging to the image of $\barD$. These values nonetheless share a similar form as the values of $\barD$ on certain reasonable elements of $\CN$ such as cluster variables for instance. 
   
    Let us provide a couple of examples of such new rational fractions. The cluster structure of $\bAQ$ allows us to mutate in the direction of $ x_4 := x_4^{\wiQ} = [X_{2,-2}]$. As recalled in Section~\ref{remindCN}, this mutation produces a new seed consisting in a new quiver $Q'$ (which we do not display here) and with the same cluster variables, except $x_4$ which is replaced by $x'_4$ given by the exchange relation~\eqref{exchangerelation}. As $\barD_Q$ is an algebra morphism, it is then straightforward to compute $\barD_Q(x'_4)$. We find 
  $$ \barD_Q(x'_4) = \frac{\alpha_1 + 2\alpha_2 + \alpha_3}{\alpha_1 \alpha_2 \alpha_3 (\alpha_1+\alpha_2+\alpha_3)} . $$
We can perform similar computations starting with the same initial seed as above and mutating in the direction of $x_5 := x_5^{\wiQ} = [X_{1,-3}]$ or $x_6 := x_6^{\wiQ} = [X_{3,-3}]$. We respectively obtain 
$$ \barD_Q(x'_5) = \frac{\alpha_2 + 2\alpha_3}{\alpha_1 \alpha_2 (\alpha_1 + \alpha_2)}\enspace \text{and} \enspace \barD_Q(x'_6) = \frac{2\alpha_1+ \alpha_2}{\alpha_2 \alpha_3 (\alpha_2 + \alpha_3)} .  $$ 
   It is not hard to check that the rational fractions $ \barD_Q(x'_4), \barD_Q(x'_5), \barD_Q(x'_6)$ do not belong to the image of $\barD$. Nonetheless, these fractions share a similar form as the values taken by $\barD$ on the cluster variables of $\CN$, which belong to the MV basis. Recalling Theorem~\ref{thmBKK} it is therefore natural to ask the following:
   
    \begin{Question} \label{finalquestion}
 Is it possible to construct a basis $\mathcal{B} = (b_Y)$ of $\bAQ$ indexed by a family of closed varieties $Y$, such that 
  \begin{itemize}
  \item The cluster variables $x_t^{\wiQ} , 1 \leq t \leq 2N-n$ belong to $\mathcal{B}$.
  \item The elements of the MV basis of $\CN$ are sent onto elements of $\mathcal{B}$ under the natural injection $\CN \longrightarrow \bAQ$.
  \item For every $Y$, there exists $p \in Y$ such that $\barD_Q(b_Y)$ is equal to the equivariant multiplicity $\epsilon_p^{T}(Y)$ of $Y$ at $p$ with respect to the action of some torus $T$.
  \end{itemize}
     \end{Question}
     
      \begin{remark}
   It is a general fact (see for instance {{\cite[Theorem 4.2]{Brion}}}) that if $X$ is a closed projective scheme with an action of a torus $T$ and if $p$ is a non-degenerate point in a $T$-invariant closed subvariety $Y \subset X$ such that $Y$ is smooth at $p$, then one has $\epsilon_p^{T}(Y) = 1/P$ where $P$ is the product of the weights of the action of $T$ on the tangent space $T_pY$. Therefore Proposition~\ref{AforDtilde}  suggests to investigate possible smoothness properties of the varieties that would correspond to the cluster variables $x_t^{\wiQ} , 1 \leq t \leq 2N-n$ via the first part of Question~\ref{finalquestion}. 
   \end{remark}




\begin{thebibliography}{10}

\bibitem{BKK}
P.~Baumann, J.~Kamnitzer, A.~Knutson, A. with an appendix by A.~Dranowski,
  J.~Kamnitzer, and C.~Morton-Ferguson.
\newblock {The Mirkovic-Vilonen basis and Duistermaat-Heckman measures}.
\newblock {\em {Acta. Math.}}, {227}({no. 1}):{1--101}, {2021}.

\bibitem{BFZ}
A.~Berenstein, S.~Fomin, and A.~Zelevinsky.
\newblock {Cluster algebras III. {U}pper bounds and double {B}ruhat cells}.
\newblock {\em {Duke Math. J.}}, {126}({1}):{1--52}, {2005}.

\bibitem{Lea21}
L.~Bittmann.
\newblock {A quantum cluster algebra approach to representations of simply-laced quantum affine algebras}.
\newblock {\em {Math. Z.}}, {298}({no. 3-4}):{1449--1485}, {2021}.

\bibitem{Brion}
M.~Brion.
\newblock {Equivariant {C}how groups for torus actions}.
\newblock {\em {Transform. Groups}}, {2}({3}):{225--267}, {1997}.

\bibitem{BKM}
J.~Brundan, A.~Kleshchev, and P.J. McNamara.
\newblock {Homological properties of finite-type {K}hovanov-{L}auda-{R}ouquier
  algebras}.
\newblock {\em {Duke Math. J.}}, {163}({no. 7}):{1353--1404}, {2014}.

\bibitem{Casbi}
E.~Casbi.
\newblock {Dominance order and monoidal categorification of cluster algebras}.
\newblock {\em {Pacific J. Math.}}, {305}({2}):{473--537}, {2020}.

\bibitem{Casbi3}
E.~Casbi.
\newblock {Equivariant multiplicities of simply-laced type flag minors}.
\newblock {\em {Represent. Theory}}, {25}:{1049--1092}, {2021}.

\bibitem{CP95a}
V.~Chari and A.~Pressley.
\newblock {Quantum affine algebras and their representations}.
\newblock In {\em {Representations of groups ({B}anff, {AB}, 1994)}},
  volume~{16} of {\em {CMS Conf. Proc.}} {Amer. Math. Soc., Providence, RI},
  {1995}.
  
\bibitem{CDFL} W.~Chang, B.~Duan, C.~Fraser, J.-R.~Li. 
\newblock {Quantum affine algebras and Grassmannians}. 
\newblock {\em Math. Z.}, {296}({no. 3--4}):{1539--1583}, {2020}.

\bibitem{Drin87}
V.~Drinfeld.
\newblock {A new realization of Yangians and of quantum affine algebras}.
\newblock {\em {Dokl. Akad. Nauk SSSR}}, {296}({no.1}):{13--17}, {1987}.


\bibitem{FZJAMS}
S.~Fomin and A.~Zelevinsky.
\newblock {Double Bruhat cells and total positivity}.
\newblock {\em {J. Amer. Math. Soc.}}, {12}({2}):{335--380}, {1999}.

\bibitem{FZ1}
S.~Fomin and A.~Zelevinsky.
\newblock {Cluster algebras I, Foundations}.
\newblock {\em {J. Amer. Math. Soc.}}, 15:497--529, 2002.

\bibitem{FZ4}
S.~Fomin and A.~Zelevinsky.
\newblock {Cluster algebras IV: Coefficients}.
\newblock {\em {Compos. Math.}}, 143:112--164, 2007.

\bibitem{FR1}
E.~Frenkel and N.~Reshetikhin.
\newblock The $q$-characters of representations of quantum affine algebras and
  deformations of $W$-algebras.
\newblock {\em {Contemp. Math.}}, 248:163--205, 1999.

\bibitem{Fujita}
R.~Fujita.
\newblock {Affine highest weight categories and quantum affine Schur-Weyl
  duality of Dynkin quiver types}.
\newblock {\em {Represent. Theory}}, {26}:{211-263}, {2022}.

\bibitem{FHOO}
R.~Fujita, D.~Hernandez, S.-J. Oh and H.~Oya.
\newblock {Isomorphims among quantum Grothendieck rings and propagation of
  positivity}.
\newblock {\em {J. reine angew. Math.}}, {785}:{117--185}, {2022}.

\bibitem{FOComm.Math.Phys.}
R.~Fujita and S.-J. Oh.
\newblock {$Q$-data and representation theory of untwisted quantum affine algebras   }.
\newblock {arXiv:2007.03159v3}, 
\newblock {\em {Comm.Math.Phys.}}, {384}({no. 2}):{1351--1407}, {2021}. 

\bibitem{GLSAENS}
C.~Geiss, B.~Leclerc, and J.~Schr{\"o}er.
\newblock {Semicanonical bases and preprojective algebras}.
\newblock {\em Ann. Sci. Ec. Norm. Sup.}, {38}:{193--253}, {2005}.


\bibitem{GLS}
C.~Geiss, B.~Leclerc, and J.~Schr{\"o}er.
\newblock {Cluster structures on quantum coordinate rings}.
\newblock {\em Selecta Math.}, 19(2):337--397, 2013.

\bibitem{HeAdv}
D.~Hernandez.
\newblock {Algebraic approach to $q,t$-characters}.
\newblock {\em {Adv. Math.}}, {187}:{1--52}, {2004}.

\bibitem{H06}
D.~Hernandez.
\newblock {The {K}irillov-{R}eshetikhin conjecture and solutions of
  {$T$}-systems}.
\newblock {\em {J. reine angew. Math.}}, {596}:{63--87}, {2006}.

\bibitem{HL}
D.~Hernandez and B.~Leclerc.
\newblock {Cluster algebras and quantum affine algebras}.
\newblock {\em Duke Math. J.}, 154(2):265--341, 2010.

\bibitem{HL15}
D.~Hernandez and B.~Leclerc.
\newblock {Quantum Grothendieck rings and derived Hall algebras}.
\newblock {\em {J. reine angew. Math.}}, {701}:{77--126}, {2015}.

\bibitem{HLJEMS}
D.~Hernandez and B.~Leclerc.
\newblock {A cluster algebra apporach to $q$-characters of Kirillov-Reshetikhin
  modules}.
\newblock {\em {J. Eur. Math. Soc.}}, {18}({5}):{1113--1159}, {2016}.

\bibitem{Jim85}
M.~Jimbo.
\newblock {A $q$-difference analogue of {$U({\mathfrak{g}})$} and the
  {Y}ang-{B}axter equation}.
\newblock {\em {Lett. Math. Phys.}}, {10}({1}):{63--69}, {1985}.

\bibitem{Joseph}
A.~Joseph.
\newblock {On the variety of a highest weight module}.
\newblock {\em {J. Algebra}}, {88}:{238--278}, {1984}.

\bibitem{Kashicrystal}
M.~Kashiwara.
\newblock {On crystal bases of the $q$-analogue of universal enveloping
  algebras}.
\newblock {\em {Duke Math. J.}}, {63}:{465--516}, {1991}.

\bibitem{KKKInvent18} S.-J. Kang, M.~Kashiwara, M.~Kim.
\newblock {Symmetric quiver Hecke algebras and R-matrices of quantum affine algebras}. 
\newblock {\em {Invent. Math.}}, {211}({2}):{591--685}, {2018}. 

\bibitem{KKKOSelecta}
S.-J. Kang, M.~Kashiwara, M.~Kim, and S.-J. Oh.
\newblock {Symmetric quiver {H}ecke algebras and {$R$}-matrices of quantum
  affine algebras {IV}}.
\newblock {\em {Selecta Math. (N.S.)}}, {22}({4}):{1987--2015}, {2016}.


\bibitem{KKKO}
S.-J. Kang, M.~Kashiwara, M.~Kim, and S.-j. Oh.
\newblock {Monoidal categorification of cluster algebras}.
\newblock {\em J. Amer. Math. Soc.}, 31(2):349--426, 2018.


\bibitem{KKOPlocalization} M.~Kashiwara, M.~Kim, S.-j.~Oh, and E.~Park. 
\newblock {Localizations for quiver Hecke algebras}.
\newblock {\em {Pure and Applied Mathematics Quarterly}}, {17}({no. 4}):{1465--1548}, {2021}. 

\bibitem{KKOP21} M.~Kashiwara, M.~Kim, S.-j.~Oh, and E.~Park. 
\newblock {Monoidal categorification and quantum affine algebras II}. 
\newblock {arXiv:2103.10067}.

\bibitem{KL}
M.~Khovanov and A.~D. Lauda.
\newblock {A diagrammatic approach to categorification of quantum groups I}.
\newblock {\em Represent. Theory}, 13:309--347, 2009.

\bibitem{KRhom}
A.~Kleshchev and A.~Ram.
\newblock {Homogeneous representations of {K}hovanov-{L}auda algebras}.
\newblock {\em {J. Eur. Math. Soc.}}, {12}({5}):{1293--1306}, {2010}.

\bibitem{KR}
A.~Kleshchev and A.~Ram.
\newblock {Representations of Khovanov-Lauda-Rouquier algebras and
  combinatorics of Lyndon words}.
\newblock {\em {Math. Ann.}}, 349(4):943--975, 2011.

\bibitem{Leclerc}
B.~Leclerc.
\newblock {Dual canonical bases, quantum shuffles and q-characters}.
\newblock {\em Math. Z.}, 246(4):691--732, 2004.

\bibitem{Lusztig}
G.~Lusztig.
\newblock {Canonical bases arising from quantized enveloping algebras}.
\newblock {\em {J. Amer. Math. Soc.}}, {3}({2}):{447--498}, {1990}.

\bibitem{McNamarafinite}
Peter~J. McNamara.
\newblock {Finite dimensional representations of {K}hovanov-{L}auda-{R}ouquier
  algebras {I}: {F}inite type},.
\newblock {\em {J. reine angew. Math.}}, {707}:{103--124}, {2015}.

\bibitem{MV}
I.~Mirkovi\'{c} and K.~Vilonen.
\newblock {Perverse sheaves on affine {G}rassmannians and {L}anglands duality}.
\newblock {\em {Math. Res. Lett.}}, {7}({1}):{13--24}, {2000}.

\bibitem{Muthiah}
D.~Muthiah.
\newblock {Weyl group action on weight zero {M}irkovic-{V}ilonen basis and
  equivariant multiplicities}.
\newblock {\em {Adv. Math.}}, {385}:{Paper No. 107793, 40 pp}, {2021}.

\bibitem{Nakada}
K.~Nakada.
\newblock {Colored hook formula for a generalized {Y}oung diagram}.
\newblock {\em {Osaka J. Math.}}, {45}({4}):{1085--1120}, {2008}.


\bibitem{NakaAnn.Math}
H.~Nakajima.
\newblock {Quiver varieties and $t$-analogs of $q$-characters of quantum affine
  algebras}.
\newblock {\em {Ann. of Math. (2)}}, {160}({no. 3}):{1057--1097}, {2004}.

\bibitem{P2}
R.~A. Proctor.
\newblock {Minuscule elements of {W}eyl groups, the numbers game, and
  {$d$}-complete posets}.
\newblock {\em {J. Algebra}}, {213}({1}):{272--303}, {1999}.

\bibitem{Rossmann}
W.~Rossmann.
\newblock {Equivariant multiplicities on complex varieties}.
\newblock {\em {Ast{\'e}risque}}, {173-174}:{313--330}, {1989}.

\bibitem{R}
R.~Rouquier.
\newblock {Quiver Hecke algebras and 2-Lie algebras}.
\newblock {\em Algebra. Colloq.}, 19(2):359--410, 2012.

\bibitem{Sage} SageMath, the Sage Mathematics Software System (Version 9.0),
The Sage Developers, 2020, https://www.sagemath.org.

\bibitem{Stem}
J.~R. Stembridge.
\newblock {Minuscule elements of {W}eyl groups}.
\newblock {\em {J. Algebra}}, {235}({2}):{722--743}, {2001}.

\bibitem{VV}
M.~Varagnolo and E.~Vasserot.
\newblock {Canonical bases and {K}{L}{R} algebras}.
\newblock {\em {J. reine angew. Math.}}, {659}:{67--100}, {2011}.

\bibitem{Zelikson}
S.~Zelikson.
\newblock {{A}uslander-{R}eiten quivers and the {C}oxeter complex}.
\newblock {\em {Alg. Rep. Theory}}, {8}:{35--55}, {2005}.

\end{thebibliography}
\end{document}